\documentclass[12pt]{article}

\usepackage{amsfonts}

\def\a{\mathfrak{a}}

\def\i{\mathfrak{i}}

\def\s{\mathfrak{s}}
\def\t{\mathfrak{t}}

\def\si{\sigma}
\def\ep{\varepsilon}

\def\CC{{\cal C}}
\def\O{\cal O}

\def\LL{\Lambda}
\def\l{\lambda}
\def\aa{\alpha}

\def\cu{C^{\infty}(U_0\backslash G, \psi)}

\def\lu{L^2(U_0\backslash G, \psi)}
\def\ccu{C_c^{\infty}(U_0\backslash G, \psi)}
\def\cmu{C^{\infty}(U_0\cap M \backslash M, \psi)}
\def\ccmu{C_c^{\infty}(U_0\cap M \backslash M, \psi)}
\def\lgp{Wh(P, \si)\otimes  i_P^G E}

\def\sgpss{sous-groupe parabolique semi-standard}

\def\DD{\Delta}
\def\Si{\Sigma}

\def\igps{i^G_P\sigma}

\def\N{\mathbb{N}}

\def\Z{\mathbb{Z}}
\def\R{\mathbb{R}}
\def\C{\mathbb{C}}

\def\fd {\hspace{0.35cm} \raise
-0.5mm\hbox{$\blacksquare$}\\}

\def\qed{{\null\hfill\ \raise3pt\hbox{\framebox[0.1in]{}}}\break\null}
\newtheorem{theo}{Th\'eor\`eme}
\newtheorem{prop}{Proposition}
\newtheorem{lem}{Lemme}
\newtheorem{cor}{Corollaire}
\newtheorem{rem}{Remarque}
\newtheorem{defi}{D\' efinition}

\def\ste{\par\smallskip\noindent}

\def\dem{ {\em D\' emonstration: \ste }}

\def\beq{\begin{equation}}
\def\eeq{\end{equation}}

\def\E{\cal E}

\def\O {\mathcal O}

\newenvironment{res}
               {\begin{equation}\begin{minipage}{0.85\textwidth}}
               {\end{minipage}\end{equation}}
\def\ber{\begin{res}}
\def\eer{\end{res}}

\catcode`\Ž=\active\def Ž{\'e}
\catcode`\=\active\def {\`e}
\catcode`\ˆ=\active\def ˆ{\`a}
\catcode`\=\active\def {\`u}
\catcode`\=\active\def {\^{e}}
\catcode`\"=\active\def "{\^{i}}
\catcode`\™=\active\def ™{\^{o}}
\catcode`\ž=\active\def ž{\^{u}}
\catcode`\=\active\def {\c{c}}

\title{Thorme de Paley-Wiener pour les fonctions de Whittaker sur un groupe rŽductif $p$-adique}

\author{ Patrick Delorme}
\date{}

\begin{document}

\maketitle

\section{Introduction}
Soit $F$ un corps local non archimŽdien et soit $G$ le groupe des points sur $F$ d'un groupe rŽductif connexe dŽfini sur $F$. Soit  $(P_0, P_0^-)$ un couple de sous-groupes paraboliques minimaux opposŽs de $G$. On note $M_0$ leur sous-groupe de LŽvi commun.  Soit $U_0$ le radical unipotent de $P_0$ et $A_0$ le plus grand  tore dŽployŽ de $M_0$. Soit $K$ un bon sous-groupe compact maximal  de $G$ relativement ˆ $A_0$.  Soit $\psi$  un caractre unitaire lisse de $U_0$ non dŽgŽnŽrŽ, i.e. tel que pour toute    racine $P_0$-simple de $A_0$, $\aa$,  sa restriction au  sous-groupe radiciel $({U_0})_\alpha$ soit non triviale.  
 \\On note $\cu$ l'espace des fonctions de Whittaker lisses sur $G$, i.e. des fonctions, $f$, sur $G$,  invariantes ˆ droite par un sous-groupe compact ouvert de $G$ et telles que:
  $$f(u_0g)= \psi(u_0) f(g), g\in G, u_0 \in U_0.$$ 
  On note $\ccu$ l'espace des ŽlŽments de $\cu$ qui sont ˆ support compact modulo $U_0$.
  \\ Le but de cet article est de dŽfinir une transformŽe de Fourier pour cet espace et d'en caractŽriser l'image. Ce ThŽorme est l'analogue pour les fonctions de Whittaker du ThŽorme de Paley-Wiener pour les fonctions sur le groupe, du ˆ Joseph Bernstein [B2] et dont Heiermann [H] a fourni une autre preuve. Ici, c'est le travail d'Heiermann qui sert  de fil conducteur ˆ notre  preuve.\\
  Noter que pour les fonctions de Whittaker sur un groupe rŽductif rŽel, la formule de Plancherel a ŽtŽ Žtablie (Harish-Chandra, non publiŽ,  Wallach [Wall], Chapitre 15). La formule de Plancherel dans le cas $p$-adique a Žglement ŽtŽ Žtable par Haris-Chandra (non publiŽ), comme nous l'a indiquŽ Laurent Clozel. \\
 Si $(\pi, V)$ est une reprŽsentation lisse de $G$, on note $Wh(\pi, U_0)$ ou $Wh(\pi)$ l'espace des formes linŽaires, $\xi$,  sur $V$ telles que:
$$ \langle \xi, \pi(u_0)v \rangle  = \psi(u_0)  \langle \xi, v\rangle  , v\in V, u_0\in U_0.$$
Si $\pi$ est de longueur finie,  $Wh(\pi) $ est de dimension finie (cf.[BuHen], ThŽorme 4.2) et [D3] pour une autre dŽmonstration).   Notez que si le groupe n'est pas quasi-dŽployŽ, cette dimension n'est pas toujours infŽrieure ou Žgale ˆ  1, comme expliquŽ dans l'introduction de [BuHen] (voir aussi le ThŽorme \ref{jacquetint}, appliquŽ ˆ un sous-groupe parabolique minimal de $G$). 
\\ Si $v\in V$, on note $c_{\xi, v}$ le coefficient gŽnŽralisŽ dŽfini par:
$$ c_{\xi, v}(g):=  \langle \xi, \pi(g) v\rangle  , \>g \in G.$$
C'est un ŽlŽment de $\cu$. \\
 Dans la suite la phrase ``Soit $P= MU$ un sous-groupe parabolique semi-standard de $G$'' voudra dire que $P$ contient $M_0$, que $M$ est son sous-groupe de LŽvi contenant $M_0$ et que $U$ est son radical unipotent. On note 
$X(M)$ le groupe des caractres non ramifiŽs de $M$. C'est un tore complexe. On note $\delta_P$ la fonction module de $P$. Soit    
$(\si, E)$ une  reprŽsentation lisse de  $M$. On note
$(i^G_P  \si, i^G_P E)$ l'induite parabolique de $(\si, E)$. On suppose $\si$ irrŽductible et on note $\O$
 l'ensemble des classes d'Žquivalences des reprŽsentations 
 $\si_\chi:= \si\otimes \chi$,  $ \chi \in X(M)$, qui est un tore complexe. On utilise dans la suite la notion de fonction sur $\O$, qui dŽpendent des objets concrets $\si_\chi$, avec des rgles de transformations pour tenir compte des Žquivalences de reprŽsentations. On note que 
$Wh(\si _\chi ):= Wh(\si_ \chi, U_0\cap M)$ est indŽpendant de $\chi \in X(M) $, car $\chi$  est trivial sur  $ U_0\cap M$. Par restriction des fonctions ˆ $K$, les reprŽsentations $i^G_P  \si_\chi$ admettent une rŽalisation dans un espace indŽpendant de $\chi$, la rŽalisation compacte.\\ 
\\ {\bf ThŽorme: Fonctionnelles de Jacquet pour les sous-groupes paraboliques anti-standard}. 
{ \em Soit $P=MU$ un sous-groupe parabolique anti-standard de $G$, i.e. contenant $P_0^-$,  $P^-=MU^-$ le sous-groupe parabolique opposŽ relativement ˆ $M$.  On note $(\si, E)$ une reprŽsentation lisse de longueur finie de $M$. \\
Il y a un isomorphisme naturel $Wh(\si) \to Wh (i^G_P  \si)$ notŽ  $\eta \mapsto \xi(P, \si, \eta)$ (Rodier [R], Casselman-Shalika [CS], Shahidi [Sh], Proposition 3.1)tel que:
\ste (i) Pour tout  $v$ dans l'espace de la rŽalisation compacte et $\eta \in Wh(\si) $, $ \langle \xi (P, \si_\chi, \eta), v\rangle  $  est polynomiale en  $\chi \in X(M)$. \ste  (ii) De plus,  si  $\chi \in X(M)$ est suffisamment $P$-dominant, le  vecteur distribution $\xi(P, \si, \eta)$ est donnŽ  par la fonction  sur $G$ ˆ valeurs dans  $E^*$, $\tilde{ \xi} (P, \si, \eta)$, dŽfinie par: 
\ste  $ \langle \tilde { \xi} (P, \si, \eta)(um u^-), e\rangle   =  \psi(u^-)^{-1} \langle \eta,  Ê\si (m^{-1})\delta_P^{1/2} (m)e\rangle  , \> e \in E$, 
\ste $\tilde { \xi} (P, \si, \eta) (g)= 0$ isi $g \notin UMU^-= PU_0$.}\\ \\
On dŽfinit les fonctionnelles de Jacquet  pour un sous-groupe parabolique semi-standard de $G$,  $P=MU$, par 
 transport de structure. \ste Plus prŽcisŽment,  soit $K$ un bon sous-groupe compact maximal  de $G$ relativement ˆ $A_0$. On note $\overline{ W}^G$ (resp. $\overline{ W}^M$) le groupe de Weyl de $G$ (resp. $M$) relativement ˆ $M_0$.  On fixe un ensemble $W^G$ de reprŽsentants dans $K$  de  $\overline{ W}^G$. Pour un bon choix de $w\in W^G$  
 tel que $Q:= wPw^{-1}$ est anti-standard, on dŽfinit 
 $$ \xi(P, \si, \eta): = \xi (Q, w\si, \eta) \circ \lambda(w),  \eta \in  Wh(P, \si):= Wh(w\si),$$
o $\lambda (w) $ est la translation ˆ gauche par $w$, qui entrelace $i^G_P\si$ et $i^G_Qw\si$.
Ce choix de $w$ conduit malheureusement ˆ quelques complications et ˆ quelques  calculs pŽnibles, mais que nous n'avons pas su Žviter.\\
Soit $P$,  $Q$  des sous-groupes paraboliques semi-standard de $G$,   possŽdant  le m\^eme sous-groupe de LŽvi,  $M$, contenant $M_0$.\ste On suppose $\si$ de longueur finie. On introduit les intŽgrales d'entrelacement, $A(Q, P, \si)$, qui, lorsqu'elles sont dŽfinies, entrelacent  $i^G_P \si$ et  $i^G_Q \si$. \\  
Les intŽgrales d'entrelacement transforment les fonctionnelles de Jacquet en des fonctionnelles de Jacquet, ce qui permet d'introduire les matrices $B$: \\
{\em Il existe une unique fonction rationnelle sur $X(M)$ ˆ valeurs dans $End (Wh (Q, \si), Wh (P, \si))$, $\chi \mapsto B(P, Q, \si_\chi)$,  telle que:}
$$\xi(Q, \si_\chi, \eta) \circ A(Q, P, \si_\chi)= \xi (P,\si_\chi,  B(P, Q, \si_\chi) \eta). $$ 
On dŽfinit les intŽgrales de Jacquet par: $$ E^G_P(\si, \eta, v) = c_{\xi, v}  \in  \cu , $$ o  $\xi= \xi(P, \si, \eta)$, $v \in i^G_P \si.$\\
Soit $ (\pi, V)$ {une reprŽsentation lisse unitaire irrŽductible et cuspidale de  $G$. On note $A_G$ le plus grand tore dŽployŽ du centre de $G$. Gr\^ace au Lemme de Schur, on voit  qu'il  existe un unique produit scalaire sur $Wh(\pi)$ tel que:
$$\int _{A_G U_0\backslash G}  c_{\xi, v}(g) \overline{ c_{\xi', v'}(g)    }dg= (\xi, \xi') (v, v'), \>\xi, \xi' \in Wh(\pi), v,v' \in V.$$
On notera que $\pi$ cuspidale implique  que $c_{\xi,v}$ est ˆ support compact modulo $A_GU_0$ donc l'intŽgrale est bien dŽfinie. \ste
On dŽfinit maintenant la transformŽe de Fourier-Whittaker de  $f\in \ccu$, $\hat{f}$,   comme suit. Pour tout sous-groupe parabolique semi-standard de $G$,  $P=MU$, et toute reprŽsentation lisse unitaire irrŽductible cuspidale de $M$,  $(\si, E)$, le ThŽorme de reprŽsentation de Riesz montre qu'il existe un  unique ŽlŽment de $ Wh(\si)\otimes i^G_PE$,   $\hat{f}(P, \si)$,  tel que:  $$(\hat{f}(P, \si), \eta \otimes  v) = \int_{U_0\backslash G } f(g) \overline {E^G_P(\si,  \eta ,  v) (g)}dg, \eta \in Wh(P, \si), v\in i^G_PE.$$
 Notons $F$ au lieu de $\hat{f}$. Alors $F$  vŽrifie les propriŽtŽs suivantes:\ste
 1) a) L 'application $\chi \mapsto F(P, \si_\chi) $,  dŽfinie pour $\chi$ ŽlŽment du groupe des caractres non ramifiŽs  unitaires de $M$, $X(M)_u$,  s'Žtend en une fonction polynomiale   sur $X(M)$. En particulier, dans la rŽalisation compacte, cette application est ˆ valeurs dans un espace vectoriel de dimension finie.\\   
  b) Si $(\si, E)$ et $(\si_1, E_1) $ sont unitairement Žquivalentes, $F(P, \si)$ et $F(P, \si_1)$ vŽrifient une relation de  transport de structure.\\
  c) On peut dŽfinir pour $g\in G$, $\rho_\bullet (g) F(P, \si)$, en posant $(\rho_\bullet (g)F)(P, \si) = (Id \otimes i^G_P \si(g)) F(P, \si)$. 
  Alors  il existe un sous-groupe compact ouvert de $G$, $H$, tel $\rho_{\bullet}(h)F= F$ pour tout $h\in H$. \\
  On appelle orbite inertielle unitaire, $\O_u$,  d'une reprŽsentation lisse unitaire irrŽductible et cuspidale, $(\si, E)$,  l'ensemble des classes d'Žquivalence unitaire des reprŽsentations  $\si\otimes \chi$, $\chi\in X(M)_u$, qui est muni d'une mesure non nulle et $X(M)_u$-invariante, convenablement normalisŽe.  \\ On rŽsume les propriŽtŽs a), b), c) en disant que pour toute orbite inertielle unitaire, $\O_u$,   d'une  reprŽsentation lisse unitaire irrŽductible cuspidale de $M$,  $(\si, E)$, $F(P, )$ dŽfinit un ŽlŽment de 
  $Pol(\O_u, Wh(P, )\otimes   i^G_P)$. Alors $\rho_\bullet $ dŽfinit une reprŽsentation lisse de $G$ sur cet espace.  \ste 2) Il existe de plus un sous-groupe compact ouvert, $H$, de $G$ fixant les diverses applications $F(P,.)$, quand $P$ et $\O$ varient.  \\
  3) Si $P$ est semi-standard et $w$ est choisi comme plus haut, on note $w.P:= wPw^{-1}$ et l'on a:
$$F(w.P, w \si)= (Id \otimes \lambda (w) ) F(P, \si).$$
  4) Si $P$ est un sous-groupe parabolique  anti-standard et  $Q$ est un  sous-groupe parabolique semi-standard  de $G$ admettant $M$ pour sous-groupe de LŽvi: $$(B(Q, P, \si) \otimes Id) F(P, \si) = (Id \otimes A(Q, P, \si) ) F(Q, \si).$$
 Cette relation est loin d'\^etre triviale car elle utilise l'ŽgalitŽ, pour $\si$ reprŽsentation  lisse unitaire irrŽductible et cuspidale:  
  \beq \label{B=B} B(Q, P, \si)^*= B(P,Q, \si),\eeq que l'on Žtablit beaucoup plus loin. \\\\
  {\bf ThŽorme principal}\\
 {\em Si  $F$ satisfait ces propriŽtŽs,  elle est de la forme $\hat{f} $ pour un unique $f \in \ccu$.}\\\\
  Donnons une idŽe de notre  preuve.\\ L'unicitŽ rŽsulte de l'adaptation (cf. \ref{app}) d'un rŽsultat de Joseph Bernstein (cf. [B1]). \\ Pour l'existence, on  commence par introduire les paquets d'ondes. Avec les notations ci-dessus et en supposant $P$ anti-standard, soit   $\phi \in Pol(\O_u, Wh(P,) \otimes  i^G_P)  $.  On dŽfinit le paquet d'ondes $f_\phi$ par: 
$$f_\phi(g):= \int_{\O_u} E(P, \si, \phi(\si))(g)d\si, g \in G,$$
qui est un ŽlŽment de $\cu$. 
On dŽfinit ensuite la notion de $\phi$ rŽgulire (resp.  trs rŽgulire), que nous ne dŽtaillerons pas dans cette introduction. Elle fait intervenir les intŽgrales d'entrelacement et les matrices $B$. Il  y a beaucoup de $\phi$ trs rŽgulires. En effet, si $\phi$ est non nulle, il existe un ŽlŽment, $z$,  du centre de Bernstein, $ZB(G)$, tel que 
$\rho_\bullet (z)  \phi$ soit trs rŽgulire et non nulle. \\
On montre alors: \\
{\bf Proposition}: {\em Si $\phi$ est  rŽgulire, alors  $f_\phi$ est un ŽlŽment de $ \ccu$.}\\
La preuve passe par le terme constant des ŽlŽments de $\cu$ le long d'un sous-groupe parabolique standard de $G$, $P$, qui donne lieu ˆ une application(cf. [D3]): 
$$\cu \to \cmu,   \>\>f \mapsto f_P. $$ 
Sa dŽfinition utilise le ThŽorme de deuxime adjonction de J. Bernstein (cf. [B2], [B3], [Bu]).
Le terme constant des intŽgrales de Jacquet se calcule (cf. ThŽorme \ref{thct}), comme le terme constant faible  des coefficients ordinaires (cf [W], section V),  avec quelques variantes nŽammoins. Plus prŽcisŽment, soit $P=MU$ (resp. $P'= M'U'$) un sous-groupe parabolique anti-standard (resp. standard) de $G$.    Soit  $(\si, E)$ une reprŽsentation  lisse  cuspidale irrŽductible de $M$ et $\phi \in Wh(P, \si)\otimes i^G_P  \si$. \\Soit :
$$\overline{W}(M'\vert G \vert M)= \overline{W}^{M'} \backslash  \{ s \in\overline{W}^G, sMs^{-1} \subset M'\}.$$On rŽalise les reprŽsentations $i^G_P  \si_\chi$ dans un espace indŽpendant de $\chi \in X(M)$. Alors, pour $\phi$ dans cet espace,  on a une identitŽ de fractions rationnelles sur $X(M)$:
$$E(P, \si_\chi, \phi)_{P' }= \sum _{s \in \overline{W}(M'\vert G\vert M) } E_s(\si_\chi, \phi),$$
o  $E_s(\si_\chi, \phi)$ fait intervenir les intŽgrales d'entrelacements et les matrices $B$. La dŽfinition de $ \phi \in  Pol(\O_u, Wh(P,) \otimes  i^G_P)  $ trs rŽgulire comporte notamment le fait que $E_s (\si_\chi, \phi(\si_\chi))(g)$ est polynomiale en $\chi \in X(M)$, pour tout $g \in G$. \\
Par des arguments de translation ˆ droite par des ŽlŽments de $A_0$, la preuve de la Proposition se ramne ˆ  montrer que pour tout $\phi$  rŽgulire, la restriction ˆ la chambre de Weyl nŽgative,  $A_0^-$,  de $A_0$ relativement ˆ $P_0$, est ˆ support compact. \\
On procde ensuite par rŽcurrence sur la dimension de $G$.  On utilise ensuite une partition de $A_0^-$, $(X_{P'})$ indexŽe par les sous-groupes paraboliques standard, $P'$,  de $G$   telle que sur $X_{P'}$, $E(P, \si_\chi, \phi(\si_\chi) )$ est Žgale ˆ $ E(P, \si_\chi, \phi(\si_\chi))_{P'}$, pour tout $\chi \in X(M)$. Puis on utilise la formule pour le terme constant des intŽgrales de Jacquet. L'hypothse que $\phi$ est  rŽgulire permet d'utiliser l'hypothse de rŽcurrence pour finir la preuve de la Proposition. \\
Ensuite on calcule (cf. ThŽorme  \ref{fphihat}) $\hat{f}_\phi$ pour $\phi $ trs rŽgulire en utilisant la formule du terme constant des intŽgrales de Jacquet. 
On procde comme dans [W], section VI, en introduisant la transformŽe unipotente de $f \in \ccu$ relative ˆ $P=MU$, sous-groupe parabolique anti-standard de $G$, $f^P$ dŽfinie par 
$$f^P(m):= \delta_P^{1/2}(m) \int _Uf(mu)du, m\in M.$$
On en dŽduit une formule pour le produit scalaire $L^2$ de deux paquets d'ondes $f_\phi$ et $f_{\phi'}$ (cf. Proposition \ref{prodscalpaq}):$$(f_\phi, f_{\phi'}):= \int_{U_0\backslash G}f_{\phi}(g)\overline{ f_{\phi'}(g)}dg.$$
En Žchangeant le r™le de $\phi$ et $\phi'$, on obtient une autre expression de ce produit scalaire. En faisant varier $\phi$ et $\phi'$ dans l'ŽgalitŽ de ces 2 expressions de $(f_\phi, f_{\phi'})$, on en dŽduit la formule d'adjonction(\ref{B=B}) (cf. ThŽorme \ref{adjB}).\\Pour poursuivre, on utilise le rŽsultat suivant d'Heiermann (cf. [H],  Proposition 0.2). Soit $e_H$ la mesure de Haar normalisŽe d'un sous-groupe compact ouvert $H$ contenu dans $K$.  Soit $P=MU$ un sous-groupe parabolique anti-standard de $G$ et $\O_u$ l'orbite inertielle unitaire d'une  reprŽsentation  lisse unitaire irrŽductible cuspidale de $M$. Alors, il existe 
 $\zeta(P,.) \in Pol(\O_u,Hom(i_P^G, i^G_{P^-}))$  tel que,  pour $\si\in \O$: 
$$(i^G_P\si) (e_H) = \sum _{w \in W^G, w{\cal  O} = \cal O} A(P,  w^{-1}P^-w, \si) \lambda(w^{-1}) \zeta(P, w\si) 
\lambda(w)A(w^{-1}Pw, P, \si).$$
Pour finir la preuve du ThŽorme principal, soit 
 $F$ satisfaisant les conditions de celui-ci et $H$-invariant. On dŽfinit:   
 $$ \Phi_{\O}(\si):= (Id\otimes A(P^-, P, \si)^{-1}) \zeta(P, \si) F(P, \si), \si)\in \O_u$$
 et on introduit le paquet d'onde dŽcalŽ $f^{sh}_{\Phi_{ \O}}$ dŽfini comme le paquet d'ondes ordinaire mais en intŽgrant sur ${\cal O}_u\Lambda$, o $\Lambda \in X(M)$ est suffisamment $P$-antidominant . On  montre de manire analogue ˆ [H], Proposition 2.1:\\
{\bf Proposition} 
{\em  Pour  $P$ anti-standard, le paquet d'ondes dŽcalŽ $ f_{ \O}:= f^{sh}_{\Phi_{ \O}}$ est ŽlŽment de  $\ccu$}.\ste
 Maintenant on dŽfinit $f$ comme Žtant la somme finie sur les classes de conjugaison de couples $(M, \O)$ des $f_\O$ non nulles. Alors le ThŽorme rŽsulte du fait que la transformŽe de Fourier-Whittaker de $f$ est Žgale ˆ $F$. Pour montrer  cela,  on se rŽduit au cas o $\Phi_\O$ est trs rŽgulire,   auquel cas les paquets d'ondes dŽcalŽs sont Žgaux aux paquets d'ondes ordinaires, par holomorphie. Pour effectuer cette rŽduction on remarque que,  si  $z $ est   un ŽlŽment du centre de Bernstein de $G$ et $\rho$ est la reprŽsentation rŽgulire droite de $G$ sur $\ccu$, on a: 
$$(\rho(z)f)\hat {}=  \rho_\bullet (z ) \hat{f}, \rho(z)f^{sh} _\Phi= f^{sh}_{\rho_\bullet (z)\Phi}. $$
On conclut gr\^ace au calcul de la transformation de Fourier des paquets d'ondes mentionnŽ plus haut. \\
L'unicitŽ rŽsulte de l'injectivitŽ de la transformation de Fourier-Whittaker. Celle-ci est une consŽquence de l'adaption ˆ notre contexte de rŽsultats de Bernstein sur les espaces homognes [B1], que nous donnons en appendice.\\
Comme suggŽrŽ par le referee, donnons quelques prŽcisions supplŽmentaires   dans le cas de $G=SL(2, F)$, $A_0$ Žtant le tore diagonal. \\On s'intŽresse aux transformŽes de Fourier-Whittaker de l'espace $\ccu ^K$ des fonctions $K$-invariantes ˆ droite de $\ccu$. \\On note $P$ le sous-groupe parabolique   de $G$ opposŽ au sous-groupe parabolique standard $P_0$. On notera aussi $P^-$ au lieu de $P_0$. Pour $\chi$ ŽlŽment du groupe des caractres non ramifiŽs de $A_0$, $X(A_0)$, on note $V_\chi$ l'espace de la reprŽsentation $\pi_\chi:=i_P^G \chi$ et $V$ l'espace de sa rŽalisation compacte. On note $v$ le vecteur $K$-invariant de $V$ dont la valeur en l'ŽlŽment neutre de $G$ est Žgale ˆ 1 et $v_\chi$ l'ŽlŽment correspondant de $V_\chi$. L'espace $Wh(\chi)$ est Žgal ˆ $\C$.
Avec les notations du dŽbut de l'introduction   soit
$\xi_\chi: = \xi(P, \chi, 1)\in Wh(\pi_\chi)$ o $1\in \C= Wh(\chi)$ et soit $E_\chi $ le coefficient gŽnŽralisŽ $c_{\xi_\chi, v_\chi}$. 
  La transformŽe de Fourier Whittaker de $f\in\ccu^K$ est entirement dŽterminŽe par la fonction polynomiale sur $X(A_0)$, $F$, dŽfinie par: $$ F(\chi):= \int_{U_0\backslash G } f(g) \overline{ E_\chi (g) } dg.$$
  On notera parfois ${\hat f}$ au lieu de $F$.\\
Soit $w$ un reprŽsentant dans $K$ de l'ŽlŽment non trivial du groupe de Weyl de $A_0$. On a $w\chi= \chi^{-1}$. Par transport  de structure,  les intŽgrales d'entrelacement dŽterminent un opŽrateur d'entrelacement $A(w, \chi)$ entre $\pi_\chi$ et $\pi_{\chi^{-1}} $.\\ On dŽfinit une fonction rationnelle sur $X(A_0)$, $a$ (resp. $b$),  ˆ valeurs complexes,  par les ŽgalitŽs:  $$A(w, \chi) v_\chi= a(\chi) v_{w\chi},\>   \xi_{\chi^{-1}} \circ A(w, \chi)= b(\chi) \xi_\chi.  $$
On en dŽduit l'ŽgalitŽ de fonctions rationnelles sur $X(A_0)$, dite Žquation fonctionnelle pour $E_\chi$:
$$a(\chi) E_{\chi^{-1}}(g)= b(\chi) E_{\chi}(g), g \in G.$$ 
Utilisant l'entrelacement donnŽ ˆ l'aide de $w$ entre $i_P^G\chi$ et $i_{P^-}^G \chi^{-1}$, 
on a les ŽgalitŽs: \ber \label{AABB} $$ B(P, P^-, \chi)= b(\chi), B(P^{-}, P, \chi)=b(\chi^{-1})$$
$$A(P^-, P, \chi)v_0= a(\chi) v , A(P, P^-, \chi )v= a(\chi^{-1}) v, $$ \eer
les dernires Žtant Žcrites dans la rŽalisation compacte.
\\ Par ailleurs, les fonctions $a$ et $b$ satisfont les relations, pour $\chi $ unitaire:
$$Ê\label{bezout) }a(\chi)= a(\chi^{-1})^*, b(\chi)= b(\chi^{-1})^*,$$
la deuxime Žtant corollaire de notre Žtude du produit scalaire de paquets d'ondes (cf. ThŽorme \ref{fphihat}). 
Gr\^ace ˆ ses relations et ˆ l'Žquation fonctionnelle pour $E_\chi$,  on montre aisŽment l'identitŽ de fonctions rationnelles sur $X(A_0)$:
\beq \label{relF} a(\chi) F(\chi)= b(\chi) F(\chi^{-1}). \eeq
On va voir que si $F$ est une fonction polynomiale sur $X(A_0) $ qui vŽrifie cette relation,  elle provient de la transformŽe de Fourier-Whittaker d'un ŽlŽment $f$ de $\ccu^K$. 
Le rŽsultat d'Heiermann mentionnŽ plus haut signifie ici qu'il existe une fonction, $\zeta$,  polynomiale sur $X(A_0)$ telle que l'on ait l'identitŽ:
\beq \label{bezout} a( \chi^{-1}) \zeta (\chi)+ a(\chi)  \zeta(\chi^{-1})= 1. \eeq On note $d\chi$ la mesure invariante de masse totale 1 sur le groupe $X(A_0)_u$ des caractres unitaires non ramifiŽs de $A_0$. Si $\Lambda \in X(A_0)$, on en dŽduit une mesure sur $ X(A_0)_u\Lambda $. On fixe un tel $\Lambda$ suffisamment $P$-antidominant.
On dŽfinit $$\Phi(\chi):= a(\chi)^{-1} \zeta(\chi) F(\chi), \> f_{\Phi}(g):= \int_{ X(A_0)_u \Lambda} \Phi(\chi)E_{\chi}(g) d\chi, g\in G.$$
Indiquons sommairement pourquoi  $f_\Phi \in \ccu$  et que $F$ provient de $f_\Phi$. 
Comme $f_\Phi$ est $K$-invariante et que $G=U_0A_0 K$, il suffit, pour montrer  qu'elle est ˆ support compact modulo $U_0$,  de voir que sa restriction ˆ $A_0$ est ˆ support compact. D'aprs des propriŽtŽs gŽnŽrales des fonctions de Whittaker,  $f_\Phi(a_0)=0$ pour $a_0\in A_0$ suffisament $P_0$-dominant. Il reste ˆ voir que pour $a_0$ suffisamment $P_0$-antidominant  $f_\Phi(a_0)=0)$. Pour un tel $a_0$, $E_{\chi}(a_0)$ est Žgal ˆ son terme constant. Le calcul  de celui-ci (cf.Proposition  \ref{cteis}) et les relations  conduit ˆ l'ŽgalitŽ:
$$E_\chi(a_0)= b(\chi)  \chi(a_0)  + a(\chi)  (\chi^{-1}) (a_0).$$
de sorte que, tenant compte de la dŽfinition de $\Phi$, on a 
$$ f_\Phi(a_0) =  \int_{\Lambda X(A_0)_u} b(\chi)a(\chi)^{-1}\zeta(\chi) F(\chi)  \chi(a_0) d\chi+  \int_{\Lambda X(A_0)_u} \zeta(\chi) F(\chi) \chi^{-1}(a_0)d\chi. $$
Comme $\Lambda$ a ŽtŽ supposŽ suffisamment $P$-antidominant, un dŽplacement de contour est possible sans rencontrer de p™les pour le premier terme du membre de gauche. Alors pour $a_0$ suffisamment $P_0$-antidominant, on en dŽduit,  par un passage ˆ la limite,  que ce premier terme est nul. Pour le deuxime, on intgre une fonction polynomiale, de sorte qu'on peut dŽplacer le contour et se ramener ˆ intŽgrer sur $X(A_0)_u$. En utilisant le fait que la transformŽe de Fourier d'une fonction polynomiale sur $X(A_0)$ est ˆ support compact, on voit encore que pour $a_0$ suffisamment $P_0$-antidominant, ce deuxime terme est nul. Ainsi $f_\Phi$ est bien ˆ support compact.  \\ 
Etudions la transformŽe de Fourier-Whittaker de $f_\Phi$ et notons $F'$ la fonction polynomiale sur $X(A_0)$ associŽe ˆ celle-ci comme ci-dessus. 
D'abord tenant compte du fait que $f_\Phi$ est ˆ support compact modulo $U_0$, le calcul se ramne, par l'utilisation du centre de Bernstein au cas o $\Phi$ est polynomiale. Alors on peut se ramener ˆ une intŽgrale sur $X(A_0)_u$ dans la dŽfinition de $f_\Phi$. Le calcul de la transformŽe de Fourier-Whittaker
est alors trs semblable ˆ un calcul de Waldspurger dans sa preuve de la formule de Plancherel pour les groupes [W], gr\^ace ˆ notre introduction des sous-groupes paraboliques semi-standard. 
Tenant compte des relations (\ref{AABB}), on trouve (cf. ThŽorme \ref{fphihat}):
$$F'(\chi)= a(\chi) a( \chi^{-1})\Phi(\chi) + a(\chi^{-1}) b(\chi) \Phi( \chi^{-1}).$$
Soit encore, en tenant compte de la dŽfinition de $\Phi$:
$${\hat f}_\Phi=a(\chi^{-1}) \zeta (\chi) F(\chi) + b(\chi) \zeta(\chi^{-1})F(\chi^{-1}).$$
Tenant compte de  (\ref{relF}),  puis de la relation (\ref{bezout}), on trouve finalement $F=F'$.
\\Ce traitement de ce  cas particulier nous a conduit ˆ une reformulation du ThŽorme principal ne faisant intervenir que les sous-groupes paraboliques anti-standard (cf. Corollaire du ThŽorme \ref{PW}).\\ Les sous-groupes paraboliques semi-standard sont utilisŽs dans les dŽmonstrations  qui sont de ce fait 
 proches  des preuves de Waldspurger pour la formule de Plancherel d'Harish-Chandra pour les groupes (cf. [W]). \\
Outre le rŽsultat principal, l'un des intŽr\^ets de ce travail  est  l'adaptation  des techniques de  [W] aux  fonctions de Whittaker. Au delˆ, il offre un cadre pour l'Žtude de l'analyse harmonique sur des espaces homognes comme les espaces symŽtriques (cf. par exemple [D2] pour le cas rŽel). Concernant ceux-ci, le calcul du terme constant des intŽgrales d'Eisenstein semble \^etre l'un des obstacles majeurs ˆ surmonter. \\ Nous avons Žgalement inclus (cf. ThŽorme \ref{LM}) une rŽponse positive ˆ une conjecture de Lapid et Mao(cf. [LM], conjecture 3.5)   sur le spectre discret des fonctions de Whittaker.  Nadir Matringe (cf. [Ma], Corollaire 3.1),  a obtenu indŽpendamment une rŽponse positive ˆ cette conjecture pour certains groupes.   
\\ {\bf Remerciements} \\
Je remercie chaleureusement  Paul Mezo qui m'a donnŽ l'impulsion nŽcessaire pour lever l'hypothse de la caractŽristique nulle.\\
Je remercie vivement le referee pour ses suggestions,  notamment concernant l'introduction, l'addition de la section  sur la conjecture de Lapid-Mao   et l'appendice.  \section{Notations, dŽfinitions}
 \subsection{} 
 \setcounter{equation}{0}
Voici un systŽme de notations empruntŽ a [W] et [A]. 
Soit $F$ un corps local non archimŽdien. 
On considre divers
groupes algŽbriques dŽfinis sur
$F$ et on utilisera des abus de terminologie du type suivant: ``soit $A$ un tore
dŽployŽ" signifiera ``soit
$A$ le groupe des points sur $F$ d'un tore dŽfini et dŽployŽ
sur
$F$". Avec ces conventions, soit $G$ un groupe  algŽbrique linŽaire rŽductif et
connexe.
On fixe un tore dŽployŽ maximal, $A_0$, de $G$ et on note $M_0$ son centralisateur dans $G$.  On fixe $P_0$ un sous-groupe parabolique minimal de $G$ qui admet $M_0$ comme sous-groupe de LŽvi.  On notera $U_0$ le radical unipotent de $P_0$.\ste
Si $P$ est un sous-groupe parabolique de $G$, on dit que $P$ est semi-standard (resp. standard)  si $M_0\subset P$ (resp. $P_0\subset P$). Si $P$ est semi-standard, il  possde un unique sous-groupe de LŽvi, $M$,  contenant $M_0$.  On dit que $M$ est un sous-groupe de LŽvi semi-standard.  
 \ste L'expression  ``$P=MU$ est sous-groupe parabolique semi-standard de $G$'' signifiera que $P$ est un tel sous-groupe, que $M$ est son sous-groupe de LŽvi semi-standard et que $U$ est son radical unipotent.  On notera  $P^-=MU^-$ le sous-groupe parabolique opposŽ ˆ $P$ de sous-groupe de LŽvi $M$. Un sous-groupe 
parabolique semi-standard de $G$, $P$,  sera dit anti-standard si $P^-$ est
standard. \ste 
 Si $H$ est un groupe algŽbrique, on note $Rat (H)$ 
le groupe des  caractres algŽbriques de $H$ dŽfinis sur $F$.\\ Si $V$
est un espace vectoriel, on note $V'$ son dual et,  s'il est rŽel,  on note  $V_{\C} $
son complexifiŽ. \\ On note $A_{G}$ le plus grand tore dŽployŽ dans le centre
de $G$.  On note
$\a_{G}= Hom_{\Z} (Rat (G),\R)$. La restriction des caractres algŽbriques 
de $G$ ˆ $A_{G}$  induit un isomorphisme: 
\beq Rat(G)\otimes_{\Z}  \C \simeq Rat(A_{G})\otimes_{\Z} \C. \eeq  
On dispose de l'application canonique: \beq H_{G}: G\rightarrow
\a_{G}, \eeq
 dŽfinie par: 
\beq  \label{HG} e^{ \langle H_{G}(x), \chi\rangle   }= \vert \chi (x)\vert_{F}, \> x\in G, \chi \in
Rat (G ) .\eeq o
$\vert. \vert_F$ est la valuation normalisŽe de $F$. Le noyau de 
$H_{G}$, qui  est notŽ   $G^{1}$,  est l'intersection des noyaux des
caractres de $G$ de la forme $\vert \chi \vert_{F}$, $\chi \in Rat (G)$.
On notera
$X(G)= Hom (G/G^{1}, \C^{*})$. C'est le groupe des caractres non ramifiŽs de $G$. \ste   On a des notations similaires pour des
sous-groupes de LŽvi semi-standard. Si $P$ est un sous-groupe
parabolique semi-standard de $G$,  on notera
$\a_{P}= \a_{M_{P}}$, $H_{P}=H_{M_{P}}$. On note $\a_{0}= \a_{M_{0}}$,
$H_{0}=H_{M_{0}}$.
   On note $\a_{G, F}$, resp. ${\tilde
\a}_{G, F}$,  l'image de $G$, resp. $A_{G}$, par $H_{G}$. Alors
$G/G^{1}$ est un rŽseau isomorphe ˆ
$\a_{G,F}$.
Soit $M$ un
sous-groupe de LŽvi semi-standard.  Alors  les inclusions $A_{G}\subset
A_{M}\subset M\subset G$ dŽterminent   un morphisme  surjectif
 $\a_{M, F}\rightarrow  \a_{G, F}$, resp.  un morphisme injectif
 ${\tilde
\a}_{G, F}   \rightarrow {\tilde
\a}_{M, F}$, qui se prolonge de manire unique  en une  application linŽaire
surjective entre $\a_{M}$ et $\a_{G}$, resp. injective entre $\a_{G}$ et
$\a_{M}$.
La deuxime  application permet d'identifier $\a_{G}$ ˆ un sous-espace de
$\a_{M}$ et le noyau de la premire, $\a^{G}_{M}$, vŽrifie:
\beq \a_{M}= \a^{G}_{M}\oplus \a_{G}.\eeq
Soit $P=MU$ un \sgpss. On note $\Sigma(A_M)$ (resp. $\Sigma(P)$)
l'ensemble des racines de
$A_M$ dans l'algbre de Lie de $G$ (resp. $P$) qui s'identifie ˆ un
sous-ensemble de $\a_M'$. On note
$\Delta(P)$ l'ensemble des racines simples de $\Sigma(P)$. Si $\aa\in
\Sigma(A_M)$, on note $U_\aa$ le sous-groupe radiciel de  $U$ correspondant ˆ $\aa$.   On peut associer ˆ tout $\aa \in \Sigma(A_M) $ une coracine $\check{\aa}\in \a_M$ (cf. [A], section 3).\ste

\ber \label{lchi}   Soit  $\l \mapsto \chi_\l $ l'application $(\a_{G}')_{\C}\rightarrow X(G)\rightarrow 1$  qui, en utilisant la dŽfinition (1.1)  de $\a_{G}$ avec $Rat(G)$,
associe, ˆ
$\chi\otimes s$,  le caractre $g\mapsto \vert \chi (g)\vert
^{s}$. \eer  Le noyau est un rŽseau et cela dŽfinit sur $X(G)$ une structure de
variŽtŽ algŽbrique complexe pour laquelle $X(G)\simeq \C^{*d}$, o
$d=dim_{\R}\a_{G}$. Pour $\chi \in X(G)$, soit $\lambda \in
\a_{G,\C}^{*}$ un ŽlŽment se projetant sur $\chi $ par l'application
(\ref{lchi}). La partie rŽelle $Re\>\lambda \in \a_{G}^{*}$ est indŽpendante du
choix de $\lambda$. Nous la noterons $Re \> \chi$. Si $\chi \in Hom (G,
\C^{*})$, le caractre $\vert \chi \vert$ appartient ˆ $X(G)$. On pose
$Re\> \chi=  Re\> \vert \chi\vert $. De mme, si $\chi \in Hom(A_{G},
\C^{*})$, le caractre $\vert \chi \vert$ se prolonge de faon unique en un 
ŽlŽment de $X(G)$  ˆ valeurs dans $\R^{*+}$, que l'on note
encore  
$\vert
\chi \vert$  et on pose $Re \> \chi=  Re\> \vert \chi
\vert$.\ste Soit 
$
\> X(G)_u:=
\{
\chi
\in X(G)\vert Re
\>  \chi =0\}$ l'ensemble des ŽlŽments unitaires de $X(G)$. \\
Les notations ainsi dŽfinies s'appliquent ˆ tous les sous-groupes de LŽvi semi-standard de $G$.
\\ On choisit $K$ un sous-groupe compact maximal de $G$, dont on suppose qu'il est le fixateur d'un point spŽcial de l'appartement associŽ ˆ $A_0$ dans l'immeuble de $G$. Pour le  rŽsultat suivant, cf.  [C], Prop. 1.4.4:
\ber \label{iwa} Il existe une suite dŽcroissante de sous-groupes compacts ouverts de $G$, $K_n$, $n\in \N$,  telle que pour tout  $n\in \N^*$, $H= H_n$
est normal  dans  $K= H_0$ et pour tout sous-groupe parabolique standard de $G$, $P$,  on a: \ste 
1) $ H= H_{U^-} H_M H_U$ o $H_{U^-} = H \cap U^-$, $H_M=
H\cap M$,  $H_U= H \cap U. $\ste
   2) Pour tout $a\in A_M^-:=\{a \in A_M\vert \vert \alpha(a) \vert_F \leq 1, \aa \in \Sigma(P) \}$, $aH_Ua^{-1}\subset H_U$, $a^{-1}H_{U^-}a \subset 
K_{U^-}.$ \ste
3) Le groupe $H_M$ vŽrifie 1) et 2) relativement aux sous-groupes paraboliques de $M$ contenant $P_0\cap M$.
 \ste 
4) La suite  $H_n$ forme une base de voisinages de l'identitŽ de  l'ŽlŽment neutre de  $G$. \eer 
 On dit que  $H$  possde une factorisation  d' Iwahori par rapport ˆ  $(P,P^-)$ si 1) et 2) sont satisfaits.
 
  \subsection{Choix de mesures}
 On munit $G$ (resp. $K$) d'une mesure de Haar, $dg$, (resp. de la mesure
de Haar de masse totale 1, $dk$).  Pour tout sous-groupe fermŽ, $H$,  de $G$, on note $dh$ une mesure de Haar ˆ gauche sur $H$, dont le choix sera Žventuellement spŽcifiŽ et 
on note
$\delta_H$ la fonction module de
$H$. \ste Pour tout espace totalement discontinu $Z$, on note $C_c^{\infty}(Z)$ (resp.  $C(Z)$, resp. $C_c(Z)$)
l'espace des fonctions localement constantes, ˆ support compact  (resp. continues, resp. continues ˆ support compact) sur  $Z$ ˆ valeurs dans  
 $\C$.\ste
Soit $P=MU$ un sous-groupe parabolique semi-standard
de $G$. Alors $G=PK$ et le groupe $M\cap K$ vŽrifie relativement ˆ $M$ les mmes propriŽtŽs que $K$ relativement ˆ $G$. Pour $g\in G$, on choisit $u_P(g)\in U$, $m_P(g) \in
M$, $k_P(g) \in K$,  de sorte que $g=u_P(g)m_P(g)k_P(g)$. 
 On note $du^-$ la mesure de Haar sur $U^-$ telle que:
\beq \label{fintu} \int_{U^-} \delta_P(m_P(u^-)) du^-=1\eeq
et de mme pour $U$. On note $dk$ la mesure de Haar sur $K$ de masse totale 1.
Il existe une
unique mesure de Haar sur
$M$,
$dm$, telle que pour tout $f
\in C_c(G)$ on a (cf. [W], I (1.2)): 
 \ber \label{iwas} $$\int_{G}f(g)dg= \int_{U\times M\times K} f(umk) \delta_P(m)^{-1} dkdmdu. $$
  $$\int_{G}f(g)dg= \int_{U\times M\times U^-} f(umu^-)
\delta_P(m)^{-1}dudmdu^-$$\eer
  Soit $f$ une fonction  continue sur $K$ et invariante ˆ gauche par  $K\cap P$. Alors (cf. par exemple [K], ch. V, section 6, consŽquence 7, pour la version rŽelle): 
\ber \label{intu-}$$ \int_Kf(k)  dk= \int_{U^-}
f(k_P(u^-))\delta_P(m_P(u^-) ) du^-. $$
\eer
o l'intŽgrale du membre de droite est absolument convergente.
En particulier si $f$ est une fonction continue sur $G$   telle que
$f(umg)= 
\delta_P(m)  f(g)$, $u\in U, m\in M, g\in G $
\ber \label{intu--}$$ \int_Kf(k)  dk=  \int_{U^-}
f(u^-)du^-. $$ \eer
o l'intŽgrale du membre de droite est absolument convergente. \\
On fixe  la  mesure de Har de mase totale 1 sur  $A_G\cap K$ et sur $X(A_G)_u$,  qui s'identifie au dual unitaire de   $A_G/A_G \cap K$. On fixe sur $A_G/A_G \cap K$ la mesure de Haar duale de celle sur  $X(A_G)_u$. Des mesures ainsi choisies,  on dŽduit une mesure sur $A_G$, notŽe $da_G$. L'homomorphisme de restriction dŽtermine un morphisme surjectif de $X(G)_u$ sur  $X(A_G)_u$. On fixe sur $X(G)_u$ une mesure de Haar telle que ce morphisme prŽserve localement les mesures de Haar choisies.\\
On choisit un ensemble de reprŽsentants dans $K$, $W^G$,  du quotient du  normalisateur de $M_0$
dans $G$ par son centralisateur, $\overline{W}^G$, qui existe parce que $K$ est le fixateur d'un point spŽcial de l'appartement associŽ ˆ $A_0$. On le choisit de telle sorte qu'il contienne l'ŽlŽment neutre de $G$,  qu'on notera $1_G$ ou seulement $1$ s'il n'y a pas de confusion.
\ber \label{xpoint}Si $x\in G$ et $Y$ est une partie  de $G$,  on notera $x.Y:=\{xyx^{-1}\vert y \in Y\}$. 
De mme si $H$ est un sous-groupe  de $G$, et $(\si, E)$ est une reprŽsentation de $H$, on notera $x\si$ la reprŽsentation de $x.H$ dans $xE:=E$ dŽfinie par $x\si(xh x^{-1})= \si(h), h\in H$. \ste
On rappelle que le dual d'un espace vectoriel $V$ est notŽ $V'$. Si $T$ est une application linŽaire entre deux espaces vectoriels complexes, on note $T^t$ sa transposŽe. Si $(\si,E) $ est une reprŽsentation de $H$, on note $\si'$ la reprŽsentation de $H$ dans $E'$ dŽfinie par $\si' (h)= \si(h^{-1})^t, \>h \in H$. \ste 
Pour toute application dŽfinie sur $H$, $f$, on note:  $$(\l(h)f)(h')= f(h^{-1} h'), \> (\rho(h)f)(h')= f(h'h), \> h, h'\in H.$$\eer
 \subsection{ReprŽsentations induites} 
 Les reprŽsentations lisses de $G$ et de ses sous-groupes fermŽs seront toujours ˆ coefficients complexes. \\
  Soit $P=MU$ un sous-groupe parabolique semi-standard de $G$. 
  Si $(\si, E)$  est une reprŽsentation lisse de $M$, on l'Žtend en une reprŽsentation de $P$ triviale sur $U$. Soit $\chi \in X(M)$. On note $E_\chi$ l'espace de la reprŽsentation $\si\otimes \chi$, qu'on notera  $\si_\chi$, et on note
$i^G_{P}E_\chi$ l'espace des fonctions $v$ de $G$ dans $E$, invariantes ˆ droite par 
un sous-groupe compact ouvert de $G$ et telles que $v(mug)= \delta_P(m)^{1/2}\si_\chi (m) f(g)$ pour tout
$m\in M$, $u\in U$, $g\in G$. On note $i^G_{P}(\si_\chi)$ la
reprŽsentation de $G$  dand 
$i^G_{P}E_\chi$ par translations ˆ droite.  On note $i^K_{P\cap K}E$ l'espace des fonctions
$v$ de $K$ dans $E$, invariantes ˆ droite par un sous-groupe compact ouvert de $K$ et telles
que $f(pk)= \si(p)f(k)$ pour tout $k\in K$ et $p\in P\cap K$. La restriction des fonctions ˆ
$K$ dŽtermine un isomorphisme de $i^G_{P}E_\chi$ sur $i^K_{P\cap K}E$. On notera ${\overline
i}^G_{P}(\si _ \chi)$ la reprŽsentation de $G$ dans $i^K_{P\cap K}E$ dŽduite de
$i^G_{P}(\si _\chi)$ par transport de structure. Cette reprŽsentation sera appelŽe  la rŽalisation compacte de 
$ i^G_P(\si_\chi)$ dans cet espace indŽpendant de
$\chi$. 
Si $v \in i^K_{P\cap K}E$, on note
$v_\chi$ l'ŽlŽment de  $i^G_{P}E_\chi$ dont la restriction ˆ $K$ est Žgale ˆ $v$.   
\ste Si $\si$ est unitaire,  on munit $i^K_{P\cap K}E$ du produit scalaire dŽfini par: 
  \beq (v, v')= \int_{K}(v(k),v'(k))dk, \> v, v' \in i^K_{P\cap K}E .\eeq
 Alors, muni de ce produit scalaire, la reprŽsentation $\overline{ i}^G_P (\si_\chi)$ est unitaire pour $\chi$ unitaire et par consŽquent, par transport de structure, $i^G_P(\si _\chi)$ Žgalement.  On a alors, gr\^ace ˆ (\ref{intu--}):
\beq \label{prodscalu-}(v_\chi, v'_\chi) =\int_{U^-} (v_\chi(u^{-}), v'_\chi (u^-)) du^-. \eeq
Soit $H\subset H'$ deux sous-groupe fermŽs de $G$. Soit $(\theta, V)$ une reprŽsentation de $H$. On note $ind_H^{H'}\theta$ la reprŽsentation de $H'$ par reprŽsentation rŽgulire droite dans l'espace $ind_H^{H'}V$ des fonctions
$v$ de $H'$ dans $V$, invariantes ˆ droite par un sous-groupe compact
ouvert de $H'$, ˆ support compact modulo $H$,  et telles que $v(hh')=
\theta(h)v (h')$ pour tout $h\in H$ et $h' \in H'$.  \ste
On remarque que:
\ber\label{otimes} Si $(\chi, \C_\chi) $ est une reprŽsentation lisse de dimension 1 de $H'$, $(ind_H^{H'} \theta )\otimes \chi$ est naturellement isomorphe ˆ $ind_H^{H'} (\theta \otimes \chi_{\vert H})$, l'entrelacement Žtant donnŽ par $v \mapsto \chi v$, $v \in ind_H^{H'} V$.  \eer
 
On note $H_0(H,V)$ le quotient de $V$ par le sous-espace engendrŽ par l'ensemble des $\rho(h)v-v$, $h\in H, v \in V$. 
\ber\label{ex} Si $H$ est union de sous-groupes compacts, le foncteur qui ˆ $V$ associe $H_0(H,V)$ est exact.  \eer 
On note  $\C_{\delta_{H'}\delta_H ^{-1}}$ le $H$-module de dimension 1 donnŽ par $\delta_{H'}\delta_H^{-1}$. Alors, (cf.[ BD], Lemme 1.14 et Proposition 1.15, o, dans la preuve et les ŽnoncŽs,  les fonctions modules se sont ŽgarŽes):
\ber\label{shapiro}
$H_0(H', ind _H^{H'}V)$  est naturellement isomorphe ˆ $H_0(H, V\otimes \C_{\delta_{H'}\delta_H ^{-1}})$ .\eer 
Explicitons l'isomorphisme, en supposant en outre $H$ et $H'$ unimodulaires. Si $v\in V$ on note $p(v)$ son image dans $H_0(H,
V)$. 
 Si $v\in ind_H^{H'}V$,  l'application de $H'$ dans $H_0(H, V)$, 
$h'\mapsto p(v(h'))$ passe au quotient en une  application
localement constante ˆ support compact sur $H\backslash H'$ et l'on a: \ber \label{sha}L'application $v\mapsto \int_{H\backslash H'} p(
v(h'))dh'$ de $ind_H^{H'}V$ dans $H_0(H, V)$ passe au quotient en
un isomorphisme de $H_0(H', ind_H^{H'}V)$ avec $H_0(H, V)$.
\eer
\ste 
On remarque que,  avec les notations prŽcŽdentes:
 \beq \label{indgland} i^G_P \si= ind ^G_P(\si\otimes  \delta_P^{1/2}) .\eeq
\subsection{Fonctions dŽpendant  de reprŽsentations}
 Soit $(\pi, V)$ une reprŽsentation lisse irrŽductible de $G$ (resp. lisse unitaire irrŽductible).  On rappelle que si $\chi \in X(G)$, $\pi_\chi$ dŽsigne la reprŽsentation $\pi \otimes \chi$. On notera $[(\pi,V)]$,  ou simplement $[\pi]$,  la classe d'Žquivalence de
$\pi
$. Soit $\O=\{ [\pi_ \chi]Ê\vert \chi \in X(G)\}$ (resp. $\O_u=\{ [\pi_ \chi]Ê\vert \chi \in X(G)_u\}$. Les orbites
inertielles (resp. orbites inertielles unitaires)  de reprŽsentations  lisses  irrŽductibles (resp. unitaires irrŽductibles) sont par dŽfinition les 
ensembles du type $\O$ (resp. ${\cal O}_u$).\ste On note $\CC$ (resp. $\CC_u$)  la catŽgorie 
 dont les objets sont les reprŽsentations lisses (resp. reprŽsentations lisses  unitaires) de $G$   Žquivalentes ˆ l'une des reprŽsentations $\pi_
\chi$, $\chi\in X(G)$ (resp.  $ X(G)_u$)  et dont les flches sont les entrelacements  bijectifs (resp. entrelacements unitaires).  On appellera, par abus de langage,  objet de $\O$ (resp. $\O_u$) les objets de $\CC$ (resp. $\CC _u$).\ste
On veut dŽfinir les ``fonctions sur $\O$ (resp. $\O_u$)'', relativement ˆ
 un foncteur $\Psi$ entre $\CC$ (resp. $\CC_u$)  et la catŽgorie des espaces
vectoriels. 
Commencons par donner une dŽfinition formelle, avant de
 lui donner un sens plus concret. \ste
 On note $\E$ la catŽgorie dont les objets sont les couples $(E,e)$ o $E$ est un espace
vectoriel et
$e$ est un ŽlŽment de
$E$, avec les morphismes induits par les morphismes
 d'espaces vectoriels. Une ``fonction sur $\O$ (resp. $\O_u$)'', ou 
``fonction sur  $\O$ (resp. $\O_u$) ˆ valeurs  dans $\Psi$''  est la donnŽe d'un foncteur $\Phi$  de $\CC$ (resp. $\CC_u$) dans $\E$ tel que, pour
tout objet $\pi_1$ de $\O$  (resp. $\O_u$), 
$\Phi(\pi_1)$ est de la forme
$(\Psi(\pi_1), \phi(\pi_1))$  et  tel que, pour tout  $T$ morphisme de ${\cal C}$ (resp. ${\cal C}_u$), le  morphisme
$\Phi(T)$ soit induits par $\Psi(T)$. Dans toutes nos applications le
foncteur $\Psi$ sera tel que, si $T$ et $ T'$ sont des  isomorphismes  entre deux
objets de $\CC$ (resp. $\CC_u$), $\Psi(T)=\Psi(T')$. Par exemple, cette propriŽtŽ est
vŽrifiŽe pour le foncteur  qui ˆ $(\pi_1, V_1)$ associe $Hom_\C(V_1,
V_1)$.  On suppose cette propriŽtŽ vŽrifiŽe par $\Psi$ dans la suite. 
 Alors si $T$ est un isomorphisme (resp. entrelacement unitaire)   entre $\pi_1$ et $\pi_2$,  on a la
relation:
 \beq \label{fonct} \phi(\pi_2)=  \Psi(T) \phi(\pi_1). \eeq
 On vŽrifie aisŽment que la donnŽe  de $\Phi$ Žquivaut ˆ la donnŽe d'une fonction qui ˆ tout $\chi\in X(G)$ (resp. $X(G)_u$) associe
 $f(\pi\otimes \chi) \in\Psi(\pi\otimes \chi)  $  telle que si $T$ est un isomorphisme  (resp. entrelacement unitaire)    entre $\pi\otimes \chi$ et
$\pi\otimes
\chi_1$: 
 \beq \label{fonct'}f(\chi_1)= \Psi (T) f(\chi). \eeq 
 MalgrŽ son cotŽ trs formel,  l'intŽrt de notre dŽfinition est  de montrer que $\Phi(\pi_1) $ et $\phi(\pi_1)$ sont  dŽfinis pour toute  reprŽsentation $\pi_1$ Žquivalente 
ˆ l'une des reprŽsentations $\pi_ \chi$, $\chi \in X(G)$ (resp. $\chi \in X(G)_u$).\ste
 En pratique on ne dŽfinit pas le foncteur $\Psi$,  qui est implicite,   mais on dŽfinit $\phi$ ou $f$ et  on vŽrifie la relation (\ref{fonct}) (resp. (\ref{fonct'})). On Žcrira alors que cette relation dŽfinit  $\phi$ comme fonction sur $\O$, ou, si l'on prŽcise $\Psi$, comme fonction sur $\O$   ˆ valeurs dans $\Psi$. 
 \ste Si pour tout $\pi$ objet de $\O$ (resp $\O_u$), l'espace $\Psi(\pi_ \chi)$ s'identifie canoniquement ˆ  un espace $E_\pi$,  indŽpendant de $\chi$, on dira que: 
 \ber   \label{polo}
 La fonction sur $\O$ (resp. $\O_u$) ˆ valeurs dans $\Psi$, $\phi$,  est polynomiale  si pour tout $\si$ objet de $\O$ (resp. $\O_u$), l'application $\chi \mapsto \phi(\pi_\chi),  \chi \in
X(G)$ (resp. $X(G)_u$),  est ˆ  valeurs dans un sous-espace de dimension finie de $E_\pi$ et polynomiale (resp. 
se prolonge en une application polynomiale sur $X(G)$). On notera $Pol(\O, \Psi)$ (resp. $Pol(\O_u, \Psi)$) l'espace vectoriel des
fonctions polynomiales sur $\O$ (resp. $\O_u$) ˆ valeurs dans $\Psi$. \\
On dŽfinit de m\^eme les applications rationnelles.
 \eer
  \section{\label{termeconstant}Terme constant de fonctions de Whittaker}
  \setcounter{equation}{0}
  Un homomorphisme, $\psi$,  de $U_0$ dans $\C^*$ est un  caractre lisse non dŽgŽnŽrŽ,
 si et seulement si son noyau est ouvert et  pour tout $\alpha$ ŽlŽment de l'ensemble des racines simples de $A_0$ dans $P_0$, 
$\Delta(P_0)$, la restriction de $\psi$ au sous-groupe radiciel   $(U_0)_\alpha$ est non triviale. On fixe un
tel caractre $\psi$ dans la suite. \ste  
  Si $(\pi, E)$ est un module lisse pour  $U_0$, on notera
$(\pi_{\psi^{-1}}, E_{\psi^{-1}})$ la reprŽsentation $\pi\otimes \psi^{-1}$
de $U_0$  dans $E$.\\
  On note $\cu$  l'espace des fonctions $f$ sur $G$, invariantes ˆ droite par un
 sous-groupe compact ouvert  et telles que $f(ug)=\psi(u) f(g)$ pour $g \in G$, $u \in U_0$. On
note $\ccu$ le sous-espace de $\cu$ formŽ des ŽlŽments de $\cu$ ˆ support compact modulo $U_0$.\ste
  Soit $(\pi, V)$ une reprŽsentation lisse de $G$. 
On note $({\check \pi},{\check V})$ la reprŽsentation
 de $G$ dans le dual lisse ${\check V}$ de $V$. On note $(\pi', V')$ la
reprŽsentation contragrŽdiente de $(\pi, V)$. On note $Wh(\pi)$ l'espace
des formes linŽaires $\xi$ sur $V$ telles que $\pi'(u_0)\xi=
\psi^{-1}(u_0)\xi$ pour tout $u_0\in U_0$. En d'autres termes: 
\beq \label{wh} Wh(\pi)= (V_{\psi^{-1}})'^{U_0} \> ou \> encore \> Wh(\pi)=
(H_{0}(U_0,V\otimes \C_{\psi^{-1}}))' \eeq o $H_0$ dŽsigne l'homologie  et
$\C_{\psi^{-1}}$ dŽsigne l'espace de  la reprŽsentation de dimension 1 de  $U_0$ donnŽe par
$\psi^{-1}$. 
 On appelle les ŽlŽments de $Wh(\pi)$ les fonctionnelles de Whittaker de
$\pi$.
\ste  Si $\xi\in Wh(\pi)$ et $v\in V$, on note $c_{\xi, v}$ le
coefficient gŽnŽralisŽ dŽfini par: \beq \label{coef} c_{\xi, v}(g)=  \langle \xi,
\pi(g) v\rangle   , g\in G.  \eeq Alors
$c_{\xi,v} $ est un ŽlŽment de $\cu$. \ste
On a le rŽsultat suivant du ˆ Bushnell et Henniart [BuHen], ThŽorme 4.2 (cf. [D3], ThŽorme 5.7 pour une autre preuve): 
\ber Si $(\pi, V) $ est de longueur finie,  $Wh(\pi)$ est de dimension
finie. \eer 
Soit $H$ un sous-groupe compact ouvert de $G$. On  note $e_H$  la mesure de Haar normalisŽe de $H$, qu'on regarde comme un ŽlŽment de l'algbre de Hecke de $G$.  On fera dans cette partie rŽfŽrence ˆ [D3]. Il faut toutefois tenir compte, ˆ chaque fois,  de notre changement de point de vue: relations de covariance ˆ  gauche  et action ˆ droite dans l'induction.  
On note, pour $\varepsilon >0$,  $$A_0^-(\varepsilon):=\{ a\in A_0\vert \vert \alpha(a)\vert_F \leq \varepsilon, \aa \in \DD(P_0)\} \> \> et \>\> A_0^-= A_0^- (1).$$ Comme dans [D3], Lemme 3.1,  on voit \ber \label{H'} Pour tout sous-groupe compact ouvert $H$ de $G$, il existe un
sous-groupe compact ouvert $H'$ tel que pour toute reprŽsentation lisse de $G$,
$(\pi,V)$, tout $\xi
\in Wh(\pi)$ et tout
$v\in V^H$ on ait:
$$ \langle \xi, \pi(a)v\rangle  =  \langle  e_{H'}\xi, \pi(a) v\rangle   , a \in A^-_0$$
o $e_{H'}\xi$ est l'ŽlŽment de ${\check V}$ dŽfinit par $ \langle e_{H'}\xi, v\rangle  : =
 \langle \xi, \pi(e_{H'}) v\rangle  , \> v\in V$.\eer 
En effet prenons $H'$ comme dans (\ref{iwa}), contenu dans $H$ avec $H'_{U_0}\subset Ker \psi$. Soit  $h'\in H'$. On Žcrit $h'= u m u^-$  avec $u \in H'_{ U_0}$, $u^- \in H'_{U_0^-}$, $m \in H'_{ M_0}$. Alors $$\langle \xi, \pi(h'a)v\rangle  = \langle \xi, \pi(a  ma^{-1}u^- a)  v\rangle  $$
Mais $a^{-1} \in A_0^+$  normalise   $H'_{U_0^-}$ et $H'$ fixe $v$.  Donc
$$\langle \xi, \pi(h'a)v\rangle  = \langle \xi, \pi(a)v\rangle  , h' \in H', a \in A_0^{-}.$$
Ce qui prouve (\ref{H'}). \\
 On rappelle (cf. [D3],  Lemme 5.4):
\ber \label{fa-} Soit $H$ un sous-groupe compact ouvert de $G$. Il existe $C>0$ 
tel que, pour tout $f\in \cu$, invariante ˆ droite par $H$,  $f(a)=0$ si $\vert a^\aa \vert_{F} >C$ pour au moins un ŽlŽment $\aa$ de  $\Delta(P_0)$,  i.e. pour  $a$ ŽlŽment du complŽmentaire de $A_0^-(C).$
\eer 
Rappelons la caractŽrisation du terme constant des ŽlŽments de $\cu$ et
des fonctionnelles de Whittaker.
Soit $P=MU$ un sous-groupe parabolique standard de $G$.
Avec les notations ci-dessus, on note $(\pi_P, V_P)$ le produit tensoriel entre, d'une part 
 la reprŽsentation de $M$ dans   le  quotient de  $V$ par le 
$M$-sous-module engendrŽ  par les  $\pi(u)v-v$, $u\in U, v\in V$, et d'autre part 
la reprŽsentation  de  $M$ sur  $\C$ donnŽe par  $\delta_P^{-1/2}$. On appelle $(\pi_P, V_P)$ module de Jacquet normalisŽ de $V$ relatif ˆ $P$.   On 
note, pour $v \in V$, $j_P(v) $ ou $v_P$ sa projection naturelle dans  
$V_P$.
\ste Soit $\Theta_P$ l'ensemble des ŽlŽments de  $\DD(P_0)$  qui  sont racines de $A_0$ dans l'algbre de Lie de $M$. 
On note, pour $\varepsilon >0$:
$$A_0^-(P, < \varepsilon):=\{ a \in A_0^-\vert \vert \alpha(a)\vert_F < \varepsilon, \aa \in \DD(P_0)\setminus \Theta_P\}.$$
D'aprs [De] ThŽorme 3.4, Remarque  3.5 et  Proposition 3.6, on dispose d'une unique   application linŽaire  $Wh(\pi)\mapsto Wh(\pi_P)$, $\xi \to
j_{P^-}(\xi)$, notŽe aussi $\xi \mapsto \xi_P$ pour plus de commoditŽ,  qui vŽrifie:  
\ber  \label{jpxi}Pour tout sous-groupe compact ouvert $H$ de $G$, il existe
$\varepsilon_H>0$, indŽpendant de $P$,  avec les propriŽtŽs suivantes:
\ste Pour toute reprŽsentation  lisse
$(\pi, V)$ et 
$\xi
\in Wh(\pi)$, on a:   
$${\delta_P^{1/2}}(a)\langle \xi_P, \pi_P(a) v_P\rangle  _P= \langle \xi, \pi(a)
v \rangle, \> a \in A_0^{-}(P,<\varepsilon_H), \> v\in V^H.$$\eer 
Le terme constant  le long de $P$ d'un ŽlŽment $f$ de $\cu$ a
ŽtŽ dŽfini dans [D] DŽfinition 3. C'est un ŽlŽment
$f_P$ de
$\cmu$.   Il vŽrifie, pour tout $f$ invariant ˆ
droite par un sous-groupe compact ouvert $H$ de $G$: 
\ber \label{ffp} $$\delta_P^{1/2}(a)  f_P(a)= f(a), \> a \in
A^-_0(P, <\varepsilon_H). $$\eer

\ber \label{covct} L'application $f\mapsto f_P$ est un morphisme de $P$-modules entre $\cu$ et $\cmu$,
 o $P$ agit par reprŽsentation rŽgulire
droite sur le premier espace et $M$ (resp. $U$) agit par reprŽsentation rŽgulire
droite tensorisŽe par $\delta_P^{1/2}$(resp. trivialement)  sur le second. 
\eer \label{ctcxiv}
 Soit $ (\pi, V)$ une reprŽsentation lisse de $G$ et $\xi\in Wh(\pi)$. Avec les notations de (\ref{coef}), on a:
\beq \label{tctcoef} (c_{\xi,v})_P= c_{ \xi_P, v_P}, \> v\in V. \eeq 
\begin{lem} \label{borne}
\ste (i)   Si $(\pi, V)$ est une reprŽsentation lisse cuspidale
de $G$, pour tout $v\in V$, le coefficient gŽnŽralisŽ
$c_{\xi,v}$ est ˆ support compact modulo l'action ˆ gauche   de $U_0A_G$ sur
$G$.\\
(ii) On suppose ici $\psi$ unitaire. Si $(\pi, V)$ est une reprŽsentation
lisse et unitaire de
$G$ et
$\xi\in Wh(\pi)$, pour tout $v\in V$, le coefficient gŽnŽralisŽ
$c_{\xi,v}$ est bornŽ en module sur $G$.
\end{lem}
\dem 
(i) est donnŽ par [D3], ThŽorme 4.4 (ii). 
\\(ii) L'ŽgalitŽ $G=P_0K$ montre qu'il existe un ensemble fini, $I$,
d'ŽlŽments de $G$ tel que $G=U_0 A_0 I K$.
Pour dŽmontrer (i), il suffit donc de montrer que pour tout $v\in V$,
$c_{\xi,v}$ est bornŽ sur $A_0$. Mais si $a_0\in A_0$, 
$$c_{\xi, v}(a_0a)= c_{\xi,\pi(a_0)v}(a), a \in A_0.$$
Tenant compte de (\ref{fa-}), on voit que pour $a_0$ bien choisi,
$c_{\xi,\pi(a_0)v}$ est nul sur $A_0\setminus A_0^-$. Ainsi, on est
rŽduit ˆ prouver que pour tout
$v\in V$,
$c_{\xi, v}$ est bornŽ sur $A_0^-$.
Mais cela rŽsulte de (\ref{H'}) et du fait que tout coefficient d'une
reprŽsentation unitaire est bornŽ.\qed 

\section{Fonctionnelles et intŽgrales de Jacquet}
\setcounter{equation}{0}
\subsection{Sous-groupes paraboliques anti-standard}
Soit $P=MU$ un sous-groupe parabolique anti-standard de $G$.  Soit $(\si, E)$ une reprŽsentation lisse  de $M$. 
On  rappelle que 
pour
$\chi \in X(M)$, la reprŽsentation $({\si}_{\chi}, E_\chi)$ dŽsigne  la reprŽsentation  $ (\si \otimes \chi,
E)$ de
$M$. On note $B(M)$ ou $B$,  l'algbre des fonctions rŽgulires sur $X(M)$: c'est
l'algbre de fonctions sur $X(M)$  engendrŽe par les applications $b_m$,
$m\in M$, dŽfinies par $b_m(\chi):=\chi(m), \chi \in X(M)$. On dŽfinit
Žgalement une structure de
$(M,B)$-module sur
$E_B=  E
\otimes B$ en faisant agir $B$ par multiplication sur le deuxime facteur
et $m \in M$ par le produit tensoriel de $\si(m)$ avec la
multiplication par l'ŽlŽment 
$b_m$ de $B$.\ste
Soit encore:
$$\si_B(m)(e\otimes b) = (\si(m) e)\otimes b_mb, \> e \in E, b\in B. $$ 
On Žtend l'action de $M$ ˆ $P$ en la prenant triviale sur $U$.\\ On notera
$\bullet $ ˆ la place de $\chi$ ou $B$.  On considre $(i^G_{P}\si_\bullet, i^G_{P} E_\bullet)$  
On  notera aussi  $I_{\bullet}$ au lieu de $ i^G_{P} E_{\bullet}$ qui est un $(G,B)$-module.\ste   
 Les points (i) ˆ (iv) du Thorme suivant sont dus ˆ Rodier [R] et Casselman-Shalika [CS] (voir aussi [Sh], Proposition 3.1).  Nous en donnons une preuve qui permet de montrer  la polynomialitŽ du point (v). 
\begin{theo}\label{jacquetint}
\ste (i) On note
$J_\bullet$ = $\{
v \in I_\bullet \vert \> v \hbox{  est ˆ support contenu dans }
 PU_0
\}$ qui est un sous-$U_0$-module lisse de $I_\bullet$. Alors 
$H_0(U_0,J_\bullet \otimes \C_{\psi^{-1}})$ est naturellement isomorphe ˆ
$H_0(M
\cap U_0, E \otimes \C_{\psi^{-1}}) $ si $\bullet $ est Žgal ˆ $\chi$, et
ˆ
$H_0(M
\cap U_0, E \otimes \C_{\psi^{-1}})\otimes B $ comme $B$-module si $\bullet
$ est Žgal ˆ $ B$.
\ste (ii)   Passant au dual, cela dŽtermine un isomorphisme  entre
$Wh(\si) $ et  $Wh(J_\chi )$, $\eta \mapsto \xi^0(P,\si_\chi,  \eta) $. \\
La restriction des ŽlŽments
de $J_\chi$ ˆ $U_0$ sont ˆ support compact modulo $U_0\cap M$. Alors
on a:  
$$ \langle \xi^0 (P,\si_\chi,  \eta), v \rangle   =  \int_{ U^- } \langle 
\eta,
v (u^-)\rangle   \psi(u^-)^{-1} du^-,  
v\in J_\chi.$$  \ste    
(iii) 
L'injection naturelle de
$J_\bullet $ dans $I_\bullet$ dŽtermine un isomorphisme:
$$H_0(U_0,J_\bullet\otimes \C_{\psi^{-1}} )
\simeq H_0(U_0,I_\bullet \otimes \C_{\psi^{-1}}).$$ 
(iv) Passant aux
duaux dans (iii) et tenant compte de  (ii),
  on dispose d'un isomorphisme:
$$Wh(\si)\to Wh(i^G_{P}\si_\chi)$$ qu'on note: 
$$\eta \mapsto \xi (P,\si_ \chi,\eta).$$ La restriction de  $\xi
(P,\si_\chi, \eta)$ ˆ $J_{\chi}$ est  Žgale ˆ $\xi^0(P,\si_\chi, 
\eta)$ et elle dŽtermine entirement  $\xi
(P,\si_\chi, \eta)$.
\ste (v) On dispose de la rŽalisation compacte de $i^G_{P}\si
_ \chi$  dans un espace indŽpendant de
$\chi$, $I$. 
 On note ${\overline \xi }
(P,\si_\chi, \eta)$ la forme linŽaire obtenue sur $I$ dŽduite de $\xi(P, \si_\chi, \eta)$  par transport de structure. Alors
pour tout $v \in I$, l'application 
$\chi\mapsto \langle  {\overline \xi }(P,\si_\chi,\eta),v\rangle  $ est 
  une fonction polynomiale sur $X(M)$, i.e. un ŽlŽment de $B$.  En d'autres termes, pour tout $v \in I$ , notant pour $\chi \in X(M)$, $v_\chi$ l'ŽlŽment de $I_\chi$ dont la restriction ˆ $K$ est Žgale ˆ $v$, l'application $\chi \mapsto \langle \xi(P, \si_\chi, \eta), v_\chi\rangle $ dŽfinit un ŽlŽment de $B$. 
\end{theo} 
\dem \ste 
On commence par Žtudier $I_\bullet $ comme $U_0$-module. D'aprs la dŽcomposition de Bruhat, il n'y a qu'un 
un nombre fini de
$(P, U_0)$-doubles classes  et on peut choisir un ensemble de reprŽsentants de celles-ci dans $W^G$,   $\Omega =\{x_0,x_1,\dots, x_n\}$  ,  contenant $1$. On introduit les ensembles   $O_0=P U_0\subset O_1
\subset...  \subset O_n = G$  tels que
$O_{i+1}
\backslash O_i = P x_{i+1} U_0$. Un bon choix de l'ordre des $x_i$ permet de supposer les $O_i$ ouverts. On note $I_i = \{v \in I_\bullet
\vert supp
v \subset O_i \}$  de sorte que
$I_0 = J_\bullet$ et
$\{0\} \subset I_0 \subset I_1... \subset I_n= I_\bullet $.  On montre
(cf. par exemple [BlD], Proposition 1.17,  voir aussi [BZ], ThŽorme 5.2) que:\ber \label{filt} $$I_i / I_{i-1}
\simeq ind_{U_0 \cap {x_i.P }}^{U_0} x_i\si_{\bullet \vert
U_0
\cap x_i. P}$$ o, pour
$i=0,...,n$, 
$x_i\si  _{\bullet}$ est la reprŽsentation de $x_i.P$
dans
$E_\bullet$  dŽfinie par:
$$ x_i\si_{\bullet}(x_ipx_i^{-1})  = {  \si} _{\bullet}(p), p \in
P .$$  \eer 
Remarquons que notre  dŽfinition  des induites paraboliques diffŽre de celle de [BlD] (cf. (\ref{indgland})).  Cela implique qu'il  devrait plutot apparaitre 
 ${\tilde \si}= \si \otimes \delta_{P}^{1/2}$ au lieu de $\si$ dans le second membre de la dernire Žquation. Mais $\delta_P^{1/2}$ Žtant triviale sur les sous-groupes compacts et $U_0\cap x_i. P$  Žtant rŽunion de tels sous-groupes, on peut ignorer ce facteur. 
 \ste  Ce qui prŽcde montre en particulier:  \ber \label{Jind} La restriction des fonctions ˆ $U_0$ dŽtermine un isomorphisme de $U_0$-modules entre $J_\bullet $ et $ ind_{U_0 \cap P}^{U_0}
(E_{\bullet})$. \eer
Comme $U_0$ et ses sous-groupes fermŽs sont  rŽunion de sous-groupes compacts, ils sont unimodulaires.  Le  lemme de Shapiro (cf. (\ref{sha})) et (\ref{otimes})  impliquent donc: \ber Il  existe un isomorphisme canonique $T$ entre $H_0(U_0,J_\bullet \otimes \C_{\psi^{-1}}) $ et 
$H_0(U_0\cap P, E_\bullet \otimes  \C_{\psi^{-1}})  .$\eer Mais  $U_0 \cap P =U_0\cap M$  car $P$ est
anti-standard. D'autre part si $\chi \in X(M)$ et $u_0\in U_0\cap M$,
$\chi(u_0)=1$ et $b_{u_0}=1$. On en dŽduit (i) immŽdiatement. \ste
Pour (ii), l'assertion sur le support des restrictions des ŽlŽments de $J
$ ˆ $U_0$ rŽsulte de l'explicitation de (\ref{Jind}). Il en rŽsulte que l'intŽgrale de (ii) est bien dŽfinie,
et qu'elle  dŽfinit un ŽlŽment de $Wh(J_\chi)$. On remarque que $(U_0\cap M)\backslash U_0$ s'identifie canoniquement ˆ $U^-$.   L'explicitation de
l'isomorphisme de (i),  gr\^ace ˆ  (\ref{otimes}), (\ref{sha}), (\ref{Jind}), conduit ˆ (ii).
\ste  Maintenant on choisit
$x = x_i$ avec
$i>0$. On pose
$P_1=x. P$,  $M_1=x. M$,  $\chi_1=x\chi$ etc....  On note $(E_1)_\bullet $ l'espace de 
$(\si_1)_\bullet :=x\si_\bullet $. 
On veut calculer $H_0(U_0,I_i/I_{i-1}\otimes \C_{\psi^{-1}})$. D'aprs (\ref{sha}) et (\ref{filt}), 
 cet espace est isomorphe ˆ $H_0(U_0 \cap P_1, (E_1)_{\bullet }
\otimes \C_{{\psi^{-1}}})$. Comme $x_i$ normalise $M_0$, $P_1$ est un sous-groupe parabolique semi-standard et $U_0 \cap P_1$ est Žgal ˆ $(U_0 \cap M_1) (U_0 \cap U_1)$. Alors, d'aprs [BlD], Proposition 1.12,  $H_0(U_0 \cap P_1, (E_1)_{\bullet }
\otimes \C_{{\psi^{-1}}})$ est
isomorphe ˆ $H_0 (U_0\cap M_1, H_0({U_0\cap U_1, (E_1)_{\bullet } \otimes
\C_{{\psi^{-1}}}}))$.
\ste Montrons que $H_0({U_0\cap U_1,  (E_1)_{\bullet } \otimes \C_{{\psi^{-1}}}})$ est
rŽduit ˆ zŽro.
L'action de $U_1$ Žtant triviale sur $ (E_1)_{\bullet }$, on est rŽduit ˆ prouver que
\beq \label{nul} H_0(U_0\cap U_1, \C_{\psi^{-1}})=\{0\} .\eeq
Comme $x=x_i$ avec $i>0$, $U_0P_1$ n'est pas
ouvert et en particulier $P_1$ n'est pas anti-standard. Montrons qu'il existe
$\aa\in \DD(P_0)$ telle que
$(U_0 )_\aa$ soit contenu dans
$U_1$. Si c'Žtait faux, tous les $(U_0)_\aa$ serait contenu dans ${P_1}^-$, donc ${P_1}^-$
contiendrait $U_0$, ce qui voudrait dire que ${P_1}^-$ est standard, donc que $P_1$ est
anti-standard. On a donc prouvŽ notre affirmation. Ceci prouve que la restriction de
$\psi$ ˆ $U_0\cap U_1$ est non triviale. L'assertion (\ref{nul})  en rŽsulte
immŽdiatement en utilisant la dŽfinition. 
On a donc montrŽ que pour tout $i>0$, $H_0(U_0, I_i/I_{i-1}\otimes
\C_{\psi^{-1}}) $ est nul.
\ste Par ailleurs, $U_0$ Žtant rŽunion de sous-groupes compacts ouverts, le foncteur qui, 
ˆ tout $U_0$-module lisse $V$, associe $H_0(U_0, V\otimes \C _{\psi^{-1}}) $ est
exact (cf. (\ref{ex})). Alors, un argument de suite exacte montre (iii).
\ste(iv) est obtenu par passage au dual dans (iii).\ste
 On va prouver  (v) essentiellement comme dans [BlD], ThŽorme 2.8 (iv). A l'aide (i) et (iii), on voit que $Hom_B(H_0(U_0,I_B), B)$ est canoniquement isomorphe ˆ $Hom_B(H_0(U_0\cap M, E \otimes \C_{\psi^{-1}} )\otimes B, B)$.  
 A $\eta\in Wh(\si)$, il  correspond un unique ŽlŽment de cet espace qui   ˆ $e\otimes b\in H_0(U_0\cap M, E \otimes \C_{\psi^{-1}} )\otimes B$  associe $\eta(e)b\in B$.  On note  $\xi(P, \si_B, \eta)$  l'ŽlŽment correspondant  de $Hom_B(H_0(U_0,I_B), B)$ dans l'isomorphime ci-dessus. \ste 
Pour tout $B$-module, ou morphisme de $B$-module, on dispose de la spŽcialisation en tout ŽlŽment $\chi $ de $X(M)$, qu'on
note avec $\chi$ en un indice infŽrieur. Le spŽcialisŽ de $I_B$ est $I_\chi$, notamment.
Montrons que   le spŽcialisŽ en $\chi$ de $\xi(P, \si_B, \eta)$ est $\xi(P, \si_ \chi,
\eta)$.  D'aprs le point (iii), il suffit pour prouver cette ŽgalitŽ, d'Žtudier la restriction de ces formes linŽaires ˆ $J_\chi$.
On conclut en explicitant les isomorphismes gr\^ace ˆ (i) et (ii).  
\ste La rŽalisation compacte de $I_B$ se fait dans l'espace ${\overlineÊI}\otimes B$. Si
$v\in I$,  $\xi(P, \si_B, \eta)(v\otimes 1)$ est un
ŽlŽment de $B$ dont la valeur en $\chi$  est Žgale ˆ $\xi(P, \si_
\chi, \eta)(v)$, d'aprs ce que l'on vient de voir. Cela prouve (v).\qed
\begin{prop} \label{xicont}  On suppose $\si$ unitaire. Soit $\eta \in Wh(\si).$\\
(i) Si $\chi \in X(M)$ est tel que  $Re (\chi \delta_P^{-1/2})$ soit strictement $P$-dominant,  la
fonction sur $G$ ˆ valeurs dans  $E'$, 
$\chi\mapsto \tilde{\xi}(P, \si_\chi, \eta)$,    dŽfinie
par:
$$\langle \tilde{\xi}(P, \si_\chi, \eta)(umu^-),e \rangle = \psi(u^-)^{-1}\langle \eta, 
(\chi^{-1}  \delta_P^{1/2})(m)\si(m^{-1})e\rangle  , e \in E, u\in U, m\in M, u^- \in U^-$$
et $$ \tilde{\xi}(P, \si_\chi, \eta)(g)=0, g\notin UMU^-, $$
est faiblement continue, i.e. pour tout $e\in E$, l'application $g \mapsto \langle \tilde{\xi}(P, \si_\chi, \eta)(g),e \rangle$ est continue sur $G$.
\\ (ii) Pour tout $v\in I_\chi$ et $ \chi$ comme en (i), on a:  $${\xi}(P,
\si_\chi, \eta)(v)= \int_K \langle \tilde{\xi}(P, \si_\chi, \eta),
v(k)\rangle   dk$$
ainsi que:
$${\xi}(P,
\si_\chi, \eta)(v)=  \int_{ U^- } \langle 
\eta,
v (u^-)\rangle   \psi(u^-)^{-1}du^- ,$$ 
l'intŽgrale Žtant absolument convergente
\end{prop}
\dem Nous aurons besoin du Lemme suivant pour prouver (i).
\begin{lem} 
 Si $umu^-\in UMU^-$ tend vers un  
ŽlŽment  du complŽmentaire de 
$UMU^-$ dans $G$,  alors  
pour tout ŽlŽment $\nu$ de $\a_M^*$ strictement $P^-$-dominant, 
$e^{\nu(H_M(m))}$ tend vers zŽro.

\end{lem}
\dem  
Nous allons utiliser des reprŽsentations rationnelles de $G$.
Fixons quelques notations supplŽmentaires.\\
On note $\underline{G}$ le groupe algŽbrique dont $G$ est le groupe des points sur $F$. On utilisera des notations similaires pour les sous-groupes de $G$. Soit $T$ un $F$-tore  maximal de $\underline{G}$ contenant $\underline{A}_0$, $B$ un sous-groupe de Borel de $\underline{G}$, contenant $T$ et contenu dans $\underline{P}^-_0$.  On note $\Sigma(T)$ l'ensemble des racines de $T$ dans $\underline{G}$. On note $\Lambda(T)$ (resp.  $\Lambda(T)_{rac}$) le rŽseau des poids (resp. le rŽseau des racines) de $T$ relatif ˆ $G$. On note $\Gamma$ le groupe de Galois de $F$ qui agit sur ces rŽseaux. On note $\Lambda ^+ $ l'ensemble des poids dominants relatifs ˆ $B$.  \\On note $\Lambda^+_M$ l'ensemble des ŽlŽments $\lambda$ de $\Lambda^+$ tels que $G$ admette une reprŽsentation rationnelle irrŽductible, de dimension finie,  de plus haut poids $\lambda$, $(\pi_\l, V_\l)$,  telle que , notant $v_\l$ un vecteur non nul de poids $\l$ on ait:
\ber  Le vecteur $v_\l$ se transforme sous un caractre rationnel de $M$ notŽ ${ \tilde \lambda}$. De plus, il existe  $v'_\l \in V'_\l$ invariant par $U$ et tel que le coefficient $c_\l (g)= \vert \langle \pi_\l(g)v_\l, v'_\l\rangle \vert_F$ soit nul sur le complŽmentaire dans $G$ de $UMU^-$.  On note $l $ l'Žlement de $\a_M^*$ tel que pour tout $m\in M$, $\vert { \tilde \lambda}(m)\vert _F = e^{l(H_M(m))}$. \eer
On remarque que:
\ber\label{mult} $\LL ^+_M$ est stable par multiplication  par la restriction ˆ $T$ des caractres rationnels  de $G$.  \eer
Soit $\beta$ un ŽlŽment de l'ensemble des racines de $A_0$ dans $G$, $\Sigma(A_0)$. 
On note $$\tilde{\beta}: = \sum_{ \aa \in \Sigma(T), \aa_{\vert A_0}= \beta} \aa. $$
On voit facilement qu'il existe $n\in \N^*$ tel  que,  pour tout $\beta \in \Sigma(A_0)$,  il existe $n_\beta\in \N^* $ tel que $ n_\beta\tilde{\beta}_{\vert A_0}= n \beta$. Notons $n.\Lambda_{rac}(A_0)$ le rŽseau engendrŽ  par les $n_\beta \tilde{\beta}$, $\beta \in \Sigma(A_0)$. 
   Par construction, les ŽlŽments de $n.\Lambda_{rac}(A_0)$ sont invariants sous $\Gamma$ et ŽlŽments de $\Lambda_{rac}(T)$. On note que  tout ŽlŽment de $n.\Lambda_{rac}(A_0)$ est invariant par le groupe de Weyl de $\underline{M}_0$ relatif ˆ $T$, $W(\underline{M}_0 , T)$.  Notons $\LL(A_0)$ le rŽseau des poids  de $A_0$.  Comme l'ensemble des restrictions ˆ $A_0$ des ŽlŽments de  $n.\Lambda_{rac}(A_0)$ contient $n \LL_{rac}(A_0)$, o $\LL_{rac}(A_0)$ dŽsigne le rŽseau des racines de $A_0$ dans $G$,   il 
 existe $n'\in \N^*$ tel que l'ensemble des restrictions ˆ $A_0$ de $n.\Lambda_{rac}(A_0)$ contienne $n'\LL(A_0)$. \\ 
   Par ailleurs un ŽlŽment de $n.\Lambda_{rac}(A_0)$ est dans $\Lambda^+$ si et seulement si  sa restriction ˆ $A_0$ est un ŽlŽment de l'ensemble, $\LL ^+(A_0)$,  des poids $U_0^-$-dominants de $A_0$.
   Donc :
\ber \label{rest} L'ensemble des restrictions ˆ $A_0$ des ŽlŽments de $n.\Lambda_{rac}(A_0))^+:= n.\Lambda_{rac}(A_0) \cap \LL^+$ contient  $n'\LL^+(A_0)$ \eer
  Soit $\l \in (n.\Lambda_{rac}(A_0))^+\subset \LL(T)_{rac} \cap \LL^+$. On dŽduit  de [T],  ThŽorme 3.3 et Lemme 3.2, l'existence d'une reprŽsentation rationnelle de $G$ , irrŽductible de plus haut poids $\l$, $(\pi_\l, V_\l)$.
 Soit $v_\l\in V_\l$ un vecteur non nul de poids de poids $\l$. Montrons qu'il se transforme sous un caractre rationnel de $M_0$. On peut pour cela passer ˆ la cl\^oture algŽbrique. L'invariance de $\lambda$ par $W(\underline{M}_0 , T)$, le fait que l'espace de poids $\l$ soit de dimension 1 (cf. [Hu], Proposition 31.2) et la dŽcomposition de Bruhat de $\underline{M}_0$ permet de conclure. 
 \ber On note $ \LL^{++}_M$ l'ensemble des ŽlŽments de $(n.\Lambda_{rac}(A_0))^+$  dont la restriction ˆ $A_0$ est orthogonale aux racines de  simples, pour $U_0^-$, de $A_0$ dans $M$, ce qui est automatique, et non orthogonale aux autres racines simples de $A_0$ dans $U^-_0$  . \eer 
 On va montrer:
 \beq \label{inclus}  \LL^{++}_M \subset \LL ^+_M .\eeq
 Soit $\lambda \in \LL^{++}_M$. Il faut d'abord montrer que $v_\l$ se transforme sous un caractre rationnel de $M$, qu'on notera ${\tilde \lambda}$. On note $\l_0$ la restriction de $\l$ ˆ $A_0$. On  prouve de manire analogue ˆ la Proposition  31.2 de [Hu], en utilisant ici la densitŽ de $U_0^-M_0U_0$ dans $G$,  que  le sous-espace de poids $\l_0$ sous $A_0$ dans $V$est de dimension 1.  On voit de mme que les poids de $A_0 $ dans $V_\l$ sont de la forme $\mu = \l_0 + \sum_{\beta \in \Delta(P_0) }c_\beta \beta$, o les $c_\beta$ sont ŽlŽments de $\N$ (on rappelle que $B$ est contenu dans $\underline{P}^-_0$). Par ailleurs le goupe de Weyl,  $\overline{W}^{M}$,  de $M$ relatif ˆ $A_0$ fixe $\l_0$ car $\l\in \LL^+_M$. On achve de prouve  notre assertion sur $v_\l$ en utilisant la dŽcomposition de Bruhat de $M $ relative ˆ $P_0^- \cap M$. \\
 On considŽre maintenant l'hyperplan de $V_\l$ engendrŽ par les sous-espaces de poids sous $A_0$ pour des poids distincts de $\l_0$.  Gr\^ace ˆ l'information obtenue ci-dessus sur les poids de $A_0$ dans $V_\l$ et ˆ la Proposition 27.2 de [Hu], on voit que cet hyperplan   est stable par  $U$. En consŽquence la forme linŽaire sur $V_\l$, $v'_\l$,  nulle sur cet hyperplan et valant 1 sur $v_\l$ se transforme sous un caractre rationnel du groupe unipotent $U$. Elle est donc invariante par $U$. \\
ConsidŽrons la fonction $c_\lambda$ sur $G$, ˆ valeurs rŽelles,
dŽfinie par:
$$c_\lambda(g)= \vert \langle  \pi_\lambda(g)v_\lambda, v'_\lambda\rangle  \vert _F, \> g \in G. $$
Montrons que  l'hypothse sur $\lambda$ implique que   $c_\lambda $ est nulle
en dehors de $UMU^-$. En effet, d'aprs la dŽcomposition de Bruhat, un
ŽlŽment de ce complŽmentaire s'Žcrit $g=u wmu^-$,  o $w\in W^G$  reprŽsente un
ŽlŽment du groupe de Weyl  $\overline{W}^G$  qui n'appartient pas  au groupe de Weyl
$\overline{W}^M$.
Alors $c_\lambda(g)$ est  un multiple de $\vert \langle 
\pi_{\lambda}(w)v_\lambda, v'_\lambda\rangle   \vert_F$ .  Mais
$\pi_\lambda(w)v_\lambda$ est de poids $w  \l_0$ sous $A_0$. Ce poids est
distinct de $\lambda_0$   car $\lambda \in \LL^{++}_M$. On en dŽduit que $\vert \langle 
\pi_{\lambda}(w)v_\lambda, v'_\lambda\rangle    \vert_F$ est nul.  Donc  $c_\lambda(g)$ est nul, comme dŽsirŽ.  Ceci achve de prouver (\ref{inclus}). 
\ste Soit  $ \l \in \LL^+_M$.  La fonction $c_\lambda$ est continue sur $G$. On en dŽduit  que, si $umu^-$ tend vers un un ŽlŽment du
complŽmentaire de $UMU^-$ dans $G$, $\vert { \tilde \lambda}(m)\vert _F= e^{l(H_M(m)}$ tend vers zŽro. 
Soit $\nu $ comme dans l'ŽnoncŽ. Notons, pour $\beta$ racine simple de $A_0$ dans l'algbre de Lie de $P_0^-$, $\delta_\beta$ le poids fondamental de $A_0$ relatif ˆ $U_0^{-}$ correspondant. D'aps (\ref{inclus}) et (\ref{rest}), on voit que toute combinaison linŽaire ˆ coefficients dans $\N^*$  des $n'\delta_\beta$, o $ \beta$ n'est pas une racine de $A_0$ dans $M$, est Žgale ˆ un $l$ comme ci-dessus. En prenant tous les coefficient Žgaux ˆ un sauf celui correspondant ˆ un $\beta$,  que l'on prend trs grand,  et en faisant varier $\beta$ puis en tenant compte de (\ref{mult}), on voit que $\nu$ s'Žcrit comme combinaison linŽaire ˆ coefficients  positifs  d'ŽlŽmŽnts $l$ pour $\l  \in \LL^+_M$.  Le Lemme en rŽsulte. \qed
Alors le point (i) de la Proposition rŽsulte immŽdiatement du Lemme
prŽcŽdent et du fait que pour tout $e\in E$, la fonction sur $M$,
$m\mapsto \langle \si'(m) \eta, e\rangle  $ est bornŽe car $\si$ est unitaire(cf. Lemme \ref{borne} (ii)).\\
Prouvons (ii). D'abord l'ŽgalitŽ: 
$$\int_K \langle \tilde{\xi}(P, \si_\chi, \eta),
v(k)dk=   \int_{ U^- } \langle 
\eta,
v (u^-)\rangle   \psi(u^-)^{-1}du^-$$
rŽsulte de la formule intŽgrale (\ref{intu--}) et de la la continuitŽ montrŽe en (i), les intŽgrales Žtant de plus absolument convergentes. Le deuxime membre de cette ŽgalitŽ, lorsque $v$ varie, dŽfinit  une
fonctionnelle de Whittaker sur $I_\chi$, comme le montre de
simples changement de variables, en  Žtudiant  sŽparŽment la transformation
selon $U_0\cap M$ et $U^-$ de cette forme linŽaire. Celle-ci possde les
propriŽtŽs caractŽristiques de $\xi(P, \si_\chi, \eta)$ (cf. ThŽorme
\ref{jacquetint} (i) et (ii)). (ii) en rŽsulte.\qed 
\begin{prop} \label{xisisi1} Soit  $P$ un sous-groupe parabolique anti-standard de $G$. \\
(i) Soit $\O$ l'orbite inertielle d'une reprŽsentation lisse irrŽductible de $M$.  Si $(\si, E)$ est  un objet de $\O$   et $\eta \in Wh(\si)$, on dŽfinit 
  $\xi(P, \si, \eta)$ en idenfiant $\si$ ˆ  $\si \otimes 1$. La correspondance $\eta \mapsto \xi(P, \si, \eta)$ est une application linŽaire bijective entre  $Wh(\si) $ et  $Wh (i^G_P \si)$, que l'on note $\phi(\si)$.  \\  
Soit $(\si, E)$ et $(\si_1, E_1)$ deux objets de $\O$  Žquivalents. Soit   
$T: E\mapsto E_1$ un opŽrateur d'entrelacement bijectif entre $\si$ et $\si_1$. La transposŽe de $T$, $T^t$,  dŽtermine une bijection notŽe encore $T^t$, de $Wh(\si_1)$ sur $Wh(\si)$.  On note $ind T$,
l'opŽrateur d'entrelacement induit par $T$ entre
$i^G_{P}\si$ et $i^G_{P}\si_1$, qui est simplement la composition des
applications de $G$ dans $E$ avec $T$. Alors on a:
\beq\label{fonctxi} \xi(P, \si_1, \eta_1)= \xi(P, \si, T^t \eta_1)\circ (ind T)^{-1}, \eta_1 \in Wh(\si_1) . \eeq 
Cela dŽfinit $\phi$ comme fonction sur $\O$ ˆ valeurs dans $Hom(Wh(.), Wh(i^G_P)(.))$, notŽ aussi $Hom(Wh, Wh(i^G_P))$ (cf.  (\ref{fonct})).\\
(ii) En particulier, on a:
\ber \label{ximsi} Si $\si_1=m\si$ avec $m\in M\cap K$, l'application  $T=
\si(m^{-1})$ entrelace $\si$ et $\si_1$, $T^t =\si'(m) $ est une   bijection
de $Wh(\si_1)$ sur $Wh(\si)$. De plus 
$ind T$ est la multiplication par $\si(m^{-1})$, donc est Žgal ˆ  $\l(m)$  (cf. (\ref{xpoint})  pour les notations), car $\delta_P^{1/2} $ est Žgale ˆ 1 sur les ŽlŽments de $M\cap K$. La formule  ci-dessus se lit donc dans ce cas:
$$\xi(P, m\si, \eta_1) = \xi(P,\si,  \si'(m) \eta_1) \circ
\l(m^{-1}), \eta_1\in Wh(\si_1) $$
ou bien 
$$\xi(P, m\si, \si'(m^{-1})\eta) = \xi(P,\si,   \eta) \circ
\l(m^{-1}), \eta \in Wh(\si).$$\eer
\end{prop}
\dem
Cela  rŽsulte de  la caractŽrisation de
$\xi(P,\si, \eta)$ donnŽe par le ThŽorme \ref{jacquetint} (iv).
\qed \ste 
Soit $(\si, E)$ un objet de $\O$ et $e\in E$. Soit $H$ un sous-groupe compact ouvert de 
$G$ contenu dans $K$ possŽdant  une factorisation d'Iwahori par rapport ˆ $(P, P^-)$ (cf. (\ref{iwa}))  et tel que $e$ soit invariant par $H_M$.  On suppose en outre que $H$ est
assez petit, de sorte  que  $H_{U^-}$  soit contenu dans $Ker \psi$. 
On
dŽfinit une application de $G$ dans $E$, 
 $v_{e, \si}^{P,H}$, par: \ber \label{vsi} $v _{e, \si}^{P,H}
(umh_{U^-})=\delta_P^{1/2}(m)\si(m)e $ si
$h_{U^-}\in H_{U^-}$ et $m\in M, u\in U$.
\\ $v _{e, \si}^{P,H}(g)= 0 $ si $g\notin PH= PH_{U^-}.$  \eer Comme $e $ est $H_M$-invariant, 
$v_{e, \si}^{P,H}$ est invariante ˆ droite par $H$. C'est un ŽlŽment de
$i^G_{P}E$ et mme de $J_1$, avec les notations du ThŽorme 1.  On remarque que: 
\ber \label{vsik} Pour $\chi\in X(M)$,  la restriction de $v_{e, \si_\chi}^{P,H}$ ˆ
$K$ ne dŽpend pas de $\chi$.   \eer 
Notons $$vol(H_{U^-})= \int_{H_{U^-}} du^-$$
o $du^-$ est la mesure de Haar sur $U^-$ choisie en Ê(\ref{fintu}). En utilisant le ThŽorme \ref{jacquetint} (ii) et
(iv), on voit que: 
\ber  \label{xiv} $$\langle  \xi(P, \si,  \eta), v_{e, \si}^{P,H}\rangle  = vol(H_{U^-}) \langle \eta,e\rangle  , $$\eer
\begin{lem} \label{indetxi}{\em Induction par Žtage pour les fonctionnelles de Jacquet}\\ 
 Soit $P_1=M_1U_1$  sous-groupe parabolique de $G$ contenant le sous-groupe parabolique anti-standard $P$. On note $(\si_1, E_1)=(i^{M_1}_{P\cap M_1} \si,
i^{M_1}_{P\cap M_1}E)  $. 
  On identifie $i^G_P\si $ ˆ $i^G_{P_1}\si_1$. Alors si $\eta \in Wh(\si)$, $\eta_1:=\xi(P\cap M_1, \si, \eta)$ est un
ŽlŽment de $Wh (\si_1)$ et on a:
  $$\xi(P, \si, \eta)= \xi(P_1, \si_1, \eta_1).$$
  \end{lem}
\dem D'aprs le ThŽorme \ref{jacquetint} (iv) et (ii),  il suffit de voir que pour tout $v \in i^G_PE$ ˆ support dans $PU^-$ on a:
\beq \label{xixi1}\langle \xi(P, \si, \eta),v\rangle  = \langle \xi(P_1, \si_1, \eta_1),v_1\rangle  , \eeq o 
$v_1$ est  l'ŽlŽment de $i^G_P E_1$ correspondant ˆ $v$ dans l'identification de  $i^G_PE$ ˆ $i^G_{P_1}E_1$.
Le premier membre de l'ŽgalitŽ ˆ dŽmontrer est donnŽ  par le ThŽorme \ref{jacquetint} (ii) et (iv):
$$\langle  \xi(P, \si, \eta), v\rangle   =  \int_{ U^- } \langle 
\eta,
v (u^-)\rangle   \psi^{-1}(u^-)^{-1}du^- .$$  
De m\^eme $v_1$  est ˆ support dans $ P_1U_1^-= P_1U_0$.  Cela permet de calculer le
deuxime membre:
$$\langle \xi(P_1, \si_1, \eta_1), v_1\rangle  = \int_{U_1^-}  \langle \eta_1, v_1(u_1^-) \rangle  \psi^{-1}(u_1^-) ^{-1}du_1^- .$$
Mais $v_1(u_1^-)$, comme ŽlŽment de $E_1= i_{P\cap M_1} ^{M_1}E$ est ˆ support dans $(P\cap M_1)(U^-\cap M_1) $. En utilisant ˆ nouveau le  ThŽorme \ref{jacquetint} (ii) et (iv),   on exprime $\langle \eta_1, v_1(u_1^-)\rangle  $.
Le ThŽorme de Fubini permet de conclure ˆ l'ŽgalitŽ (\ref{xixi1}).  \qed
 \subsection{ Sous-groupes paraboliques semi-standard}
On rappelle que $W^G$ dŽsigne un ensemble de reprŽsentants dans $K$ du groupe de Weyl, ${\overline W}^G$,  de $G$ par rapport ˆ $M_0$.
 Si $M$ est un sous-groupe de LŽvi d'un sous-groupe parabolique semi-standard de $G$, $P$, on note $W^M= W^G\cap M$ qui est un ensemble de reprŽsentants dans 
$M$ de ${\overline W}^M$.  La longueur des ŽlŽments de $\overline{W}^G$ est dŽterminŽe par le choix de $P_0$. \\On suppose en outre ici $P$ anti-standard. 
 Il existe un  ensemble de reprŽsentants de $\ \overline{ W}^G / \overline{ W}^M $,  $\overline{ W}_M$, dans $\overline{W}^G$ tel que (cf. [War], Proposition 1.1.2.13):  \ber   \label{wM}Tout ŽlŽment $w$ de $\overline{W}^G$ s'Žcrive  sous la forme $w_Mw^M$, avec  $w_M \in \overline{ W}_M$, $w^M \in \overline{ W}^M$,   et tel que la longueur de $w$ soit Žgale ˆ la somme des longueurs de $w_M$ et $w^M$. \eer 
 \ber \label{wP} Soit $P=MU$ un sous-groupe parabolique semi-standard de $G$. On notera
$w_P$ ou parfois seulement $w$,  s'il n'y a pas d'ambiguitŽ,  l'ŽlŽment $w_P$ de $G$ 
tel que
$w_P^{-1}\in W^G$,  $P'= w_P.P$ soit anti-standard   et tel  que $w_P^{-1}$
reprŽsente l'ŽlŽment du groupe de Weyl  de longueur minimum dans
$w_P ^{-1}{\overline W}^{M'}= {\overline W}^{M}w_P^{-1}$. L'unicitŽ de $w_P$ rŽsulte du fait que deux sous-groupes paraboliques anti-standard de $G$ conjuguŽs sont Žgaux.  \eer Alors, avec les notations de  (\ref{xpoint}), $w_P.M$
est le sous-groupe de LŽvi de
$w_P.P$.\\ Soit  $(\si,
E)$ une reprŽsentation lisse  de
$M$ et $w \in W^G$.  On dispose de l'isomorphisme $\l(w): i^G_{P}E\mapsto
i^G_{w.P}w E$ entre les reprŽsentations $ i^G_{P}\si$ et $
i^G_{w.P}w\si$  qui ˆ 
$v$ associe $v_{w}$, o $v_w(g)= v(w^{-1}g)$ pour $g\in G$. 
Notons, que comme $w\in K$, pour $v \in i^K_{P\cap K}E$ et tout $\chi \in X(M)$, la restriction de 
$\l(w) v_\chi$ ˆ $K$ est Žgale ˆ $\l(w)v$. On voit aussi que si $\si$ est
unitaire,
$\l(w)$ est unitaire.  \ste
\begin{defi}  \label{whp} On dŽfinit: 
  $$Wh(P, \si):= Wh(w_P\si)$$
$$\xi(P, \si, \eta):= \xi(w_P .P, w_P\si, \eta)\circ \l (w_P),  \eta \in Wh(P, \si).$$
 \end{defi}
Alors: 
\ber  \label{bijetaxi}L'application $\eta \mapsto \xi(P, \si,  \eta)$ est une bijection entre $Wh(P, \si)$  et  $Wh(i^G_{P}\si)$($\subset  (i^G_{P}E)'$) . \eer 
On a:
\ber \label{xipolbis}Pour tout $v \in i^K_{P\cap K}E$, l'application $\chi\mapsto \langle  \xi(P, \si _ \chi, \eta), v_\chi\rangle  $
est polynomiale en $\chi \in X(M)$, o $v_\chi$ est l'ŽlŽment de l'espace de $i^G_P\si_\chi$ dont la restriction ˆ $K$ est Žgale ˆ $v$. \eer  En effet cela rŽsulte de la DŽfinition \ref{whp} et  du  ThŽorme \ref{jacquetint} (iv).\\
Soit $(\si, E)$ et $(\si_1, E_1)$ des reprŽsentations  lisses  Žquivalentes de  $M$. Soit   
$T: E\to E_1$ un opŽrateur d'entrelacement bijectif entre $\si$ et $\si_1$. Comme $T$ entrelace aussi $w_P \si$ et $w_P \si_1$, la transposŽe de $T$, $T^t$ dŽtermine une bijection, notŽe encore $T^t$, de $Wh(P, \si_1)$ sur $Wh(P, \si)$. Alors,  on dŽduit de la DŽfinition \ref{whp} et de la Proposition \ref{xisisi1} que:
\beq\label{fonctxip} \xi(P, \si_1, \eta_1)= \xi(P, \si, T^t \eta_1)\circ (ind T)^{-1}, \eta_1 \in Wh(P, \si_1).\eeq 
Reformulons la DŽfinition \ref{whp} en posant $s=w_P^{-1}$, $Q=  w_P.P$, de sorte que $Q$ est anti-standard et $P=s.Q$,  et en changeant $ \si$ en $s\si$:
\ber  \label{whps} Si $Q$ est un sous-groupe parabolique anti-standard de $G$,  de sous-groupe de LŽvi $M_Q$, et si $s \in W^G$ est de longueur minimum
dans $sW^{M_Q}$, on a: 
$$Wh(s.Q, s\si)= Wh(Q, \si)$$ 
$$\xi(s.Q,s\si, \eta)= \xi(Q,\si, \eta)\circ\l(s^{-1}).$$\eer
\subsection{ IntŽgrales de Jacquet}
\begin{defi} \label{eis} Soit $P=MU$ un sous-groupe parabolique semi-standard de $G$ et $(\si, E)$ une reprŽsentation lisse de $M$.
Avec les notations prŽcŽdentes,  on dŽfinit des ŽlŽments  $E^G_P(\si,  v, \eta)$ de $ \cu$, appelŽes  intŽgrales de Jacquet,  par: 
$$E^G_P(\si,  \eta, v)(g) =\langle \xi(P, \si,  \eta), i^G_P\si (g) v\rangle  , v  \in i^G_PE, \eta \in Wh(P,\si). $$\end{defi}
On dŽduit de (\ref{xipolbis}) que: 
\ber \label{eispol} Si $\chi\in X(M)$, $v \in i^K_{K\cap P} \si$, on  rappelle que  $v_\chi$ est l'ŽlŽment de $i^G_P \si_\chi $ dont la restriction ˆ $K$ est $v$. \\Alors, pour tout $g \in G$,  l'application $\chi \mapsto E^G_P(\si_\chi,  \eta, v_\chi)(g)$ est polynomiale en $\chi \in X(M)$.\eer
\section{Fonctionnelles de Jacquet et intŽgrales d'entrelacement}
\setcounter{equation}{0}
\subsection{IntŽgrales d'entrelacement}
Soit $P=MU$, $P'=MU'$ deux sous-groupes paraboliques semi-standard de
$G$ de sous-groupe de LŽvi $M$. Soit $\O$ l'orbite inertielle d'une reprŽsentation lisse irrŽductible de $M$.
\ber \label{intertw}
Il existe une fonction rationnelle dŽfinie sur $\O$,  $A(P', P,.)$ ˆ valeurs dans $Hom_G(i^G_P., i^G_{P'}. )$
avec les propriŽtŽs suivantes:\\
 Pour tout  $(\si, E)$ objet de $\O$,  il existe $R\in \R$
tel que pour tout $\chi \in X(M)$ vŽrifiant $\langle Re \chi, \alpha\rangle    > R$ pour tout  
$\alpha\in
\Sigma(P)\cap \Sigma({P'^-})$ on ait:$$\langle  (A(P', P, \si_\chi)v)(g), {\check e}\rangle   =
\int_{ U\cap U'\backslash U'} \langle v(u'g), {\check  e}\rangle   du',\> vÊ\in i^G_{P}V_{\chi}, {\check e}\in
{\check E}, $$ l'intŽgrale Žtant absolument convergente.\\
\eer 
La rationalitŽ s'entend dans le sens suivant:
\ber \label{ratintertw}
\ste Il existe une fonction
polyn\^ome sur $X(M)$ non nulle, $b$, telle que pour
tout $v\in i^K_{K\cap P}V$, l'application qui ˆ $\chi\in X(M)$ satisfaisant la condition ci-dessus associe la restriction ˆ $K$ de 
$b(\chi)A(P', P,\si_\chi)(v_\chi)$ est ˆ valeurs dans un
espace vectoriel de dimension finie de $i^K_{P\cap K}E$ et  se prolonge de faon polynomiale en $\chi\in X(M)$ (cf. 
[W], ThŽorme IV.1.1). \\ Si $\si $ est tempŽrŽe, on peut prendre
$R=0$ (cf. [W], Proposition IV. 2.1) .\eer 

Il existe une application rationnelle sur $\O$ ˆ valeurs dans $\C$, $j$,  telle que pour tout sous-groupe parabolique semi-standard de $G$, $P$,  de sous-groupe de LŽvi
$M$, on ait (cf.  [W], IV.3 (1)):
\ber   \label{hom} Pour $\si $ objet de $\O$ tel que $A(P, P^-, \si) A(P^-,P, \si)$ soit dŽfini, cet opŽrateur est l'homothŽtie de rapport $j(\si)$.\eer 
On a (cf. [W] IV.3 (3)):
\ber \label{jw} Si $w\in W^G$, $j(w\si)=j(\si)$.Ê\eer
 D'autre part,   on obtient facilement un analogue de [W] IV.1 (11)  pour les
adjoints des intŽgrales d'entrelacement.  Cela  conduit ˆ
un analogue de l.c. IV. 3 (2) que l'on exprime sous la forme suivante: 
\ber \label{reelunit} Le nombre $j(\si)$ est rŽel si $\si$ est unitaire.Ê\eer  Si $\alpha$ est
un ŽlŽment de l'ensemble
$\Sigma_{red}(P)$ des racines rŽduites de
$\Sigma(P)$, on note $A_\aa$ la composante neutre du noyau de $\alpha$ dans $A_M$ et $M_\aa$ le centralisateur de $A_\aa$. On note  $j_\alpha(\si)$ le terme analogue ˆ $j(\si)$ obtenu en remplaant $G $ par $M_\aa$. \\
D'aprs [W] IV.3 (4), si $P, P', P''$ sont des
sous-groupes paraboliques semi-standard de
$G$ de sous-groupe de LŽvi $M$, on l'ŽgalitŽ de fonctions rationnelles sur  $\O$:
\ber  \label{aaa}$$A(P'', P', \si) A(P', P, \si)= j(P'', P', P, \si)
  A(P'',P, \si),$$
o $j(P'', P', P, \si)$ est le produit   
des $j_\alpha(\si)$ pour $\alpha \in  \Sigma_{red}(P) \cap
\Sigma_{red}(P'') \cap \Sigma_{red}(P'^-). $\eer
On a aussi:  
\ber \label{poleA} Pour $\aa \in \Sigma_{red}(P)$, les points o  l'application  rationnelle sur  $X(M)$, $\chi \mapsto  j_\aa(\si_\chi)$,  a un p\^ole ou un zŽro sont de la forme $\chi= \chi_ \lambda$  avec $ \lambda$ ŽlŽment d'un nombre fini d'hyperplans de $Ê(\a'_M)_\C$ de la forme $\langle  \lambda, \check{\aa}\rangle   = c.$ \\ Les points o   l'application  rationnelle  sur $X(M)$, $\chi \mapsto A(P', P, \si_{\chi})$  a un p\^ole ou bien o $ A(P', P, \si_{\chi})$ n'est pas inversible  sont de la forme $\chi= \chi_ \lambda$  avec $ \lambda$ ŽlŽment d'un nombre fini d'hyperplans de $Ê(\a'_M)_\C$ de la forme $\langle  \lambda, \check{\aa}\rangle  = c$,  
 avec $\aa \in \Sigma(P') \cap \Sigma(P^-)$(cf. [H],   p. 393).\eer
Par transport de structure, on a, pour $x\in G$, normalisant $M_0$:
\beq \label{xA} \lambda(x) A(P', P, \si) = A(x.P', x.P, x\si) \lambda(x) \eeq 
et 
\ber \label{jaawsi}  Si $\aa$ est une racine rŽduite de $\Sigma(P)$, $w \in W^G$, alors:
$$   j_\aa (\si)= j_{w\aa }(w \si).$$ \eer 
\subsection{Matrices $B$}
\begin{prop}\label{B} (i) Soit  $\O$ l'orbite inertielle d'une reprŽsentation lisse  irrŽductible de $M$. Il existe une unique application rationnelle  dŽfinie sur $\O$, $
B(P,P',.)$  ˆ valeurs dans $Hom_\C(Wh(P',.), Wh(P,. ))$  telle que l'on ait
l'ŽgalitŽ de fonctions rationnelles sur $\O$: 
$$\xi(P', \si,  \eta) \circ  A(P', P,\si)= \xi(P, \si, B(P,P', \si)\eta), \> \eta \in Wh(P', \si).$$
(ii) La rationalitŽ a ici le sens suivant:\\ 
Soit $\si$ un objet de $\O$. Pour tout $\chi\in X(M)$, $Wh(P,\si_\chi)= Wh(P, \si)$,  et,
avec les notations de (\ref{intertw}), pour tout $\eta\in Wh(P', \si)$, la fonction $\chi \mapsto b(\chi) B(P,P', \si_\chi)\eta$ est une fonction polynomiale sur $X(M)$ ˆ valeurs dans $Wh(P, \si)$.\\
 La relation qui dŽfinit $B(P',P, \si)$ comme fonction sur $\O$ est la suivante:  
\\ Soit $(\si,E)$, $(\si_1, E_1)$ deux objets de $\O$ Žquivalents  et $T$ un entrelacement bijectif entre $\si$ et $\si_1$. La transposŽe de $T$, $T^t$, dŽtermine une bijection entre $Wh(P, \si_1)$ et $Wh(P, \si)$ d'une part,  $Wh(P', \si_1)$ et $Wh(P', \si)$ d'autre part et l'on a:
$$B(P,P', \si_1) = (T^t)^{-1} B(P, P', \si) T^t.$$
(iii) La fonction  rationnelle sur $X(M)$, $\chi \mapsto B(P,P', \si_\chi)$ n'a de p\^oles qu'en des points ou l'application $\chi \mapsto A(P', P,  \si_\chi)$  a un  p\^ole. 
\end{prop}
\dem
Soit $(\si, E)$ un objet de $\O$. Soit $b$ comme dans (\ref{ratintertw}). D'aprs le ThŽorme \ref{jacquetint} (iv), on a:
\ber \label{etasi} Pour tout
$\eta \in Wh(P', \si)$, il existe un unique $\eta(\chi) \in Wh(P, \si)$  tel que: 
  $$\xi(P', \si_ \chi, \eta) \circ (b(\chi) A(P', P,\si_ \chi))
=
\xi(P,
\si _
\chi, \eta(\chi) ). $$\eer Montrons que l'application 
$\chi \to \eta(\chi)$ est polynomiale en $\chi\in X(M)$. Soit $w$ comme en (\ref{wP}),  de sorte que $w.P$ est anti-standard  et $Wh(P, \si)$ est Žgal ˆ $Wh(w\si)$.  Il  s'agit de montrer que:
\ber
\label{toute}  Pour tout $e\in w.E=E$, $\langle \eta(\chi), e\rangle   $ dŽpend polynomialement de $\chi\in X(M)$\eer Avec les notations de (\ref{vsi}) et (\ref{xiv}), o on change
$\si $ en
$w\si $, on a: $$\langle \eta (\chi),e\rangle   = vol(H_{w. U^-})^{-1}\langle  \xi(w.P, w\si_\chi,  \eta(\chi)), v^{w.P, H}_{w\si_\chi,e}\rangle  .$$ Mais  la DŽfinition \ref{whp} montre que:   $$ \xi(P, \si_\chi,\eta(\chi) )= \xi(w.P, w\si_\chi, \eta(\chi))\circ \l (w).$$ 
Donc $$\langle \eta(\chi),e\rangle  = vol(H_{w. U^-})^{-1}\langle \xi(P, \si_\chi,
\eta(\chi)),  \l (w^{-1}) v^{w.P, H}_{w\si_\chi,e}\rangle  .$$ Comme $w\in K$, on dŽduit de (\ref{vsik}) 
  que la restriction ˆ $K$ de  $\l (w^{-1})  v^{w.P, H}_{w\si_\chi,e}$ est
indŽpendante de  $\chi \in X(M)$.  Alors (\ref{etasi}) permet d'exprimer $\langle \eta(\chi),e\rangle  $ ˆ l'aide de $\xi(P', \si_\chi, \eta)$.  L'assertion (\ref{toute} ) rŽsulte   de (\ref{ratintertw}) et des propriŽtŽs des fonctionnelles de Jacquet (cf. ThŽorme \ref{jacquetint} (v)). \ste
On vŽrifie, gr\^ace ˆ la dŽfinition de $Wh(P, \si_1)$ et $Wh(P, \si)$, que $T^t$ dŽtermine bien un isomorphisme entre ces deux espaces et de m\^eme pour $P'$.\\
On pose $B(P',P, \si_\chi)\eta:= b(\chi)^{-1} \eta(\chi) $. On voit gr\^ace  ˆ l'unicitŽ dans (\ref{etasi})  et ˆ (\ref{toute}) que cela dŽfinit bien une fonction rationnelle sur $\O$, en utilisant par exemple (\ref{fonct'}), et  qui a toutes  les propriŽtŽs voulues. Ceci prouve (ii).
\\Prouvons (iii).  Si $A(P', P, \si_\chi)$ n'a pas de p\^ole en $\chi_0$,  pour tout $v \in i^K_{K\cap P} $, l'application $\langle  \xi( P, \si_\chi , B(P, P',  \si_\chi)\eta, v_\chi \rangle  $ est Žgalement sans pole en $\chi_0$. Mais il rŽsulte de la  DŽfinition \ref{whp} et (\ref{xiv}), que ceci suffit ˆ assurer que $B(P,P', \si_\chi) \eta$ n'a pas de p\^ole en $\chi_0$.   \qed
\subsection{Induction par Žtage pour les matrices $B$}
\begin{prop} \label{indB}
Soit $P=MU$,  $P'=MU'$  deux  sous-groupes paraboliques semi-standard de $G$ de
m\^eme sous-groupe de LŽvi,  contenus dans un mme sous-groupe parabolique anti-standard de $G$, $P_1=M_1U_1$.
Soit  $\O$ l'orbite inertielle d'une reprŽsentation lisse  irrŽductible de $M$ et $(\si, E)$ un objet de $\O$.   Alors: \\
(i)  $Wh(M_1\cap P, \si) = Wh(P, \si)$, $Wh(M_1\cap P', \si) =
Wh(P',
\si)$. \ste (ii)  On  a l'ŽgalitŽ de fonctions rationnelles sur $\O$: 
$$B(P, P', \si) = B(P\cap M_1, P' \cap M_1, \si).$$
\end{prop}
\dem (i) 
On utilise les notations de (\ref{wP}), qu'on utilise aussi  pour $M_1$. Comme $W^{M_1}$ est contenu dans $W^G$ et que $P_1$ est anti-standard,  il rŽsulte des dŽfinitions que:  
$$w_{P}= w_{P\cap M_1}, w_{P'}= w_{P'\cap M_1}
$$  Joint ˆ (\ref{wP}), ceci
prouve la premire ŽgalitŽ de (i).   On prouve la deuxime ŽgalitŽ  de manire identique. Ceci prouve (i).  \\ Prouvons (ii). 
Par rationalitŽ, il suffit de prouver l'ŽgalitŽ lorsque $T= A(P'\cap M_1, P\cap M_1, \si)$ est bijectif. Soit $\eta' \in Wh(P', \si)$ et
calculons $\xi= \xi(P, \si,B(P,P', \si) \eta')$. On a, par dŽfinition des matrices $B$: 
$$\xi= \xi(P', \si, \eta') \circ A(P', P, \si).$$
On identifie, gr\^ace ˆ l'induction par Žtages,  $i^G_P\si $ avec $i^G_{P_1} \si^-$ et $i^G_{P'}\si $ avec $i^G_{P_1} \si_1^-$, o $\si^-= i^{M_1}_{P\cap M_1}\si$ et  $\si_1^-= i^{M_1}_{P'\cap M_1}\si$. Alors, d'aprs [W] IV. 1 (14): 
  $$A(P', P, \si)= ind T. $$
On utilise la DŽfinition \ref{whp}, puis on applique le  Lemme \ref{indetxi} ˆ $w_{P'}.P$ et $w_{P'} \si$, et ˆ nouveau la DŽfinition \ref{whp} pour $M_1$   et le fait que $w_{P'} \in M_1$ pour voir que:   \ber \label{xiPxiP1}$$\xi(P', \si, \eta')= \xi (P_1, \si_1^-, \xi(P'\cap M_1, \si, \eta)).$$\eer
  Joint ˆ ce qui prŽcde et ˆ  (\ref{fonctxi}), on en dŽduit: 
 $$\xi=\xi(P_1, \si^-, T^t \xi(P'\cap M_1, \si, \eta')).$$
Mais par dŽfinition de $T$ et des matrices $B$, on a:
$$T^t \xi(P'\cap M_1, \si, \eta')= \xi(P\cap M_1, B(P\cap M_1, P'\cap M_1, \si) \eta ').$$
Finalement on a prouvŽ:
$$\xi(P, \si,B(P,P', \si) \eta')= \xi(P, \si,B(P\cap M_,P'\cap M_1, \si) \eta').$$
D'o l'on dŽduit (ii). \qed 
 
\subsection{Equation  fonctionnelle des intŽgrales de Jacquet}
\begin{lem} \label{eqfonct}
(i) Avec les notations de la DŽfinition  \ref{eis}, soit $(\si, E) $ un objet de $\O$. Soit $P=MU$, $P'=MU'$ deux sous-groupes paraboliques semi-standard de $G$ de sous-groupe de LŽvi $M$. On a l'ŽgalitŽ de fonctions rationnelles  sur $X(M)$:
$$E^G_P(
\si, B(P, P', \si_\chi)\eta, v_\chi)(g) =E^G_{P'} (\si, \eta, A(P', P, \si_\chi )  v_\chi)(g),  v \in i^K_{P\cap K}E, \eta \in Wh(P', \si).$$
(ii) Avec les notations de (\ref{wP}), on a l'ŽgalitŽ:  
 $$ E^G_P(\si, \eta, v)=E^G_{w_P.P}(w_P\si, \eta, \l (w_P) v ), \> v \in i^G_P E, \eta \in Wh(P, \si).$$
(iii) On suppose que $P$ est anti-standard. Si $(\si_1, E_1)$ est une reprŽsentation de $M$  Žquivalente ˆ $(\si, E) $  et,  si  $T$ est  un opŽrateur d'entrelacement bijectif entre $\si$ et $\si_1$, on a,  avec les notations de la Proposition  \ref{xisisi1}: 
$$E^G_{P} (\si, T^t\eta_1,  v) =E^G_{P} (\si_1, \eta_1, (ind  T) v),  v\in i^K_{P\cap K}E,
\eta \in Wh(\si_1).$$
\end{lem}
\dem
(i) est une consŽquence immŽdiate de la dŽfinition des intŽgrales de Jacquet (DŽfinition \ref{eis}) et de celle des matrices $B$ (Proposition \ref{B}).
\ste (ii) rŽsulte de (\ref{wP}) et de la dŽfinition des intŽgrales de Jacquet.\\
(iii) rŽsulte de la Proposition \ref{xisisi1}. 
 \qed
\section{EnoncŽ du ThŽorme principal}
\setcounter{equation}{0}
\subsection{Formes sesquilinŽaires  }
{\bf Hypothse supplŽmentaire}
{\em On suppose dŽsormais que $\psi$ est en outre unitaire.}
\ste
Si $E$ est un espace vectoriel complexe, on note ${\overline E}$ l'espace vectoriel
 conjuguŽ: c'est le mme groupe additif, mais la multiplication par les scalaires est conjuguŽe. Soit $(\pi,
V)$  une reprŽsentation lisse de $G$. On rappelle qu'on note $({\check \pi},
{\check  V})$ sa contragrŽdiente lisse. On note $({\overline \pi},
{\overline V})$ la reprŽsentation conjuguŽe et $(\pi^*, V^*)$ la
reprŽsentation $({\overline {\check \pi}}, {\overline {\check V}})$. Notez
que ${\overline {\check V}}$ s'identifie naturellement  ˆ l'espace des
formes antilinŽaires sur $V$ fixŽes par un sous-groupe compact ouvert de
$G$.\\
  Dans la suite produit scalaire voudra dire produit scalaire linŽaire dans la premire variable et antilinŽaire dans la seconde.
Si $\pi$ est unitaire, i.e. muni d'un produit scalaire invariant, $V^*$ s'identifie naturellement ˆ $V$, par l'application $ v\mapsto (v,.)$  et $\pi^*$ ˆ $\pi$. 
Si $\chi$ est un ŽlŽment de $X(G)$ on note $\chi^{-1}$ son inverse et ${\overline \chi}$ son complexe
conjuguŽ. Alors $(\pi\otimes \chi)^*$ est naturellement isomorphe ˆ $\pi^* \otimes {\overline \chi}^{-1}$.
Soit 
$\O$  l'orbite inertielle d'une reprŽsentation unitaire irrŽductible lisse de
$G$.  Si
$(\pi, V)$ est un objet de
$\O$, on dŽduit de ce qui prŽcde que
$(\pi^{*}, V^*)$ est aussi un objet de $\O$.
Montrons que: 
\ber Dans l'orbite inertielle $\O$ d'une reprŽsentation cuspidale lisse irrŽductible, il existe une reprŽsentation unitaire. \eer
En tensorisant par un caractre non ramifiŽ de $G$, on  trouve, dans l'orbite inertielle, une reprŽsentation,  $(\pi, V)$, telle que la restriction ˆ $A_G$ de $\pi$ soit donnŽe par un caractre unitaire.  On fixe $\check{v_0}\in \check{V}$ non nul et on dŽfinit:
$$(v,v')= \int_{A_G \backslash G}c_{\check{v_0},v}(g) \overline{ c_{\check{v_0},v'}(g)} dg, v, v' \in V.$$
On voit facilement que c'est un produit scalaire $G$-invariant, donc que $\si$ est unitaire.

\begin{lem} \label{prodscal}
On suppose que $(\pi, V)$ est une reprŽsentation  lisse, cuspidale, unitaire  et irrŽductible de $G$. Alors:\\
(i) Il existe un unique produit scalaire hermitien sur $Wh(\pi)$ tel que:
$$\int_{A_G U_0\backslash G}c_{\xi,v}(g) \overline{ c_{\xi',v'}(g)} dg= (\xi, \xi') (v, v'), \> \xi, \xi' \in Wh(\pi), v, v' \in V.$$
o la fonction sous le signe intŽgrale est ˆ support compact d'aprs  le Lemme \ref{borne} (i).\\
(ii) Si $\chi$ est un caractre unitaire  non ramifiŽ de $G$, le produit scalaire sur $Wh(\pi_ \chi)=Wh(\pi)$, ne dŽpend pas de
$\chi$.\\
(iii) Si $T$ est un opŽrateur d'entrelacement unitaire avec une autre reprŽsentation cuspidale de $G$, $(\pi_1, V_1)$, l'opŽrateur $T^t$ dŽtermine un opŽrateur unitaire entre $Wh(\pi)$ et $Wh(\pi_1)$. 
\end{lem}
\dem
(i) Il s'agit d'une simple application du Lemme de Schur. \ste
(ii) rŽsulte immŽdiatement de la caractŽrisation du produit scalaire. \\
(iii) est immŽdiat. \qed
 \ste On appliquera  ces notations aux sous-groupes de LŽvi de $G$. 
\subsection{TransformŽe de Fourier-Whittaker}
\begin{prop} \label{fpo}
(i) Soit $P=MU$ un sous-groupe parabolique  semi-standard de $G$ et  $\O$ l'orbite inertielle d'une reprŽsentation lisse unitaire irrŽductible et cuspidale  de
$M$. Soit $(\si, E)$ un objet de $\O_u$.  On munit  $Wh(P, \si)$ du produit scalaire dŽfini  gr\^ace au Lemme \ref{prodscal}  et ˆ la dŽfinition de $Wh(P, \si)$.  Par tensorisation avec le produit scalaire d'induite unitaire de $i^G_P(\si)$, on en dŽduit un produit scalaire sur $Wh(P, \si)\otimes i^G_PE $. \\ Soit $f \in \ccu$. Il existe  un unique ŽlŽment de $Wh(P, \si)\otimes i^G_PE $, $F(P, \si
)$,  tel que: 
\beq \label{fpsi}(F(P, \si), \eta \otimes v)= \int_{U_0 \backslash G} f(g) \overline{ E^G_P(\si, \eta, v)(g) } dg,
\> v \in   i^G_PE  , \eta \in Wh(P, \si).\eeq
(ii) Utilisant l'ŽgalitŽ $Wh(P, \si_ \chi)= Wh(P,  \si) $  pour tout $\chi \in X(M)$ et la rŽalisation compacte, on voit que $\chi \mapsto F(P, \si_\chi)$ s'Žtend de $X(M)_u$ ˆ $X(M)$ en une fonction  polynomiale notŽe de m\^eme. 
\\ (iii) Si  $\si_1$ est une reprŽsentation unitaire Žquivalente ˆ $\si$ de $M$ et $T$ un opŽrateur d'entrelacement unitaire entre $\si$ et $\si_1$. On utilise  les notations de la Proposition \ref{B} (ii)  et l'unitaritŽ de l'opŽrateur $ind T$.  Alors on a: 
\ber \label{fpsi1}  $$F(P, \si_1) =((T^t)^{-1} \otimes ind T)  F(P, \si ).$$ \eer
Utilisant, les notations de (\ref{polo}),   pour toute orbite inertielle d'une  reprŽsentation lisse  unitaire irrŽductible et cuspidale  de $M$, $\O$,  $\si \mapsto F(P, \si )$ est un ŽlŽment de $Pol(\O_u,Wh(P,.)\otimes  i^G_P) $. On fait agir $G$ sur ce dernier espace par une reprŽsentation notŽe $\rho_\bullet$ dŽfinie par:
$$ (\rho_ \bullet (g)F)(P, \si)= (Id \otimes i^G_P\si(g)) F(P, \si).$$

\end{prop}
\dem
Notons $\phi(v)$ le second membre de  (\ref{fpsi}). Alors $\phi$ est une forme antilinŽaire sur $i^G_PE$. Il est clair que si $ f$ est invariante ˆ droite  par un sous-groupe ouvert compact $H$ de $G$,  $\phi$ est invariante par $H$. Par le thŽorme de reprŽsentation de Riesz, en dimension finie, on en dŽduit (i). 
\\ Prouvons (ii). Il suffit de voir que pour tout $v \in   i^K_{P\cap K}E  , \eta \in Wh(P, \si)$, l'application  $\chi \mapsto (F(P, \si_\chi), \eta\otimes v_\chi)$ s'Žtend de faon polynomiale de $X(M)_u$ ˆ $X(M)$. 
Pour cela il suffit de montrer que l'application qui ˆ $\chi\in X(M)$ associe:
$$\phi_\chi(v):=\int_{U_0 \backslash G} f(g) \overline{ E^G_P((\si_\chi)^*, \eta, v_\chi)(g) } dg,
$$
est polynomiale. Choisissons un sous-groupe compact ouvert, $H$, de $G$ comme ci-dessus et   fixons  $v$. En utilisant le fait que $f$ est ˆ support compact, on voit qu'il existe des constantes $c_1, \dots, c_n \in \C$  et  $g_1, \dots, g_n \in G$ telles que:
$$ \phi_\chi(v)= \sum_{i=1, \dots, n} c_i f(g_i)\overline{ E^G_P((\si_\chi)^*, \eta, v_\chi)(g_i) }, \> \chi \in X(M).$$
Notre assertion rŽsulte alors de (\ref{eispol}).  Cela prouve (ii). La relation 
(\ref{fpsi1})  rŽsulte  de la dŽfinition de $F$ en (i), de la dŽfinition des intŽgrales de Jacquet (DŽfinition \ref{eis}), et de (\ref{fonctxip}). Le reste de (iii) rŽsulte alors de (ii). \qed
On retient les notations du  Lemme \ref{eqfonct} (ii). L'unitaritŽ de $\l(w_P)$ montre que la dŽfinition de $F$ implique que,  pour $\si$ comme dans (i):
\ber \label{fpw} $$ F(P, \si)= (Id \otimes \l(w_P)^{-1}) F(w_P.P, w_P\si).$$ \eer 
  De mme, d'aprs le Lemme \ref{eqfonct} (i),  on l'ŽgalitŽ de fonctions rationnelles sur $X(M)_u$:
\beq \label{BF}(B(P, P', \si_\chi) ^* \otimes Id )F(P, \si_\chi) = (Id\otimes A(P', P, \si_\chi)^* )F(P', \si_\chi).  \eeq
On sait [W], preuve du Lemme V.2.2,  que l'on a la formule d'adjonction, pour $\si$ unitaire:
\beq A(P', P,  \si)^*= A(P, P',  \si)  .\eeq
{\bf On admet provisoirement la relation suivante  pour $P$ anti-standard et $\si $ unitaire}
\beq \label{B*} B(P, P', \si )^*= B(P', P,  \si) .\eeq 
Nous montrerons   cette relation au prix d'un travail non nŽgligeable, 
comme consŽquence de l'Žtude des produits scalaires de paquets d'ondes.  On dŽduit alors de la relation(\ref{BF}) que:
\ber  \label{fpab} Si $P=MU$ est anti-standard, on a l'identitŽ de fonctions rationnelles sur $\O_u$:
$$(B(P', P, \si)  \otimes Id )F(P, \si) = (Id\otimes A(P, P', \si) )F(P', \si).$$  
  \eer 
{\bf Notation} On notera ${\hat f}(P, \si) $ au lieu de $F(P,\si)$ et on appellera $\hat{f}$ la transformŽe de Fourier-Whittaker de $f$, ou plus simplement sa transformŽe de Fourier. 
\subsection{EnoncŽ du thŽorme principal}
\begin{theo}\label{PW}
On suppose donnŽ pour tout  sous-groupe parabolique semi-standard  de $G$, $P=MU$, et pour toute reprŽsentation lisse unitaire, irrŽductible et cuspidale de $M$,  $(\si, E)$, un ŽlŽment  $F(P,\si)$ de $ Wh(P,\si) \otimes  i^G_P E $. 
 Alors il existe $f \in \ccu$ tel que ${\hat f}= F $ si et seulement si:\\
(i) Pour toute  orbite inertielle de reprŽsentation lisse unitaire irrŽductible et cuspidale de $M$, $\O$,  $F\in Pol(\O_u, Wh(P,.) \otimes i^G_P)$.   Notamment $F $ vŽrifie (\ref{fpsi1}). \\
(ii) Faisant agir $G$ sur $Wh(P,\si) \otimes i^G_P \si$ trivialement sur le premier facteur et naturellement sur le deuxime, il existe un sous-groupe compact ouvert fixant $F(P, \si)$ pour tout $(P, \si)$ comme ci-dessus.\\
(iii) Pour tout sous-groupe parabolique semi-standard $P$ de $G$, on a, avec les notations de (\ref{wP}):
$$ F(P, \si)= (Id \otimes \l(w_P)^{-1}) F(w_P.P, w_P\si).$$
(iv) Si $P=MU, P'=MU'$ sont deux sous-groupes paraboliques anti-standard de $G$ de m\^eme sous-groupe de LŽvi   et si $P$ est anti-standard, on a l'identitŽ de fonctions rationnelles sur $\O_u$:
$$(B(P', P, \si)  \otimes Id )F(P, \si) = (Id\otimes A(P, P', \si) )F(P', \si).$$  
  De plus $f$ est unique. \\En d'autres termes la transformation de Fourier-Whittaker dŽtermine un isomorphisme entre $\ccu$ et l'espace des fonctions $F$ satisfaisant les conditions (i) ˆ (iii) ci-dessus.
\end{theo}
Remarquons que modulo la formule d'adjonction (\ref{B*}), on a montrŽ la partie seulement si du ThŽorme.
\subsection{Reformulation du ThŽorme principal }
Soit $P=MU$, $Q=LV$ deux sous-groupes paraboliques standard et soit $s\in W^G$ tel que $s.M'=M$ et tel que $s$ soit de longueur minimum dans $s{\overline W}^{L}= {\overline W}^{M}s$. Soit $w=s^{-1}$. Alors $P':= w^{-1}. Q$ est un sous-groupe parabolique semi-standard de sous-groupe de Levi $M$. Notons, pour $\si$ objet d'une orbite inertielle, ${\cal O}$,  de  reprŽsentation irrŽductible lisse  cuspidale de $M$,  $A(w, P, \si):= \lambda(w)\circ A(P', P, \si)$, qui entrelace, lorsqu'il est dŽfini, $i_P^G \si$ et $i_Q^G w\si$. Comme dans la Proposition \ref{B}, on voit  qu'il  existe une unique fonction rationnelle sur $\O$ ˆ valeurs dans $End(Wh(w\si), Wh(\si))$,  $B(w^{-1}, Q, \si)$, telle que l'on ait l'identitŽ de fonctions rationnelles sur $\O$:  
\beq\label{xiwb} \xi (Q, w\si, \eta) \circ A(w , P,\si) =  \xi (P, \si, B(w^{-1}, Q, w\si)\eta), \eta \in Wh(\si)\eeq 
On remarque qu'avec nos dŽfinitions, $w=w_{P'}$. En utilisant (\ref{xiwb}), la DŽfinition \ref{whp} et la Proposition \ref{B}, on voit que \beq 
\label{BwBP} B(w^{-1}, Q, w\si)= B(P, P', \si) \eeq 
Si $f \in \ccu$ on note $\hat{f}_{anti}$  la restriction de $\hat{f}$ aux sous-groupes paraboliques anti-standard. \\ 
{\bf Corollaire du ThŽorme \ref{PW}}\\ {\em 
On suppose donnŽ pour tout  sous-groupe parabolique anti-standard  de $G$, $P=MU$, et pour toute reprŽsentation lisse unitaire, irrŽductible et cuspidale de $M$,  $(\si, E)$, un ŽlŽment  $F(P,\si)$ de $ Wh(P,\si) \otimes  i^G_P E $. 
 Alors il existe $f \in \ccu$ tel que ${\hat f}_{anti} = F $ si et seulement si $F$ vŽrifie (i) et (ii)   du ThŽorme \ref{PW}  et  si $F$ vŽrifie la condition suivante:\\ Pour tout $P=MU$, $Q= LV$ sous-groupes paraboliques anti-standard de $G$, pour tout  $s\in W^G$ tel que $s.M'=M$ et tel que $s$ soit de longueur minimum dans $s{\overline W}^{L}= {\overline W}^{M}s$, et   pour toute orbite inertielle de  reprŽsentation irrŽductible lisse  cuspidale de $M$, on a l'identitŽ de fonctions rationnelles sur $\O$:\beq \label{bfa}(B(w^{-1}, Q, w\si) \otimes Id)F(Q, w\si) = (Id\otimes A(w, P, \si)) F(\si),\> avec  \> w=s^{-1} \eeq 
 De plus $f$ est alors unique.}\\
 \dem  On Žtend  $F$ aux sous-groupes paraboliques anti-standard en utilisant la relation (ii) du ThŽorme \ref{PW} comme dŽfinition. La fonction ${\tilde F}$ ainsi obtenue vŽrifie toutes les conditions du ThŽorme \ref{PW}: la condition (iv) rŽsulte immŽdiatement de la relation ci-dessus satisfaite par $F$, de (\ref{BwBP}) et  (\ref{bfa}) et de la dŽfinition de ${\tilde F}$.  On en dŽduit l'existence de $f$. Par ailleurs si $f'\in \ccu$ et si  ${\hat f'}_{anti} $ est nulle, $\hat f'$ est nulle. Donc, d'aprs le ThŽorme \ref{PW}, $f'$ est nulle.  On en dŽduit que $f$ est unique. \qed
\section{Terme constant des intŽgrales de Jacquet}
\setcounter{equation}{0}
\subsection{Fonctionnelles de Jacquet et modules de Jacquet}
Soit $P= MU$ et $P'= M'U'$ deux sous-groupes paraboliques semi-standard de $G$. 
  On note $$\mathcal{W}(M'\vert G\vert M)= \{ s \in W^G \vert s.M\subset M'\}.$$
 On surligne pour indiquer l'image dans le groupe de Weyl.  Soit $$ {\overline  W} (M'\vert G\vert M)={\overline  W}^{M'} \backslash {\overline  {\cal W}} (M'\vert G\vert M). $$
On remarque que $$ {\overline  W} (M'\vert G\vert M):= {\overline  W}^{M'} \backslash {\overline  {\cal W}} (M'\vert G\vert M)/ {\overline  W}^{M}. $$
Ceci permet de choisir un ensemble de reprŽsentants $ W (M'\vert G\vert M)$ de  ${\overline  W} (M'\vert G\vert M)$   dans $W^{G}$
tel que: \ber \label{wm'gm} Pour $s\in W (M'\vert G\vert M)$,   ${\overline s} $ est de longueur minimum dans ${\overline W}^{M'}{\overline s}={\overline s} \overline{W}^M$. \eer
De plus, en utilisant un sous-groupe parabolique standard conjuguŽ ˆ $P$, $x.P$, on voit gr\^ace  ˆ [War],  Proposition 1.2.1.10, que \ber \label{P'sP}
$W(M'\vert G\vert M)$ est un ensemble de reprŽsentants de $P' \backslash G/ P$ \eer  
 \begin{defi} \label{reg} Soit $M$ le sous-groupe de LŽvi d'un sous-groupe parabolique semi-standard de $G$ et $(\si, E)$  une reprŽsentation lisse,  irrŽductible et cuspidale de $M$. On dit que $\si$ est $G$-rŽgulire si:
 \\1) Pour tout $s\in W(M \vert G\vert M)$  avec $s \not=1$, les reprŽsentations  $ s\si$ et $\si$ sont  non Žquivalentes.
 \\2) Si $P, P'$ sont des sous-groupes paraboliques de $G$ ayant $M$ pour sous groupe de LŽvi, les applications rationnelles sur $X(M)$, $\chi \mapsto A(P, P', \si_\chi)$, $\chi \mapsto B(P, P', \si_\chi)$ n'ont pas pas de poles en $\chi=1$ et leurs valeurs en $1$ sont des opŽrateurs inversibles. De plus $i^G_P \si$ est irrŽductible. 
 \\ On remarque que  si $\si$ est comme ci-dessus, l'ensemble des $\chi \in X(M)$ tel que $\si_\chi$ soit $G$-rŽgulire est un ouvert de Zariski non vide de $X(M)$ et  si $\si$ est unitaire, l'ensemble des $\chi\in X(M)_u$  tels que $\si_\chi$ soit $G$-rŽgulire est   Zariski-dense dans $X(M)$. 
 \end{defi} 
On introduit, pour $s \in W (M'\vert G\vert M)$,  les sous-groupes paraboliques de $G$:
\beq \label{ps} P_s= (M'\cap s.P)U', \> {\tilde P}_s= (M'\cap s.P) U'^-.\eeq
Soit $\si$ une reprŽsentation cuspidale $G$-rŽgulire de $M$. Posons $(\pi, V)= (i^G_P \si,
i^G_PE)$. On dŽfinit une application $\alpha$:
$$ \alpha: V \to \oplus_{s \in W (M'\vert G\vert M)} i^{M'}_{M'\cap s.P} sE,  \> v \mapsto (v_s)_{s \in
W(M'\vert G\vert M)}.$$
par $$v_s(m')= \delta^{-1/2} _{P'}(m') (A(P_s, s.P, s\si) \l (s)v)(m'), \> m'\in M'.$$  
D'aprs [W], dŽbut  de la preuve de la Proposition V.1.1, on a:
\ber  \label{alpha} L'application $\alpha$ se factorise en un isomorphisme de $M'$-modules entre le module de acquet normalisŽ de $V$ relatif ˆ $P'$,
$V_{P'}$,  et  l'espace d'arrivŽe, qui est une somme de $M'$-modules irrŽductibles non
Žquivalents, rŽduite ˆ zŽro si $W(M'\vert G\vert M)$ est vide. On identifie dans la suite $V_{P'} $ ˆ l'aide de $\alpha$ avec l'image cet isomorphisme.  \eer 
\ber \label{seta} Si $\eta \in Wh(P, \si)$, comme $\l(s)$ entrelace $i^G_P \si$ et $i^G_{s.P} s\si$,  d'aprs (\ref{bijetaxi}), il existe un unique $\s \eta \in Wh(s.P, s\si)$ tel que:
$$ \xi(s.P, s\si, \s \eta)= \xi(P, \si, \eta) \circ \l (s^{-1}).$$
\eer 
\begin{theo} \label{thct} Soit $P$ (resp.  $P'$)  un  sous-groupe parabolique semi-standard (resp. standard) de $G$. Soit $\si$ comme ci-dessus. \ste
(i) Pour $s \in W (M'\vert G\vert M)$, on a $Wh({\tilde  P}_s,s\si) = Wh (M'\cap s.P, s\si)$. \ste 
(ii)Soit $\eta\in Wh(P, \si)$.  On note $\xi= \xi(P, \si, \eta)$.  Alors,
dans l'isomorphisme ci-dessus,  
$\xi_{P'} $ qui est un ŽlŽment du dual de $V_{P'}$ (cf. (\ref{jpxi})), est    nul si $W(M'\vert G\vert M)$ est vide et  sinon Žgal ˆ $(\xi_s)_{s \in W(M'\vert G\vert M)}$, avec: 
$$\xi_s = \xi(M'\cap s.P, s\si, B({\tilde P}_s, s.P, s\si)\s \eta), $$
l'expression Žtant bien dŽfinie gr\^ace ˆ (i).
\\(iii) Si $P$ est anti-standard, $\s \eta=\eta$. 
\end{theo}
\dem
(i)  On rappelle  que, par dŽfinition de $W^{M'}$, 
$W^{M'}= W^G\cap M'$.    
Si $w^{-1}\in W^{M'}$ et  $w.(M'\cap s.P)$ est anti-standard  dans $M'$, $w. {\tilde  P}_s  $ est anti-standard
dans
$G$.  Si $w^{-1}$ est de longueur minimum dans $W^{s.M}
w^{-1}$, on  a
$w=w_{M'\cap s.P}= w_{{\tilde  P}_s}$. (i) en rŽsulte d'aprs la DŽfinition 
\ref{whp}.\ste(ii)  Si  $W(M'\vert G\vert M)$ est vide, $V_{P'}$ est  rŽduit ˆ zŽro, donc $\xi_{P'}$ est nul. Ceci prouve la premire assertion de (ii). \ste On suppose maintenant que $W(M'\vert G\vert M)$ est non  vide.  On Žcrit: 
$$ \xi_{P'}=(\xi_s)_{s \in W(M'\vert G\vert M)}, $$
o $\xi_s$ est de la forme:
\ber \label{defetas}$\xi_s= \xi(M'\cap s.P, s\si, \eta_s)$,  pour un ŽlŽment $\eta_s$ de $Wh(M'\cap s.P, s \si)$.
\eer 
ll s'agit donc de montrer:
\ber \label{etas} $$\eta_s= B({\tilde P}_s, s.P, s\si) \s \eta .$$ \eer
a) Montrons d'abord:
\ber\label{eta1} Supposons $P$ semi-standard et  $P^-\subset P'$. Alors $\eta_{1_G}= \eta . $\eer
Soit $ (\si^-, E^-)=(i_{P\cap M'}^{M'} \si, i_{P\cap M'}^{M'}E)$ de sorte
que $\pi= \igps$ s'identifie ˆ $i_{P'^{-}}^G \si^-$. On note $\eta^{-}=\xi(P\cap
M', \si, \eta)$. Comme dans la preuve de (\ref{xiPxiP1}), on voit que $\xi$ est Žgal ˆ $\xi(P'^-,
\si^-, \eta^-)$.    Soit $v \in V$ identifiŽ $i_{P'^-}GE^{-}$, ˆ support dans $P'^-U'$.  Soit 
$a\in A_{M'}$. On calcule $\langle \xi, \pi(a) v\rangle  $  en
utilisant le ThŽorme \ref{jacquetint} (ii) pour $P'$. On a:
$$\langle \xi, \pi(a) v\rangle  = \int_{U'} \langle \eta^-, v(u'a)\rangle   \psi(u')^{-1}  du'.$$En utilisant les relations de covariance satisfaites par $v$, on voit que la fonction ˆ intŽgrer est  non nulle pour $a^{-1}u'a\in Supp \>v$,  i.e.  $u'\in
a (Supp v) a^{-1}$. Pour  $\varepsilon >0$ assez petit et $a\in A_{M'} \cap A_0^-(P', < \varepsilon)$,   
$u'$ est  tel que $\psi(u')=1$ et l'on trouve: 
\beq \label{xieta_}\langle \xi, \pi(a )v\rangle  = 
\langle \eta^-,(A(P', P'^-, \si^-) v)(a) \rangle .  \eeq
Ici on a ${\tilde P}_{1_G}= P$.
Mais avec les identications de l'induction par Žtages et la dŽfinition de
${P}_{1_G}$, on a (cf.[W], IV.1(14)):
\beq \label{Asi-}A(P', P'^-, \si^-) v=A(P_{1_G}, P, \si)v .\eeq 
Alors, tenant compte de la dŽfinition de $v_{1_G}$, (\ref{xieta_}) se rŽŽcrit:
\beq \label{xiav1}\langle \xi, \pi(a) v\rangle   = \delta_{P'}^{1/2}(a) \langle  \xi(P\cap M', \si,\eta), \si^-(a)v_{1_G}\rangle .  \eeq 
Montrons: 
\ber \label{vs=0} $v_s= 0$ si $s\not =1_G $. \eer
En effet, lorsque les intŽgrales d'entrelacements sont dŽfinies par des intŽgrales convergentes, $v_s(m')$ est un multiple de l'intŽgrale sur $U' \cap s.P^-$ de $ v(s^{-1}u'm')$. Mais si $u'\in U'$ est tel que $v(s^{-1}u' m')$ est non nul, on doit avoir $s^{-1}u'm' \in PU^-\subset PP'$. Comme $s^{-1} u'm'$  appartient ˆ $ Ps^{-1} P'$ et que, d'aprs (\ref{P'sP}), $P'sP \cap P'P$ est vide, on conclut que  (\ref{vs=0}) est vrai dans ce cas. On conclut par prolongement rationnel. \\ De (\ref{vs=0}), on dŽduit que  l'image de $v$ dans $V_{P'}$ est Žgale ˆ  $v_{1_G}$. Alors, tenant compte de  (\ref{alpha}),   (\ref{xiav1}) se rŽŽcrit, pour $a$ comme ci-dessus, i.e. $a\in A_{M'} \cap A_0^-(P', < \varepsilon)$:
$$\langle \xi, \pi(a) v\rangle   = \delta_{P'}^{1/2} (a)\langle  \xi(P\cap M', \si,\eta), \pi_{P'}(a) j_{P'}(v)\rangle  . $$
Mais (cf. (\ref{jpxi})), on a, pour un $\varepsilon '>0$,  l'ŽgalitŽ:
$$\langle  \xi, \pi(a) v\rangle  = \delta_{P'}^{1/2} (a)\langle  \xi_{P'}, \pi_{P'}(a) j_{P'}(v)\rangle  ,  a \in A_{M'} \cap A_0^-({P'}, <\varepsilon' ). $$
Alors les deux membres de droite des ŽgalitŽs prŽcŽdentes sont des fonctions $ A_{M'}$-finies sur $ A_{M'}$, d'aprs les propriŽtŽs du module de Jacquet, et Žgales   sur $ A_{M'} \cap A_0(P, <inf(\varepsilon, \varepsilon'))$, donc Žgales partout. L'ŽgalitŽ en $1_G$, pour tout $v$ comme ci dessus, conduit ˆ (\ref{eta1}) en tenant compte de l'assertion suivante, que nous allons dŽmontrer.
\ber L'ensemble $\{v_{1_G} \vert v \in V, Supp\> v \subset P'^-U'\}$ est Žgal ˆ $i^{M'}_{M'\cap P} E $. \eer 
En effet, pour $H$ comme dans (\ref{vsi}) et $e \in E^-$, $v= v^{P'^-,H}_{\si^-, e}$  est ˆ support dans $P'^- U'$. La dŽfinition des intŽgrales d'entrelacement et (\ref{Asi-}) montrent  que $v_{1_G}$ est non nul en $1_G$. L'ensemble considŽrŽ est clairement un sous-$M'$-module de   $ i^{M'}_{M'\cap P} E $, qui est irrŽductible d'aprs l'hypothse de rŽgularitŽ de $\si$ (cf.  DŽfinition \ref{reg}) et qui n'est pas rŽduit ˆ zŽro. On en dŽduit l'assertion prŽcŽdente, ce qui achve de prouver  (\ref{eta1}). \\
b) Revenons ˆ la dŽmonstration de (\ref{etas}).
 On fixe  $s \in W(M'\vert G\vert M)$.
Soit $t\in W^{G}$ un reprŽsentant de l'image de $s^{-1}$ dans $\overline{W}^G$. 
  On note  $P_1$ le sous-groupe parabolique semi-standard   de $G$, ${\tilde P}_s$, de sous-groupe de  LŽvi  $M_1= t^{-1}.M$ et contenu dans $P'^-$.  Soit $(\pi_1, V_1)=
(i^G_{P_1} t^{-1}\si, i^G_{P_1}t^{-1}E)$. On affecte d'un indice $1$ tous les
objets relatifs ˆ cette reprŽsentation. Soit $v_1
\in V_1$ et $\eta\in Wh(\si)$. Posons:
  $$ v= A(P, t.P_1, \si) \l (t) v_1, \xi_1=   \xi (t.P_1,\si,   B(t. P_1, P, \si)\eta )\circ \l(t). $$
  Soit $\t B(t.P_1, P, \si) \eta$ l'ŽlŽment  de $Wh(P_1, t^{-1} \si)$  tel que ((cf. (\ref{bijetaxi})):
  \ber \label{xiPP1} $$\xi (t.P_1,\si,   B(t. P_1, P, \si)\eta )\circ \l(t)=  \xi (P_1, t^{-1} \si,   \t B(t. P_1, P, \si) \eta). $$ \eer
 On a donc: 
 \beq \label{xi1=}\xi_1 =  \xi (P_1, t^{-1} \si,   \t B(t. P_1, P, \si) \eta), 
  \eeq
  o $\t$  envoie $Wh(tP_1, \si)$ sur $Wh(P_1, t^{-1} \si)$. 
  On voit que $W(M'\vert G\vert  t^{-1}.M) $ contient $1_G$, car $s \in W(M'\vert G\vert M)$.  
  On suppose en outre que $v_{1, u}$ est nul pour tout $u \in W(M'\vert G\vert  t^{-1}. M) $  distinct de $1_{G}$.
   Alors (cf. [W], p. 291), d'aprs les propriŽtŽs des intŽgrales d'entrelacement, pour tout 
$u' 
\in W(M'\vert G\vert M)$,  $v_{u'}$
dŽpend linŽairement de $v_{1,su'}$. Donc  
$v_{u'}=0$ pour tout
$u' 
\in W(M'\vert G\vert M)$ distinct de $s$.
On a, gr\^ace a la dŽfinition des matrices $B$: 
  $$\langle \xi, \pi(g) v\rangle  = \langle \xi_1, \pi_1(g) v_1\rangle  , \> g\in G. $$ 
 D'aprs  la dŽfinition du terme constant et (\ref{tctcoef}),
  on en dŽduit:
  $$\langle \xi_{P'}, v_{P'}\rangle  = \langle \xi_{1,  P'}, v_{1,P'}\rangle  .$$
 En tenant compte de la dŽfinition de $\xi_s$ et $\xi_{1,1_G}$, on dŽduit de ce qui prŽcde:
 \beq \label{xis=xist}\langle \xi_s, v_s\rangle  = \langle \xi_{1, 1_G}, v_{1, 1_G}\rangle  . \eeq
On pose: $$m=st \in M_0\cap K.$$  
Montrons que:   \beq \label{vs=} v_s= \l(m)v_{1, 1_G}.\eeq
D'abord, on  dŽduit de (\ref{aaa}):
$$A(P_s, s.P, s\si) A(s.P, P_{1}, s\si)= A(P_s, P_{1}, s\si).$$
Tenant compte de cette ŽgalitŽ, et 
utilisant  (\ref{xA}) avec $x= s^{-1}$, puis $x=m$, on dŽduit de la dŽfinition de $v$ et de celle de $v_s$:
$$v_s(m')=\delta_{P'}^{1/2}(m')(A(P_s, P_{1}, s\si)\l(m)v_1)(m'), m' \in M'.$$
Remarquant que $P_{1, 1_G}=P_s$, la dŽfinition de $v_{1, 1_G}$ montre que: 
$$v_{1, 1_G}(m')=\delta_{P'}^{1/2}(m')(A(P_s, P_{1}, t^{-1}\si)v_1)(m'), m' \in M'.$$
Comme $m \in M_0\cap K$, ce qui prŽcde joint ˆ  (\ref{xA}) conduit ˆ (\ref{vs=}).\\ 
Etudions maintenant $\langle \xi_s, v_s\rangle  $. L'inversibilitŽ de $B({\tilde P}_s, s.P, s\si) $ (cf. DŽfinition \ref{reg}) montre que, avec les notations de (\ref{seta}), il existe un unique $\eta'\in Wh(P,\si)$  tel que:
$$\eta_s= B({\tilde P}_s, s.P, s\si, ) \s \eta' . $$ De (\ref{vs=}), de la dŽfinition de $\eta_s$  (cf. (\ref{defetas})) et de  l'ŽgalitŽ $ \tilde{P}_s= P_1$, on
dŽduit: 
 \beq\label{s=}\langle \xi_s, v_s\rangle  = \langle  \xi(M'\cap P_1, s\si, B(P_1, s.P, s\si)\s \eta'),
\l(m) v_{1,1}\rangle  .\eeq 
Par ailleurs, d'aprs (\ref{eta1}) appliquŽ ˆ $\xi_1$ et $P_1$, en tenant compte de (\ref{xi1=}),  on  voit que: 
\beq \label{xi11} \langle  \xi_{1, 1_G}, v_{1, 1_G}\rangle   = \langle \xi(M'\cap P_1, t^{-1}\si,  
 \t B(t . P_1, P, \si)\eta, v_{1, 1_G}\rangle  .\eeq
 Admettons provisoirement l'ŽgalitŽ suivante: 
 \beq \label{bm} B(P_1, s.P, s\si)\s \eta= (s\si')(m^{-1})\t B(t. P_1, P, \si )\eta.\eeq 
Alors ceci joint ˆ (\ref{xi11}), montre  que: 
 $$\langle  \xi_{1,1_G}, v_{1, 1_G}\rangle  = \langle  \xi(M'\cap P_1, t^{-1}\si, 
 (s\si')(m)B(P_1, s.P, \si)\s \eta, v_{1, 1_G}\rangle  .$$
 On utilise la deuxime relation de (\ref{ximsi}) en changeant $m$ en
$m^{-1}$, $\si$ en $s\si$ et en remarquant que
$t^{-1}= m^{-1}s$ pour trouver: 
 $$\langle  \xi_{1, 1_G}, v_{1, 1_G}\rangle  = \langle  \xi(M'\cap P_1, s\si, B(P_1, s.P, \si)\s \eta ),
 \l (m)v_{1,1_G}\rangle  .$$
 En tenant compte de (\ref{xis=xist}), de (\ref{s=}) et de l'inversibilitŽ des matrices
$B$ (cf.  DŽfinition \ref{reg}), on conclut
 que:
$$\eta= \eta', $$ ce qui Žquivaut ˆ (\ref{etas}).  \ste Il ne reste plus qu'ˆ montrer (\ref{bm}) pour achever de prouver (ii). Il s'agit simplement de transport de structure. 
D'abord, pour $\eta\in  Wh(P, \si)$, on a: 
$$\xi (P_1, s\si,  B(P_1, s.P, s\si ) \s \eta )= \xi(s.P, s\si,
\s \eta)\circ A(s.P, P_1, s\si). $$ En utilisant successivement la dŽfinition de $\s \eta$ (cf. (\ref{seta})),
  (\ref{xA}), puis  la dŽfinition des matrices $B$, on
obtient: 
$$\xi (P_1, s\si,  B(P_1, s.P, s\si ) \s \eta )= \xi(P, \si, \eta) \l(
s^{-1}) \circ A(s.P, P_1, s\si).$$
$$\xi (P_1, s\si,  B(P_1, s.P, s\si ) \s\eta )= \xi(P, \si,  \eta) \circ
A(P,   t.P_1, \si )\l(s^{-1}).$$
$$\xi (P_1, s\si,  B(P_1, s.P, s\si ) \s \eta )= \xi (t.P_1, \si,  B(t.P_1, P,
\si )\eta)  \l(s^{-1}). $$
Utilisant (\ref{xiPP1}) et   tenant compte de la relation $st=m$, on en dŽduit: 
$$\xi (P_1, s\si,  B(P_1, s.P, s\si )\s \eta )=  \xi (P_1, t^{-1}\si, \t B(t.P_1, P, \si) \eta) \l(m^{-1}). $$
On applique alors  (\ref{ximsi}) pour en dŽduire: 
 $$\xi (P_1, s\si,  B(P_1, s.P, s\si )\s \eta )=  \xi (P_1, mt^{-1}\si,  t^{-1} \si' (m^{-1})\t   B(t.P_1, P, \si) \eta).$$
 Mais $m=st$  et $t^{-1}\si(m^{-1}) = s \si (m^{-1})$.  Finalement:
 $$\xi (P_1, s\si,  B(P_1, s.P, s\si ) \eta )= \xi (P_1, s\si,   s \si'(m^{-1})  \t B(t.P_1, P, \si) \eta).$$
 Ceci prouve (\ref{bm}) et achve de prouver (ii).\\  
  (iii) rŽsulte immŽdiatement de (\ref{whps}) et (\ref{wm'gm}). \qed
  Soit $P=MU$ un sous-groupe parabolique semi-standard de $G$, $(\si, E)$ une reprŽsentation lisse de $G$. On dŽfinit,  ˆ l'aide de la DŽfinition \ref{eis} et par bilinirŽaritŽ,  une application linŽaire de $ Wh(P, \si)\otimes i_P^G E$ dans $\cu$, notŽe encore $E^G_P$ en posant:
$$E^G_P(\eta  \otimes v):= E^G_P(\si, \eta, v), \> v \in i_P^G E, \eta \in Wh(P, \si).$$
On remarquera que $\si$ a ŽtŽ omis dans  la notation. \\
C'est un entrelacement entre, d'une part, le produit tensoriel  
de la reprŽsentation triviale de $G$ sur $Wh(P, \si)$ avec $i_P^GÊ\si$ et,  d'autre part,  la reprŽsentation rŽgulire droite, $\rho$,  de $G$ sur $\cu$. \\
Soit $P'=M'U'$ un sous-groupe parabolique semi-standard de $G$, $P=MU\subset M'$  un sous-groupe parabolique semi-standard de $M'$ et soit  $Q= PU'$. Soit $(\si, E)$ une reprŽsentation lisse irrŽductible  de $M$.  De l'application: 
$$E^{M'}_P:   Wh (P, \si)\otimes i^{M'}_P  E \to  C^{\infty}(U_0\cap M' \backslash M', \psi),$$ 
se dŽduit, par fonctorialitŽ de l'induction et l'identification de $i^G_{P'}(i^{M'}_PE)$ avec  $i^G_Q E $, une application: 
 $$E_P^{P'}:Wh (P, \si)\otimes  i^G_Q E   \to  i^G_{P'} C^{\infty}(U_0\cap M' \backslash M', \psi).$$
 Si $\phi $ est un ŽlŽment de $  Wh (P, \si)   \otimes  i^G_Q E$, on le regarde comme
fonction sur $G$ ˆ valeurs dans $ Wh(P, \si) \otimes E$.  On a un isomorphisme de $G$-modules entre $i^G_{Q}E$ et 
$i^G_{P'}(i^{M'}_{P}\si ) $. Tenant compte de l'identitŽ:
$$ (\delta_Q)_{\vert M} = (\delta_P)_{\vert M} (\delta_{P'})_{\vert M},$$  
on voit facilement que cet isomorphisme, $v\mapsto {\tilde v}$, $v\in i^G_{Q}E$,  est donnŽ par: 
$$\tilde{v} (g) (m')= \delta_{P'}^{-1/2} (m') v(m'g).$$
L'Žvaluation en l'ŽlŽment neutre dans la deuxime rŽalisation de $i^G_Q E  $,  donne lieu ˆ
une application, notŽe $r_{M'}$,  de $i^G_Q E  $ dans $ i^{M'}_PE$. Soit $v \in i^G_Q E $.    On a:
$$(r_{M' }(v)) (m') =  \delta_{P'}^{-1/2}(m')  v(m'),  m' \in M'$$ de sorte que:
 \ber\label{rM'} $$ \rho(m')r_{M' }(v) = \delta_{P'}^{-1/2}(m')r_{M' }(\rho(m')v).$$ \eer On note encore
$r_{M'}$ l'application de $ Wh(P, \si) \otimes  i^G_Q ) $ dans $Wh(P, \si) \otimes  i^{M'}_{M'\cap P}E $ obtenue par tensorisation de   l'identitŽ de $Wh(P, \si)$ avec  $r_{M'}$.  On a:
\beq \label{epindep}[E^{P'}_P(\phi)(1)](m')= [E^{M'}_P(r_{M'}\phi)](m'),
m'\in M'\eeq
De plus, si $P=M'$,  $i^{M'}_{M'}E$ s'identifie naturellement ˆ $E$. Avec cette identification on a:
\beq  \label{epindepbis} (E^{P'}_{M'} \phi )(g)= E^{M'}_{M'} (\phi (g))\eeq 
On dŽfinit, pour $f\in \cu$ et pour  $P=MU$ sous-groupe parabolique standard de $G$: 
\beq \label{fpind} (f_P^{ind} (g))(m)= (\rho (g) f)_P(m), m \in M\eeq 
On vŽrifie aisŽment gr\^ace (\ref{covct}) que  $f_P^{ind} \in i_P^G\cmu$ et que l'application $f\mapsto   f_P^{ind}$  entrelace les  reprŽsentations rŽgulires droites, $\rho$,  de $G$ sur $\cu$ et $i_P^G\cmu$.\\ 

 \begin{prop} \label{cteis}  
 Soit $P=MU$(resp. $P'=M'U'$)   un sous-groupe parabolique anti-standard (resp. standard) de $G$, $\si$ une reprŽsentation lisse irrŽductible cuspidale $G$-rŽgulire de $M$.  \\
 Si  $s\in W(M'\vert G\vert M)$, on note $  C(s, P', P, \si) $ l'application linŽaire de $Wh(P, \si) \otimes i^G_PE  $  dans
$ Wh({\tilde P}_s, s \si)\otimes i^G_{P_s} sE  $ dŽfinie  par: $$ C(s, P', P, \si) =  B({\tilde P}_s, s.P, s\si)\otimes (A(P_s, s.P,  s\si) \l(s)) $$
avec l'identification de $i^G_{P_s} s E$ avec
$i^G_{P'}(i^{M'}_{M'\cap s.P}s E )$ et celle de $Wh ({\tilde P}_s, s\si)$ avec $Wh(M'\cap s.P, s\si) $ (cf. ThŽorme \ref{thct},  (i)).\\ 
Alors, pour $ \phi \in  Wh(P, \si)  \otimes i^G_P E$,
$E^G_P(\phi)_{P'}^{ind} =0  $ si $W(M'\vert G\vert M)$ est vide et sinon: 

$$E^G_P(\phi)_{P'}^{ind} = \sum_{s \in W(M' \vert G\vert
M)}E^{P'}_{M'\cap s.P} (C(s, P', P, \si) \phi), $$
 $$E^G_P(\phi)_{P'}= \sum_{s \in W(M' \vert G\vert
M)}E^{M'}_{M'\cap s.P}(r_{M'} (C(s, P', P, \si) \phi)), \> \phi \in  Wh (P, \si)   \otimes  i^G_PE.$$
\end{prop}
 \dem 
 Les deux membres de la premire ŽgalitŽ sont des fonctions sur $G\times M'$. Par Žquivariance, on se rŽduit ˆ dŽmontrer l'ŽgalitŽ en $(1, m')$ pour tout $\phi$ et tout $m' \in M'$. Cette ŽgalitŽ se rŽduit, gr\^ace ˆ (\ref{epindep}),   ˆ la seconde.  Gr\^ ace  ˆ  (\ref{covct} ) et (\ref{rM'}), on se rŽduit ˆ prouver la seconde ŽgalitŽ ŽvaluŽe en $1$.  
 Mais avec les notations du ThŽorme prŽcŽdent:
  $$E^G_{PÔ } (\phi)_{PÔ}(1) =\sum_{s \in W(M'\vert G \vert M)}\langle  \xi_s, v_s\rangle  .$$ En utilisant la dŽfinition de $\xi_s$ et $v_s$, ce ThŽorme montre la deuxime ŽgalitŽ ŽvaluŽe en $1$.   \qed 
\section{TransformŽe de Fourier-Whittaker  et produit scalaire de paquets d'ondes    } 
\setcounter{equation}{0}
\subsection{TransformŽe unipotente  de paquets d'ondes} \begin{defi}\label{tresreg}
Soit $P$ un sous-groupe parabolique anti-standard   de $G$ et $\O$ l'orbite inertielle d'une reprŽsentation lisse  cuspidale irrŽductible de $
M$.   On utilise les notations de la Proposition  \ref{fpsi} (iii). 
 \ste
On dit que $\phi \in Pol (\O_u, Wh  \otimes i^G_{P} )$ est  rŽgulire si pour tout sous-groupe parabolique standard $P'=M'U'$  de   $G$ et 
$s\in W (M'\vert G\vert M)$, l'application $\si\mapsto C(s, P', P, \si)\phi(\si)$ est polynomiale sur $\O_u$. On dit que $\phi$ est trs rŽgulire si de plus 
lorsque $M'$ et $M$ sont conjuguŽs, pour tout $g \in G$, l'application 
$\si \mapsto [(\tilde{C} (P', P, s,s\si) \phi) (s\si)](g)$
est polynomiale sur $\O_u$, o 
 \beq\label{tildeC} [\tilde{C} (s, P', P, s, s\si) \phi] (s\si): = j (P', P'^-, s.P,
s\si )     [B(P', s.P,
s\si)\otimes   (A(P', s.P,
s\si)\l(s) ] \phi(\si).\eeq  
\end{defi}
\ber Si $\phi$ est rŽgulire (resp. trs rŽgulire)  et $g\in G$, $\rho_\bullet (g)\phi$ est rŽgulire (resp. trs rŽgulire).\eer
D'aprs les propriŽtŽs de rationalitŽ des intŽgrales d'entrelacement (cf. (\ref{ratintertw}))  
et des matrices $B$ (cf. Proposition \ref{B}), on voit que: \ber \label{phitres} Il
existe une fonction polynomiale sur $\O$ non identiquement nulle, $p_\O$, 
telle que, pour tout  $\phi\in 
 Pol (\O,Wh \otimes  i^G_{P} 
 )$,   $p_\O\phi$ soit trs rŽgulire. \eer 
Par ailleurs: 
\ber\label{ptres} Si $\phi$ est rŽgulire (resp. trs rŽgulire) et $ p\in Pol(\O)$,
$p\phi$ est rŽgulire (resp. trs rŽgulire).\eer
 On note $N(D, \O)$  l'ensemble des $g\in G$ qui  normalisent $M$ et prŽservent $\O$. On note 
$W(G, \O)$ le quotient de $N(D, \O)$ par $M$ qui est fini.
Si $ p\in Pol(\O)$,   $w \in N(D,\O)$ et  $(\si, E)$ est un objet de $\O$,
$p(w^{-1}\si)$ ne    dŽpend que de la classe de $w$ dans $W(G, \O)$. Ceci permet de faire agir $W(G,\O)$ sur $Pol(\O)$ en posant $p^w(\si)=
p(w^{-1}\si)$. On remarque que si
$p\in Pol(\O)$, on a:
\beq\label{prod}\prod_{w\in W(G, \O) } p^w \in Pol(\O)^{W(G, \O)}\eeq
On dŽduit alors de (\ref{ptres}) et  (\ref{prod}) que:  
\berÊ\label{pow}
On  peut choisir  et on choisira $p_\O$ invariante par $W(G, \O)$ dans
(\ref{phitres}). 
\eer 
Montrons:  
\ber \label{phidesi} Soit $(\si_0, E_0) $ objet de  $\O_u $  et $ \phi_0  \in Wh(\si)\otimes  i^G_PE_0$
tel que
$p_\O (\si_0)$ soit non nul. Il existe
$\phi \in  Pol (\O_u, Wh  \otimes i^G_{P} )$ trs rŽgulire  telle que $\phi(\si_0)=\phi_0$ et
$\phi(w\si_0)=0$ pour  $w \in W(G, \O)$ tel que  
$w\si_0$ n'est pas Žquivalente ˆ $\si$. \eer
En effet, il existe  $\phi_1 \in Pol (\O,Wh \otimes  i^G_{P} 
 ) $   tel que $\phi_1(\si_0)=  p_\O(\si_0)^{-1} \phi_0( \si_0)$. En  multipliant $\phi_1$  par un ŽlŽment de
$Pol(\O)$ valant 1 en $\si_0$ et
$0$ en
$w\si_0$ pour  $w \in W(G, \O)$ tel que  
$w\si_0$ n'est pas Žquivalente ˆ $\si_0$, on peut supposer en outre que
$\phi_1(w\si_0)=0$ pour  $w \in W(G, \O)$ tel que  
$w\si_0$ n'est pas Žquivalente ˆ $\si_0$. En utilisant (\ref{phitres}) et
(\ref{pow}), on conclut que $\phi=p_\O\phi_1$ a les propriŽtŽs voulues.\\
D'autre part on sait (cf [DeliB]):
\ber \label{pz} Si $p\in   Pol(\O)^{W(G, \O)}$, il existe $z$ ŽlŽment
du centre de Bernstein de $G$, $ZB(G)$, tel que pour tout $(\si, E)$ objet
de $\O$,
$\igps (z)$ est la multiplication par le scalaire $p(\si)$.  \eer  
De (\ref{phitres}), (\ref{prod}) et (\ref{pz}), on dŽduit:
\begin{lem}  \label{zphi}
Si $\phi\in Pol (\O_u,  Wh  \otimes i^G_{P})$, il existe $z\in ZB(G)$ tel que
$\rho_{\bullet} (z) \phi$ soit trs rŽgulire et non nulle si $\phi$ est
non nulle.
\end{lem}

 \begin{lem} \label{suppcomp}
 Soit $\O$ l'orbite inertielle d'une reprŽsentation lisse,  cuspidale et  irrŽductible de $M$. On rappelle qu'on a choisi une mesure de Haar sur $X(M)_u$ (cf. section 1.2). On munit $\O_u$ d'une mesure $X(M)_u$-invariante et telle que pour tout objet de $\O_u$, $(\si, E)$, l'application de $X(M)_u $ dans $\O_u$, qui ˆ $\chi $ associe $[\si_\chi] $, prŽserve localement les mesures. Si  $\phi\in Pol (\O, Wh \otimes i^G_{P} )$,  on dŽfinit $f_\phi \in \cu$ par:
  $$f_\phi(g):= \int_{\O_u} E^G_P(\phi(\si))(g) d\si, \> g \in G$$
 qu'on appellera paquet d'ondes de $\phi$.\ste  
 (i) On a, pour $g \in G$, $\rho(g) f_\phi= f_{\rho_\bullet (g)\phi}$. \ste 
 (ii) Si $\phi$ est  rŽgulire,  $f_\phi $ est ŽlŽment de $\ccu$.
  \end{lem}
 \dem 
 (i) rŽsulte immŽdiatement des dŽfinitions.
 \\ Prouvons (ii). Il existe une partie compacte de $G$, $\Omega$, telle que   $G=U_0 A_0 \Omega $,   car  $M_0$ est 
compact modulo $A_0$ et $G= P_0 K$. Par ailleurs $\phi$ est invariante par un sous-groupe compact
ouvert $H$. Alors  $\Omega$ est contenu dans un nombre fini de classes ˆ droite modulo
$H$, $g_iH$. Utilisant (i) pour les $g_i$, pour  dŽmontrer ce que l'on veut il suffit donc de montrer que: 
  \ber \label{resta} Pour tout $\phi$  rŽgulire, la restriction de $f_\phi $ ˆ
$A_0$  est ˆ support compact. \eer 
Mais cela Žquivaut ˆ montrer que pour un $a_0 \in A_0$, $\rho_{a_0} f_\phi$ qui est Žgal ˆ $ f_{\rho_\bullet(a_0)\phi}$ d'aprs (i),  est ˆ support compact. D'aprs (\ref{fa-}), pour $a_0\in A_0$ bien choisi, $\rho _{a_0} f_\phi$ est ˆ support dans $A_0^-$.
On est donc ramenŽ ˆ prouver:
 \ber  Pour tout $\phi$  rŽgulire, la restriction de  $f_\phi$ ˆ $A_0^-$ est ˆ support compact.\eer 
Montrons d'abord:\ber \label{omegag} Si $\Omega_G$ est un sous-ensemble compact de
$G$ et  
$\phi \in Pol (\O_u, Wh \otimes i^G_{P})$,
$f_{\phi} $ restreinte ˆ $\Omega_G A_G$ est ˆ
support compact. \eer 
 En utilisant l'invariance de $\phi$ sous un sous-groupe compact ouvert 
de $G$ et en procŽdant comme ci-dessus, on se ramne ˆ prouver l'assertion pour $\Omega_G$
rŽduit ˆ $\{1\}$. Dans ce cas,  il suffit  de  prouver que,  si
$p$ est une fonction polynomiale sur $\O$, l'application de $A_G$ dans $\C$,  $a \mapsto
\int_{\O_u}p(\si)
\si(a)d \si$ est ˆ support compact.  Mais cela rŽsulte du fait que la transformŽe de Fourier
d'une fonction polyn\^ome sur un tore est ˆ support compact. Ceci prouve (\ref{omegag}). \\ On suppose ˆ nouveau $\phi$  rŽgulire.  
Soit $Ê\varepsilon>0$ et $Q=LV$ un sous-groupe parabolique standard de $G$. On note $\Theta_Q$ l'ensemble
des $\alpha\in \Delta (P_0)$   qui sont des racines de $A_0$ dans l'algbre de Lie de
$L$. Soit
$X_Q=
\{ a
\in A_0^-\vert  \vert \alpha(a)\vert _F  <Ê\varepsilon, \alpha \in \Delta(P_0)\setminus \Theta_Q
\>Ê\> et \>\>  
\vert
\alpha(a)\vert _F \geq Ê\varepsilon, \alpha \in \Theta_Q \}$.
Alors, les $X_Q$ forment une partition de $A_0^-$. De plus chaque $X_Q$ est
contenu dans un ensemble de la forme $A_L \Omega_Q $, o $\Omega_Q$ est un sous ensemble
compact de $A_0$. 
 Maintenant, d'aprs (\ref{ffp}),  on peut choisir  $\varepsilon >0$   tel que  pour tout
 sous-groupe parabolique standard de $G$, $Q$:
 $$f_\phi(a)= \delta_Q^{1/2}(a) (f_\phi) _Q(a), \>  a\in X_Q .$$
 On note que, d'aprs la formule du
 terme constant des paquets d'ondes (cf. [D3],   Proposition 3.17)  et, gr\^ace au fait que $ \phi$ est
 rŽgulire et  au ThŽorme \ref{thct}, 
$(f_\phi)_Q$  est une somme de  paquet d'ondes pour $L$. 
Par une application de (\ref{omegag}) ˆ chacun des termes de cette
somme,  et ceci pour tout $Q$  on voit que (\ref{resta}) est vrai.
 Le lemme en rŽsulte.  \qed 
 \begin{defi}
Si $f \in \ccu$  et $P=MU $ est un sous-groupe parabolique  anti-standard de $G$,  on dŽfinit: $$f^P(m)=
\delta_P^{1/2}(m)
\int_U f(mu) du, \> m \in M .$$
 On l'appelle la transformŽe unipotente de $f$ relativement ˆ $P$.
\end{defi} 
 Le lemme suivant rŽsulte des dŽfinitions.
 \begin{lem}
 Avec les notations ci-dessus, on a, pour $f\in \ccu$, $f^P \in \cmu$.\ste  
On dŽfinit alors
$f^{P, ind}$  par: 

\beq \label{fPind}f^{P, ind} (g, m)=(\rho (g) f)^{P}(m), \> g\in G, m\in M. \eeq
Alors $f^{P, ind}\in i_P^G\cmu$.
 \end{lem}

\begin{theo} \label{fphipind} Soit $P=MU$,   $P'=M'U'$   deux sous-groupes paraboliques anti-standard de $G$, soit $\O$ l' orbite inertielle d'une  reprŽsentation, lisse,  cuspidale irrŽductible  de $M$, et $\phi\in Pol(\O_u, Wh \otimes i^G_P) $.   Pour $\phi$ trs rŽgulire on a: \\ 
(i) Si le rang semi-simple de $M'$ est infŽrieur ou Žgal ˆ celui de $M$  et si
$M'$ n'est pas conjuguŽ ˆ $M$, $f_\phi^{P', ind}$ est nul.  
 \\ (ii)  Supposons $M$ et $M'$ conjuguŽs. 
Pour tout $m'\in M'$ et $g\in G$,  on a: $$(f_\phi^{P',ind}(g)) (m') = \sum_{s\in W(M'\vert G\vert M)} 
\int _{\O _u}    E^{M'}_{M'}  [(\tilde{C} (s, P', P,s\si) \phi)(s\si)  (g)](m')d\si. $$
 \end{theo}
\dem
 Il suffit  de prouver la formule pour  $g=1$ et $m'=1$   et pour tout $\phi$, car:
 
 $${\rho(g)}f_\phi= f_{\rho_\bullet (g)  \phi },$$ 
$$f_\phi^{P'}(m')= (\rho(m') f_\phi^{P'})(1)= \delta_{P'}^{-1/2} (m') (\rho(m') f_\phi)^{P'}(1)$$
et pour $ v' $ ŽlŽment de l'espace de  $i^G_{P'} s\si$: 
$$((i^G_{P'} s\si )(m') )v' ) (1)=   \delta_{P'}^{1/2}(m') (s\si(m') )(v'(1)).$$
On suppose maintenant $m'=g=1$. Il s'agit donc  de calculer $f_\phi ^{P'}(1)$. 
 On fixe un sous-groupe compact ouvert
$H$  comme dans (\ref{iwa}) tel que 
 $\phi$ est invariante par $\rho_\bullet (H)$. 
\ber 
Il existe une constante $C_H>0$ telle que pour tout $\varphi\in
C^{\infty}(H)$: $$\int _H
\varphi(H)dh=c_H \int _{H\cap U'^-\times H\cap M' \times  H\cap U'} \varphi(u'^-m'_1 u'
)du'^- dm' du'. $$ 
\eer 
On fixe $a \in A_{M'}$ tel que $\vert \alpha (a)\vert _F <1$ pour tout $\alpha
\in \Sigma(P')$. Pour $n \in \N$, posons $ U'_n= a^{-n} (H\cap U')a^n$. 
 Comme $f_\phi$ est ˆ support compact modulo $U_0$ (cf. Lemme \ref{suppcomp}), il existe $N\in \N$,
tel que pour tout
$n\geq N$: $$f_\phi^{P'} (1)= \int_{U'_n}f_\phi (u') du' = \delta_{P'}(a)^{-n}
\int _{H\cap U' } f_\phi (a^{-n} u'a^{n})du'.  $$ 
 On a,  pour tout $n \in \N$:  
 $$\int _{H\cap U' } f_\phi (a^{-n} u'a^{n}) du'  = \int _{H\cap U'}\int _{\O_u} 
E^G_P(\phi (\si) ) (a^{-n} u' a^n ) d\si du'.  $$ 
 On remarque que: 
 $$E^G_P(\phi (\si)) (a^{-n} u' a^n ) = E^G_P((\rho_\bullet (u'a^{n}) \phi  )
(\si) ) (a^{-n}).$$ 
 On pose alors: 
 $$\phi_n= \int _{ H \cap U'} \rho_\bullet (u'a^{n}) \phi du' \in Pol(\O_u, 
 Wh \otimes  i_P^G),$$ l'intŽgrale se rŽduisant ˆ une somme finie, puisque $\phi$ est $H$-invariante et donc $\rho_\bullet (a^n)\phi$ est invariante par $a^nHa^{-n}$.
 Alors:
  \beq \label{fp'} \delta_{P'}(a^{n}) f_\phi^{P'} (1) =\int _{H\cap U' } f_\phi
(a^{-n} u'a^{n}) du'  =
\int_{\O_u} E^G_P(
\phi_n (\si)) (a^{-n}) d \si.\eeq
   Comme $a^{-n} (H \cap M') (H\cap U'^- )a^n \subset H$, on a l'ŽgalitŽ:  
  \beq \label{phin} \phi_n = c  \int _{H\cap U' \times  H\cap M' \times H\cap  U'^-} \rho_\bullet  (
u' m' u'^-  a^n ) \phi du'^- dm' du' \eeq 
 $$ =  cc_H^{-1} \int _H \rho_\bullet (ha^n) \phi dh. $$
 o $c= vol (H\cap M')^{-1} vol (H \cap U'^-) ^{-1} $. 
 Donc $\phi_n \in Pol (\O_u,   Wh \otimes i_P^G) ^H$. \\
 D'aprs les propriŽtŽs du terme constant (cf.  (\ref{ffp})),  on
peut choisir
$N$ assez grand pour que,  pour tout
$n\geq N$, on ait:
 \beq  \label{ep'}E^G_P(\phi_n (\si ) ) (a^{-n}) = \delta_{P'^-} (a)^{-n/2}  E^G_P(
\phi_n(\si) ) _{P'^-} (a^{-n}) .\eeq 
 Si $M' $ et $M$  sont comme dans (i), $W(M'\vert G\vert M)$ est vide et $ E^G_P(\phi_n(\si) )
_{P'^-}$ est nul (cf. Proposition \ref{cteis} (i)). Cela  montre (i). 
 \\ On suppose maintenant $M'$ et $ M$  conjuguŽs. 
 Pour tout $n \in \N$,  la dŽfinition de $\phi_n$  montre que  $\phi_n$   est une combinaison linŽaire finie de fonctions du type $\rho_\bullet(g) \phi$. On dŽduit  de la dŽfinition de $f_{P'^-}^{ ind}$ et de $\phi_n$ que: 
 $$E^G_P(\phi_n (\si) )_{P'^-} (a^{-n}) =  \int _{H\cap U'} [E^G_P(
\phi(\si ) )_{P'^-} ^{ind} (u'a^{n} )](a^{-n})du'. $$ 
 En utilisant successivement  (\ref{covct}),  l'ŽgalitŽ  $\delta_{P'^-}=
\delta_{P'}^{-1}$ et  la dŽfinition de
$U'_n$ pour effectuer un changement de variable, on en dŽduit:  
$$ E^G_P(\phi_n (\si) )_{P'^-} (a^{-n}) = \delta_{P'}^{1/2} (a^{-n} )   \int _{H\cap
U'}  [E^G_P(\phi(\si)) _{P'^-} ^{ind} (a^{-n} u'a^n )] (1) du'   $$
 $$ =  \delta^{1/2}_{P'} (a^{n} ) \int _{U'_n}[E^G_P(
\phi(\si)) _{P'^-} ^{ind} (u')] (1)du' . $$
 En tenant compte de (\ref{fp'}) et (\ref{ep'}), on en dŽduit:

$$f_\phi^P(1) =\int_{\O_u}  \int _{U'_n}E^G_P(
\phi(\si)) _{P'^-} ^{ind} (u')] (1)du' d\si. $$
On va utiliser la Proposition \ref{cteis} pour donner une expression de $E^G_P(
\phi(\si)) _{P'^-} ^{ind}$. Comme $M$ et $M'$ sont 
conjuguŽs, pour tout $s \in W(M'\vert G\vert M) $, $M'\cap s.P= M'$. Dans notre utilisation du ThŽorme \ref{thct},  $P'$ est  ici remplacŽ par $P'^-$ et 
$P_s$ est ici Žgal ˆ 
$P'^-$, 
${\tilde P}_s$ est ici Žgal ˆ $ P'$. On a alors, gr\^ace ˆ la Proposition \ref{cteis}:
$$f_\phi^P(1)= \sum_{s\in W(M'\vert G\vert M)}\int _{\O_u}   \int _{U'_n} 
[E^{P'}_{M'}(C(s, P'^-, P, \si) \phi(\si))(u')](1) d\si du' $$
 et on veut passer ˆ la limite sur $n$. Mais un dŽplacement de contour d'intŽgration est nŽcessaire. Si $\chi$ est un caractre non ramifiŽ de
$M$, on note $\O_u\chi$ l'ensemble des classes d'Žquivalence des reprŽsentations
$\si _ \chi$  lorsque $\si$ dŽcrit les objets de $\O_u$.  On munit $\O_u\chi$
de la mesure obtenue par transport de structure de la mesure  sur $\O_u$. Posons:   \beq\label{phisc} \phi_s(\si):= C(s, P'^-, P, \si) \phi(\si) \in Wh(P', s\si)\otimes  i^G_{P'^-}(sE). 
 \eeq  Comme $\phi
$ est trs rŽgulire, $\phi_s$ est polynomiale en $\si$.  Pour tout choix
$\Lambda_s$,
$s\in W(M'\vert G\vert M)$,  de caractres non ramifiŽs de $M$, on a, pour des
raisons  d'holomorphie, en tenant compte de (\ref{epindepbis}): 
 \beq \label{fpoun}f_\phi^P(1)=  \sum_{s\in W(M'\vert G\vert M)}\int _{\O_u
\Lambda_s}  
\int _{U'_n} [ E^{M'}_{M'}(\phi_s(\si)(u') )](1) d\si du'.\eeq
 On veut passer ˆ la limite sur $n$ dans cette expression pour un bon choix des $\LL_s$.  
 Il faut  majorer $[ E^{M'}_{M'}(\phi_s(\si) )(u') ](1) $. \\
 Soit $(\si, E)$ un objet de $\O_u$. 
 Soit $v_s \in i^G_{P'^-} (sE)_{s \Lambda_s}$,     $\eta_s\in Wh (P', s\si)$.  On remarque que  $Wh(P', s\si)$ est Žgal ˆ $Wh(s \si)$, puisque $P'$ est anti-standard.  
Alors, il rŽsulte de la dŽfinition que:   
$$[E^{M'}_{M'}(v_s\otimes
\eta_s )(g) ](1)=\langle  \eta_s, v_s    (g)\rangle  .$$ On Žcrit:   $$u'= u^{'-} (u) m'
(u') k(u'), \> u'^-(u')\in U'^-, \> m'(u')\in M', k(u') \in K.$$
 Tenant compte des propriŽtŽs de covariance de $v_s$, on voit que: 
 $$\langle  \eta_s,v_s    (u')\rangle    = (s\Lambda_s)(m'(u')) \delta_{P'}^{-1/2}(m'(u'))\langle  \eta_s, s\si (m'(u')) 
v_s (k(u')) \rangle       .$$ 
On choisit $\Lambda_s= s^{-1} \delta_{P'}^{-1/2}$. Alors on a:
$$\langle \eta_s, v_s(u)\rangle=  \delta_{P'}^{-1}(m'(u'))\langle  \eta_s, s\si (m'(u'))  v_s
(k(u')) \rangle  , $$
Tenant compte du fait que la restriction de $v_s$ ˆ $K$ ne prend qu'un nombre fini
de valeurs, on dŽduit du Lemme  \ref{borne} (ii): 
$$ \vert \langle \eta_s, v_s(u)\rangle \vert \leq C  \delta_{P'}^{-1}(m'(u))$$
 ou $C$ est une constante qui ne dŽpend que de la restriction de $v_s$ ˆ $K$ et de
$\eta_s$ mais pas de $\si\in O_u$.
\\ On en dŽduit qu'il existe $C'>0$ tel que  pour tout $s \in W(M'\vert G\vert M)$ et $\si\in \O_u
\LL_s$: 
$$\vert E^{M'}_{M'}(\phi_s(\si))(u') ](1)\vert \leq C \delta_{P'}^{-1}(m'(u')), u'\in U'$$
Comme le second membre de cette inŽgalitŽ est une fonction intŽgrable sur $\O_u
\LL_s \times U'$, on peut appliquer le ThŽorme de convergence dominŽe et 
dŽduire de  (\ref{fpoun}): 
$$f_\phi^P(1)=  \sum_{s\in W(M'\vert G\vert M)}\int _{\O_u \Lambda_s}   \int
_{U'} [ E^{M'}_{M'}(\phi_s(\si)(u')) ](1) d\si du',$$
la fonction sous le signe intŽgrale Žtant intŽgrable pour la mesure produit. On
peut appliquer le thŽorme de Fubini et commencer par calculer l'intŽgrale sur $U'$. 
 Elle fait apparaitre l'opŽrateur $A(P', P'^-, s\si)$. On trouve,  pour $\si \in \O_u \LL_s$, en tenant compte de (\ref{phisc}) et de la dŽfininition des fonctions $C$ (cf. Proposition \ref{cteis}) et gr\^ace ˆ (\ref{epindep}) et (\ref {epindepbis}): 
 $$\int_{U'} [ E^{M'}_{M'}(\phi_s(\si))(u') ](1)du'=  
 E^{M'}_{M'}[((Id_{Wh(P', s\si) } \otimes A(P', P'^-, s\si))C(s, P'^-, P, \si)
\phi(\si))(1)](1)$$  
  Tenant compte de la dŽfinition des fonctions $C$ et de (\ref{aaa}), on voit que  
 $f_\phi^P(1)$ est Žgal ˆ la somme sur ${s \in W(M'\vert G\vert M)}$ de: $$ \int _{O_u \Lambda_s}    j (P', P'^-, s.P, s\si)  E^{M'}_{M'} 
 [ ((B(P', s.P, s\si, s\chi) \otimes  (A(P', s.P,
s\si)\l(s))
\phi(\si) )(1)](1) d\si$$ 
 En utilisant la dŽfinition des fonction $\tilde{C}$  et le fait que $\phi$ est  trs rŽgulire,   on peut remplacer l'intŽgrale sur $\O_u \LL_s$ par l'intŽgrale sur $\O_u$, pour des raisons d'holomorphie. Le ThŽorme  en rŽsulte.   \qed
 
 \subsection{ TransformŽe de Fourier-Whittaker et transformŽe unipotente}

  Soit $P=MU$ un sous groupe parabolique anti-standard de $G$ et $\O$ l'orbite inertielle  d'une reprŽsentation lisse, unitaire, irrŽductible et cuspidale de $M$.   
\begin{prop}  \label{Ehatf} 
Soit $(\si, E)$ objet de $\O_u$.  
Soit $f\in \ccu$, $g \in G$.  \\ 
(i) On a: $(\rho_\bullet (g)\hat{f}) (P, \si)= (\rho(g) f )\hat{} (P, \si)$. \ste 
(ii)  On rappelle que $Wh(\si) \otimes \igps$ est muni du produit scalaire  obtenu par produit tensoriel du produit scalaire sur $Wh(\si)$ et sur $\igps$. Pour $f \in \ccu$ et $f' \in \cu$, on note  $(f,f')_G= \int_{U_0\backslash G}f(g) \overline{f' (g) }dg$. Alors: 
$$ (f, E^G_P(\phi))_G= ({\hat f} (P, \si), \phi), \phi \in Wh(\si) \otimes i^G_PE.$$
\end{prop} 
\dem 
Les deux affirmations rŽsultent immŽdiatement de la dŽfinition de ${\hat f}$. \qed
\begin{prop}\label{findfhat}  Soit $f \in \ccu$. 
 On suppose en outre que pour tout $g\in G$, $f^{P,ind}(g)$ est ŽlŽment  de $ \ccmu$. Alors: 
 $$ {\hat f }(P,\si)(g)= (f^{P,ind}(g) {\hat )} (M, \si), g \in G,$$ 
 ŽgalitŽ qu'on rŽŽcrit:
  $${\hat f }(P,\si)=(f^{P,ind} {\hat )} (M, \si).$$
  \end{prop}
 \dem 
DŽfinissons:
 Soit $v \in i^G_PE$.
$$I:= ({\hat f}(P, \si),  \eta \otimes v).$$
 Utilisant la dŽfinition de la transformŽe de Fourier puis celle des intŽgrales de Jacquet, on voit que:
$$ I=(f,  E^G_P(\eta \otimes  v))_G, $$
c'est ˆ dire: 
  $$I= \int_{U_0\backslash G} f(g) {\overline {\langle  \xi(P, \si, \eta), i^G_P\si (g)v\rangle  } }dg.$$
  Avant de poursuivre la preuve de la Proposition,  montrons: 
  \ber  \label{tauu}
  Soit $\Omega$ un sous-ensemble compact modulo l'action ˆ gauche de $U_0$.
  Il existe une fonction continue ˆ support compact sur $G$
  telle pour tout 
  $g\in \Omega$, $\int_{U_0} \tau (ug)dg = 1$.
  \eer   
  Soit $s$ une section continue de la projection, $p$, de $G$ sur  $U_0\backslash G$ (cf. e.g. [M]).  On considŽre une fonction $\tau_1$ continue a support compact sur $U_0\backslash G$ 
  et Žgale ˆ 1 sur un voisinage  de $\Omega$, regardŽ ici  comme un sous-ensemble de $U_0\backslash G$.  On note $\tau_0$ une fonction continue sur $U_0$,  ˆ support compact   et d'intŽgrale 1.
 On pose $$\tau(g) = \tau_0(u(g)) \tau_1 (p(g)), \> avec \> u(g) = g (s(p(g)))^{-1}$$
 On vŽrifie qu'elle satisfait  toutes les propriŽtŽs voulues. Ceci prouve (\ref{tauu}).  
   \\ Appliquant ceci au support de $f$ et reprenant le calcul de $I$, on trouve:
    $$I= \int_{ G} f(g) \tau (g) \overline {\langle  \xi(P, \si, \eta), i^G_P \si(g) v \rangle  }dg.$$ 
On choisit maintenant un objet de $\O$, $(\si, E)$,  tel que  $\xi(P, \si, \eta) $  soit reprŽsentŽ par une fonction continue (cf. Proposition \ref{xicont} (ii)). On en dŽduit: 
\beq \label{xii} < \xi(P, \si, \eta), i^G_P\si (g) v>= \int_{U^-} \psi(u^{-})^{-1}< \eta, v(u^-g) >
du^- , \eeq o l'intŽgrale est absolument convergente. 
Donc:  $$ I= \int_G \tau(g) f(g) \int_{U^-}\overline { \ \psi(u^{-})^{-1}< \eta, v(u^-g) > }
du^-dg.$$
L'application  $\tau f$  est une application continue sur $G$, ˆ support  compact. Ce support est donc contenu  dans un nombre  fini de $H$-classes a droite, o $H$ est un sous-groupe compact ouvert de $G$ fixant $v$. Comme pour tout $g\in G$, 
l'intŽgrale du membre de droite de (\ref{xii}) est absolument convergente,      le ThŽorme de
Fubini s'applique. En utilisant l'ŽgalitŽ $f(u^-g) =\psi(u^-) f(g) $ et le fait que $\psi$ est unitaire,  on a:  
  $$I= \int_{U^-}  \int _G \tau(g) f(u^-g)\overline{ < \eta, v(u^-g)> }dg du^-.$$
 On change $g $ en $u^-g$ dans l'intŽgrale intŽrieure:  
  $$I= \int_{U^-}  \int _G \tau((u^-)^{-1} g) f(g) \overline{ < \eta, v(g)>} dg du^-.$$
Tenant compte de l'invariance  ˆ droite par $U_0\cap M$ de $f(g) \overline{< \eta, v(g)>
}$, on obtient: 
$$I= \int_{U^-} \int_{U_0 \cap M }\int _ {U_0\cap M\backslash G} \tau((u^-)^{-1} u_0g)
f(g)\overline{ <
\eta, v(g)> }dg du_0 dh.$$
Transformant la succession des intŽgrales sur $U^-$ et $U_0\cap M$ en une intŽgrale sur
$U_0$ et en utilisant les propriŽtŽs de
$\tau$ (cf. 
(\ref{tauu})), on en dŽduit: 
 $$I=\int _{U_0\cap M \backslash G}  f(g) \overline{< \eta, v(g)>} dg $$
 Utilisant la formule intŽgrale (\ref{iwas}) et tenant compte du fait que $v$ est invariante
ˆ gauche par
$U$, il s'ensuit: 
 $$I=   \int_{U\times (U_0\cap M \backslash M) \times U^-} f(umu^-) \overline{<
\eta, v(mu^-)>}Ê\delta_P^{-1} (m) du dm du^-.$$
Mais on a:  
 $$< \eta, v(mu^-)>=  \delta_P^{1/2}(m) E^M_M(  \eta \otimes v(u^-) )(m) $$
 et $$\int_U f(umu^-)du =  \delta_P ^{1/2} (m)   [f^{P, ind} (u^-)] ( m). $$
 Donc les fonctions modules disparaissent et l'on a:
 $$I= \int_{(U_0\cap M \backslash M)[
\times U^-} f^{P, ind} (u^-)] ( m)
\overline {E^M_M(  v(u^-) \otimes \eta)(m)} dm du^-$$
Mais $(f^{P, ind}{\hat) }(M, \si)$ est un ŽlŽment de $\lgp$ 
et $$I= ((f^{P,ind} {\hat )} (M, \si), \eta \otimes v). $$
Comme cela est vrai pour tout $v \in i_P^G E$, $ \eta\in Wh(P, \si)$, cela
prouve l'ŽgalitŽ voulue, pour $\si$ comme-ci dessus. Cette ŽgalitŽ s'Žtend ˆ  tout $(\si, E)$ objet
de $\O$ par polynomialitŽ des deux membres.\qed
 \subsection{ TransformŽe de Fourier-Whittaker  de paquets d'ondes}
 \begin{theo} \label{fphihat} Soient $P=MU, P'=M'U'$  des sous-groupes paraboliques anti-standard  de $G$ tels que  $M$ et $M'$ soient conjuguŽs, $\O$ l'orbite inertielle d'une reprŽsentation lisse cuspidale irrŽductible   de $M$.  Soit $\phi \in Pol (\O_u, Wh \otimes i^G_P)$  trs rŽgulire. Soit  $(\si_1, E_1)$ une reprŽsentation lisse, cuspidale, unitaire et  irrŽductible de $M'$. Alors 
 
 $$\hat{f}_ \phi(g) (P', \si_1)= \sum_{ s \in W(M'\vert G\vert M), s^{-1} \si_1\in
{\O_u}}[\tilde{C} (P', P, s, \si_1) \phi] (\si_1).$$
 \end{theo}
  \dem 
  Traitons d'abord le cas $M=G$. Soit $(\si_1, E_1) $ une reprŽsentation lisse, cuspidale, unitaire et irrŽductible de $G$, $\phi_1 \in
Wh(\si_1)\otimes  E_1$.
 \ste Pour toute reprŽsentation $\pi$ de $G$ admettant un caractre central, notons
$\chi_\pi$ sa restriction ˆ $A_G$. On a: 
 $$(f_\phi, E^G_G (\phi_1))_G=  \int_{U_0\backslash G} f_\phi(g) \overline{ E^G_G (\phi_1)(g) }dg =
\int_{A_GU_0\backslash G}  \varphi(g) (\overline{ E^G_G 
(\phi_1))(g)} dg$$  
 o $$ \varphi(g)= \int_{A_G} f_\phi (ag) \chi_{\si_1}(a^{-1}) da .$$
 Comme pour $(\si, E)$ objet de $\O$, $E_G^G(\phi(\si))(ag)= \chi_\si(a)E_G^G(\phi(\si))(g)$, on  a
l'ŽgalitŽ:  $$
\varphi (g) =
\int _{A_G} 
\chi_{\si_1}(a^{-1})  
\int_{\O_u} 
 \chi_\si (a)(E^G_G(\phi(\si)) (g)    d\si da. $$
La restriction de $\chi_\sigma$ ˆ $A_G \cap K$ ne dŽpend pas de $(\si, E)$ objet de
$\O$.  
 Si les  restrictions de $\chi_\si$ et $\chi_{\si_1}$ ˆ $A_G\cap K$ sont distinctes,
alors $\varphi(g)=0$.  Supposons ces restrictions Žgales. 
 On peut appliquer  la formule d'inversion de Fourier sur $A_G/ A_G \cap K$. 
La dŽfinition de l'intŽgrale sur $\O_u$  et la normalisation des mesures sur  $A_G/ A_G \cap K$, $X(G)_u$ (cf. section 2.2) et $\O_u$ (cf. Lemme \ref{suppcomp}) conduit   ˆ  l'ŽgalitŽ:   
$$\varphi= \sum_{\si \in \O_{u},  {\chi_\si }_{\vert A_M}= {\chi_{\si_1}}_{\vert A_M}}
E^G_G (
\phi(\si)) $$
D'o l'on dŽduit: $$(f_\phi, E^G_G (\phi_1))_G= \sum _{\si \in \O_{u},  {\chi_\si}_{\vert A_M}=
{\chi_{\si_1}}_{\vert A_M}}  \int_{A_GU_0 \backslash G} (E^G_G (\phi(\si))(g) \overline 
{(E^G_G (\phi_1))(g)} dg$$ 
  Le Lemme de Schur montre que  le terme correspondant ˆ $\si$ dans cette  dernire
expression est nul si
$\si$ n'est pas Žquivalente ˆ  $\si_1$. De plus si $\si=\si_1$, d'aprs le Lemme
\ref{prodscal},  ce terme est Žgal ˆ $ (\phi(\si_1), \phi_{1})_G$. 
On obtient finalement:  
\ber $(f_\phi, E^G_G (\phi_1))= 0$ si $(\si_1,  E_1)$ n'est pas objet de $\O$.
\ste  $(f_\phi, E^G_G (\phi_1))= (\phi(\si_1), \phi_1)$ si $(\si_1,  E_1)$ est un 
objet de $\O$. \eer 
En utilisant le Lemme \ref{Ehatf}, on en dŽduit le ThŽorme dans le cas $M=G$.
\\ Retournons au cas gŽnŽral. On a, d'aprs le ThŽorme \ref{fphipind}:   
$$f_\phi^{P', ind}(g)= \sum_{s \in W(M'\vert G\vert M)} f_{\phi_s}^{M'}$$
 o $ \phi_s$ est la fonction sur $s\O_u$ ˆ valeurs dans $Wh \otimes i^{M'}_{M'} $ dŽfinie par: 
$\phi_s (\si_1)=  [\tilde{C} (s, P', P, \si_1) \phi(g)] (\si_1)$. C'est un ŽlŽment de $Pol (s\O_u,  Wh \otimes i^{M'}_{M'} )$, d'aprs le fait que $\phi$ est trs rŽgulire. On utilise le Lemme  \ref{suppcomp} pour  $M'$ au lieu de $G$ et $P$, pour voir que l'on peut appliquer  
  la Proposition  \ref{findfhat}  ˆ $f_\phi$  pour
exprimer $\hat{f}_\phi (P', \si_1)$ ˆ l'aide  de $f^{P', ind}$. Joint ˆ ce que q'on
vient de dŽmontrer pour le groupe $M'$, cela implique le  ThŽorme. \qed  
\subsection{Produit scalaire de paquets d'ondes}
\begin{prop} 
 \label{prodscalpaq}
Soit $P=MU, P'=MU'$ deux sous-groupes paraboliques anti-standard de $G$. Soit $\O$ (resp. $\O_1$) l'orbite inertielle d'une reprŽsentation lisse, irrŽductible et  cuspidale de $M$ (resp. $M'$), $\phi \in Pol (\O_u, \lgp)$, $\phi_1 \in Pol(\O_{1u}, Wh \otimes i^G_{PÔ})$ trs rŽgulires. \\ 
(i) Si $M$ et $M'$  ne sont pas conjuguŽs
dans $G$, $(f_\phi, f_{\phi_1})_G$ est nul.   \ste 
(ii) Si $M$ et $ M'$  sont conjuguŽs dans $G$, $(f_\phi, f_{\phi_1})_G$ est
Žgal ˆ:
 $$ \int_ {O_{1u} }\sum_{ s \in W(M'\vert G\vert M),  s{\cal O}= {\cal
O}_1 } ([\tilde{C} (s, P', P, \si_1) \phi] (\si_1), \phi_1(\si_1))d\si_1.$$
\end{prop}
 \dem La  dŽmonstration est semblable ˆ celle de [W], Proposition VI.2.2. Nous la donnons pour la commoditŽ du lecteur.\ste 
On a l'ŽgalitŽ $$(f_\phi, f_{\phi_1})_G=  \int_{U_0\backslash G}
 f_\phi(g) \int_{{\O_1}_u} \overline{(E^G_{P'}(\phi_1(\si_1))(g)}d\si_1
dg. $$
 Comme ${\O_1}_u $ est compact et  comme $f_\phi$ est ˆ support compact 
d'aprs le Lemme \ref{suppcomp}, l'intŽgrale double est absolument convergente et l'on obtient:   
 $$(f_\phi, f_{\phi_1})_G=  
  \int_{{\O_1}_u} (f_\phi, E^G_{P'}(\phi_1(\si_1)))_Gd\si_1 $$ 
et d'aprs la dŽfinition de $\hat{f}$:
$$(f_\phi, f_{\phi_1})_G=  
  \int_{{\O_1}_u}(\hat{f_\phi}(P', \si_1), \phi_1(\si_1))d\si_1 . $$
Supposons le rang semi-simple de $M'$ infŽrieur ou Žgal ˆ celui de $M$.
En utilisant la Proposition  \ref{findfhat} et le ThŽorme \ref{fphipind} (i),
on voit que  $(f_\phi, f_{\phi_1})_G=0$ si $M$ n'est pas conjuguŽ ˆ $M'$.
Si le rang semi-simple de $M'$ est strictement plus grand que celui de
$M$, il suffit d'appliquer la relation $(f_\phi, f_{\phi_1})_G=
\overline{(f_{\phi_1}, f_\phi)_G}$, pour achever la preuve de (i).
\ste Supposons maintenant $M$ et $M'$ conjuguŽs dans $G$. 
Le ThŽorme  \ref{fphihat} calcule  $\hat{f_\phi}(P', \si_1)$, ce qui
conduit ˆ (ii).\qed

\subsection{ Adjoint de la matrice $B$}

\begin{theo} \label{adjB}
Soit $P$, $Q$ des sous-groupes paraboliques  semi-standard de $G$ possŽdant le mme
sous-groupe de LŽvi semi-standard,  $M$. On suppose $P$ anti-standard. Soit $\O$ l'orbite
inertielle d'une reprŽsentation lisse, irrŽductible et cuspidale de $M$. On a l'ŽgalitŽ de
fonctions rationnelles sur $\O_u$:
$$B(P, Q, \si)^*= B(Q, P, \si)$$
\end{theo}
\dem
Soit  $P'$ le sous-groupe parabolique anti-standard de $G$ auquel $Q$
est conjuguŽ et $\O_1$ une orbite inertielle de $M'$ conjuguŽe de $\O$ par un ŽlŽment de $W(M'\vert G\vert M)$. 
Soit $\phi\in Pol (\O_u, i^G_P\otimes Wh)$, 
$\phi_1 \in Pol (O_{1u}, i^G_{P'}\otimes Wh) $ trs rŽgulires.   Alors,  d'aprs la Proposition prŽcŽdente: 
\beq \label{ff}(f_\phi, f_{\phi'})_G =     \int_ {O_{1u}}
 \sum_{ s \in
W(M'\vert G\vert M), s{\cal O}= {\cal O}_1 } ([\tilde{C} (s, P', P, \si_1) \phi] (\si_1), \phi_1(\si_1))d\si_1.
\eeq 
Puis en utilisant $(f_\phi, f_{\phi'})_G = \overline{(f_\phi', f_\phi)_G}$, on a: 
$$(f_\phi, f_{\phi'})_G=\int_ {O_u} \sum_{ t \in W(M\vert G\vert M'),  t^{-1}  \O=  \O_1} (\phi(\si), [\tilde{C} (t, P, P', \si) \phi_1] (\si)) d \si; $$
 On pose 
$\si_1= t^{-1} \si$. D'o: 
\beq \label{barff} (f_\phi, f_{\phi'})_G= \int_ {O_{1u}} 
\sum_{ t \in W(M\vert G\vert M'),  t^{-1} 
\O=
 O_1}(\phi(t\si_1), [\tilde{C} (t, P, P', t\si_1) \phi_1] (t\si_1)) d\si_1\eeq
  A $s \in  W(M'\vert G\vert M) $ correspond un unique ŽlŽment $t$ de $W(M\vert G\vert M')$ tel que $st=m'\in
M_0 \cap K$.  
\\Montrons l'ŽgalitŽ de fonctions rationnelles  sur $\O_{1u}$:
$$\label{mumu} j
(P', P'^-, s.P,
\si_1)= \overline{j (P, P^-, t.P',
s\si_1) }$$
D'abord, d'aprs (\ref{reelunit}) et (\ref{aaa}) les deux membres de l'ŽgalitŽ ˆ
prouver sont rŽels, donc on peut ignorer la conjugaison complexe. Alors l'ŽgalitŽ
rŽsulte immŽdiatement de (\ref{aaa}) et (\ref{jaawsi}).\\
Si $\phi_1$ est trs rŽgulire, il en va de mme de
$p_1\phi_1$ pour tout $p_1\in Pol(\O_{1})$. Par ailleurs si $F\in
Pol(\O_{1u})$ est tel que:  $$\int_{\O_{1u}}p_1(\si_1) F(\si_1) d\si_1=0,
p_1\in Pol(\O_1),  $$ on en dŽduit que $F=0$.\\
 Donc, pour tout $\si_1$ objet de $ \O_{1u}$, les
expressions sous le signe intŽgrale dans les membres de droite des
ŽgalitŽs (\ref{ff}) et (\ref{barff}) sont Žgales. \\ Soit $\O'_{1u}$ l'ensemble
des 
$\si_1Ê\in
\O_1$ tel que, avec les notations de (\ref{phitres}), (\ref{pow}), 
$p_{\O_1}(\si_1)$ soit non nul et   tels que si $w\in W(M', \O_1)$, $w\si_1$ ne
soit Žquivalente ˆ
$\si_1$ que si $w=1$. On remarque que $\O '_{1u}$ est dense dans $\O_{1u}$. \\Soit $\si_1 \in \O'_{1u}$, 
$ s\in W(M'\vert G\vert M)$ et $t$ comme ci-dessus. D'aprs (\ref{phidesi}) on
peut choisir
$\phi$ tel que $\phi(t\si_1)$ soit non nul et arbitraire dans
$i^G_P(tE_1)$ et tel que
$\phi(t' \si_1)= 0$ si $ t'\in W(M\vert G\vert M')$ est distinct de $t$. De mme $\phi_1(\si_1)$ peut tre choisi arbitrairement.  Alors, tous les termes sous le signe intŽgral de (\ref{ff}) (resp.(\ref{barff})) sont nuls exceptŽ celui correpondant ˆ $s$ (resp. $t$), d'aprs la dŽfinition des fonctions $\tilde{C}$ (cf. DŽfinition \ref{tresreg}). On en dŽduit: 
$$([\tilde{C} (s, P', P, \si_1) \phi] (\si_1), \phi_1(\si_1))= (\phi(t\si_1), [\tilde{C} (t, P, P', t\si_1) \phi_1] (t\si_1))$$
Puis, utilisant la dŽfinition des fonctions ${\tilde C}$ (cf. (\ref{tildeC})), on voit que: 
$$ (B(P', s.P, \si_1) \otimes  A(P, s.P,
\si_1)\l(s) )
\phi(s^{-1} \si_1), \phi_1 (\si_1))$$  est Žgal ˆ
$$(
\phi (t\si_1),  (B (P, t.P', t\si_1) \otimes A(P, t.P', t\si_1)\l(t) )
\phi_1(\si_1) )$$
Comme $\phi_1(\si_1)$ est arbitraire, on obtient par adjonction:
\beq \label{bb*} (B(P', s.P, \si_1) \otimes A(P', s.P, \si_1)\l(s)) \phi(
s^{-1}
\si_1) = [
  B(P, t.P', t\si_1)\otimes (A(P, t.P', t\si_1)\l(t) ]^* \phi (t \si_1)\eeq
  Maintenant, des ŽgalitŽs provenant du transport de structure et la formule d'adjonction pour les intŽgrales d'entrelacement vont permettre d'achever la preuve du ThŽorme. 
Soit $m=ts\in M_0Ê\cap K$, de sorte que $t= ms^{-1}$.
On a,  d'aprs la Proposition  \ref{B} (ii),  en prenant $\si$ Žgal ˆ $s^{-1}\si_1$:
$$ B(t.P', P, t\si_1)= (s^{-1}\si_1)' (m^{-1})B(t.P', P,s ^{-1}
\si_1 )(s^{-1}\si_1') (m)$$
Comme $(s^{-1} \si_1 )(m) = \si_1(m')$, o $m'= st$,  on a: \beq \label{bmm}  B(t.P', P, t\si_1)= \si'_1(m'^{-1}) B(t.P', P,s ^{-1}
\si_1) \si_1' (m'). \eeq  
Par ailleurs, d'aprs (\ref{bm}),  dans lequel on remplace $m$ par $m'$, $P_1$ par $P'$,  $s \si$ par $\si_1$ et o $\s$ et $\t$ sont triviaux car $P$ et $P'$ sont anti-standard, on a: 
\beq \label{Bs} B(P', s.P,
\si_1)=
\si'_1(m'^{-1})B(t.P', P,s ^{-1} \si_1 ).  \eeq
Comme $\phi \in Pol (\O, Wh \otimes i^G_P  ) $,  $s^{-1} = m^{-1}t$ et $t\si_1=  m s^{-1} \si_1$,  on a, d'aprs la dŽfinition de $ Pol (\O_u, Wh \otimes i^G_P  )$ (cf. section 1.4):
$$\phi(s^{-1} \si_1)=  (t\si_1'(m) \otimes   \l(m ^{-1}))\phi (t\si_1).  $$
Donc le membre de gauche, $I$,   de (\ref{bb*}) est Žgal ˆ:  
$$ [\si_1(m'^{-1})B(t.P', P,s ^{-1} \si_1 )\otimes A(P', s.P, \si_1) \l(s)] [
 (t \si_1' (m) \otimes \l(m ^{-1})) \phi (t\si_1) ].$$
Mais:  $$\l(s)\l(m^{-1})= \l(t^{-1}), t\si_1(m)= \si_1(m').$$
Donc $$I=   [\si_1'(m'^{-1})B(t.P', P,s
^{-1}
\si' ) \si'_1(m')  \otimes A(P', s.P, \si_1) \l(t^{-1})] \phi (t\si_1). $$ 
Donc, d' aprs (\ref{bmm}):
$$I=   [B(t.P', P, t\si_1) \otimes A(P', s.P, \si_1) \l(t^{-1})] \phi(t\si_1). $$
Etudions maintenant le  membre de droite, $II$,  de (\ref{bb*}).  
L'adjoint de $A(P, t.P', t\si_1)$ est Žgal ˆ $A(t.P', P, t\si_1)$, celui
de $\l(t)$ est Žgal ˆ $\l(t^{-1})$. Enfin (cf.(\ref{xA})) $\l(t^{-1}) A(t.P', P, t\si_1)$ est Žgal ˆ 
$A(P', s.P, \si_1) \l(t^{-1})$.  Donc on a:
$$II= (B(P, t.P', t\si_1)^*    \otimes  A(P', s.P, \si_1) \l(t^{-1})) 
\phi(t\si_1). $$
Comme $\phi (t\si_1)$  peut \^etre choisi arbitrairement,  l'ŽgalitŽ de $I$ et $II$  
conduit ˆ:
$$B(t.P', P, t \si_1)= (B(P, t.P', t\si_1))^*, $$
pour notre choix de $\si_1$. Par densitŽ, cette ŽgalitŽ est vraie pour
$\si_1 \in \O_{1u}$.  Soit $s \in W(M'\vert G\vert M)$ tel que $P'=s.Q$, qui existe d'aprs notre choix de $P'$. Alors  $t. P'= Q$. En posant $\si_1= t^{-1}\si$, on obtient l'ŽgalitŽ
voulue. \qed
\section{Preuve du ThŽorme de Paley-Wiener}
\setcounter{equation}{0}
\subsection{Paquets d'ondes dŽcalŽs}
\begin{prop}\label{supcomp} Soit $P=MU$ un sous-groupe parabolique  anti-standard de $G$ et  soit $\O$ l'orbite inertielle d'une reprŽsentation lisse,  irrŽductible et cuspidale de $G$. Si 
$ \phi \in Pol(\O, Wh \otimes i^G_P)$, on note:  $$\Phi(\si):=(Id \otimes A(P^-, P, \si) ^{-1}) \phi(\si) $$
et $$f_\Phi^{sh}:= \int _{\O_u \chi_\mu, Re \mu <<_P 0} E^G_P(\Phi(\si)) d\si.$$
o $Re \mu <<_P 0$ veut dire  que $-\langle Re \mu,\check{\alpha}\rangle   $ est suffisament grand, pour tout $\aa \in \Sigma(P)$.
Alors $f_\Phi^{sh}$  ne dŽpend pas de $\mu$ et est ˆ  support compact modulo $U_0$. 
\end{prop}
On se rŽduit, comme dans la preuve du Lemme \ref{suppcomp},  ˆ dŽmontrer que pour tout $\phi$,  la restriction de $f_\Phi^{sh}$ ˆ $A_0^-$  est ˆ support compact. 
Puis on procde comme dans [H], l'analyse des p™les des fonctions rencontrŽes Žtant dŽtaillŽe ci-dessous. Les p™les potentiels ici et dans [H] vŽrifient des conditions similaires, ce qui autorise les mmes dŽplacements de contour d'intŽgration.    On note que l'on travaille ici sur $A_0^-$  et que les sous-groupes paraboliques utilisŽs pour le terme constant sont ici standard tandis que $P$ est anti-standard. Modulo l'Žtude des p™les ci-dessous, la preuve vaut mutatis mutandi, en tenant compte de notre dŽfinition diffŽrente de $H_G$ (cf. (\ref{HG})). \\
Etudions les  p™les de  de $C(s, P', P, \si) (Id \otimes A(P, P^-, \si)^{-1})$ qui est Žgal ˆ  $$  (B({\tilde P}_s, s.P, s\si)\otimes (A(P_s, s.P,  s\si) \l(s)  ) (Id \otimes A(P,  P^-, \si)^{-1}).$$ On a: 
$$A(P_s, s.P,  s\si) \l(s)A(P^-, P \si)^{-1}= \l(s) A(s^{-1}. P_s, P, \si)A(P, P^-, \si)^{-1} $$
 $$ A(s^{-1}.P_s, P^-, \si) A(P^-, P, \si)=(\prod_{ \alpha \in \Sigma_{red}(P)\cap  \Sigma_{red} (s^{-1}. P_s)} j_{\alpha}(\si))A(s^{-1}. P_s, P, \si).$$ 
 Donc $$A(P_s, s.P,  s\si) A(P^-,  P, \si)^{-1}=(\prod_{ \alpha \in \Sigma_{red}(P)\cap  \Sigma_{red} (s^{-1}. P_s)} j_\alpha^{-1}(\si) )  \l (s)  A(s^{-1}.P_s, P^-, \si). $$
 Remplaant $\si$ par $ \si_ \chi$, la fonction de $\chi\in X(M)$ ainsi obtenue a des p\^oles   pour $\chi= \chi_ \l $, avec $\l$ ŽlŽment d'un nombre fini d'hyperplans de $(\a_M)'_\C $ de la forme $\langle \l, \check {\alpha}\rangle  =c$,  $\alpha  \in \Sigma_1:=  \Sigma_{red}(P)\cap  \Sigma_{red} (s^{-1}. P_s)$ (cf. (\ref{poleA})). 
\\ La fonction sur $X(M)$, $\chi \mapsto B({\tilde P}_s, s.P, s\si_\chi)$ possde  des p\^oles     pour $\chi= \chi_ \l $, avec $\l$ ŽlŽment d'un nombre fini d'hyperplans de la forme $\langle \l, \check {\alpha}\rangle  =c$,  $\alpha \in \Sigma_2:=\Sigma_{red} (s^{-1} \tilde{P}_s) \cap \Sigma_{red} (P^-)\subset  - \Sigma_1$ ((cf. Proposition \ref{B} et (\ref{poleA})).   
 Avec les notations de Heiermann [H], qui dŽfinit $P'_s:= s^{-1} P_s$, on a   $\Sigma_1= \Sigma_{red} (P) \cap \Sigma _{red} (P'_s)$. \qed 

\subsection{Un rŽsultat d'Heiermann}
Soit $H$ un sous-groupe compact ouvert contenu dans
$K$. On note $e_H$ l'ŽlŽment de l'algbre de Hecke de $G$ dŽterminŽ par
la mesure de Haar normalisŽe de $H$. On applique la Proposition 0.2 de [H]
ˆ la famille de fonctions
$\varphi_{P,\O}$, o
$P=MU$ est un sous-groupe parabolique semi-standard de $G$,  donnŽe par
$\varphi(\si, E) =i^G_P\si(e_H)$. 
Si $(\pi, E)$ est une
reprŽsentation lisse de
$G$, on note 
$(\check{\pi},
\check{V})$ sa contragrŽdiente lisse et on identifie $V\otimes \check{V}$ ˆ un sous-espace de $EndV$. \\ Soit $P=MU$ un sous-groupe parabolique anti-standard de
$G$, $\O$ l'orbite inertielle d'une  reprŽsentation lisse,  cuspidale et  irrŽductible de $M$. 
En  transformant  la somme de Heiermann  [H] Proposition (0.2) qui porte sur l'ensemble des  $w \in W^G$ tels  que $w {\cal  O} = \cal O$ en une  somme sur  $w$ tel que $w^{-1}\in W(M\vert G\vert M)$ et $w\O=\O$,   par regroupement des  termes,  et en posant $t= w^{-1}$, on en dŽduit qu'il existe  une fonction polynomiale sur $\O$,   $\zeta (P,.)$ ˆ valeurs dans
$Hom(i_P^G, i^G_{P^-})$ telle que:
$$(i^G_P\si) (e_H) = \sum _{t \in W(M\vert G \vert M), t {\cal  O} = \cal O} A(P,  t^{-1}.P^-, \si) \l(t^{-1}) \zeta(P, t\si)  \l(t)A(t^{-1}. P, P, \si). $$
Mais, par transport de structure (cf. (\ref{xA})): $$\l(t)A(t^{-1}.P, P, \si)= A(P, t. P, t \si) \l(t).$$
Donc, on a:
\beq \label{zeta}(i^G_P\si) (e_H) = \sum _{t \in W(M\vert G \vert M), t {\cal  O}= \O} A(P, t^{-1}.P^-, \si)
\l(t^{-1})
\zeta(P, t\si)  A(P, t . P, t \si) \l(t) . \eeq
\subsection{Fin de la preuve du Thorme \ref{PW}}
On part maintenant de $F$ qui satisfait les conditions (i) ˆ (iii) du ThŽorme \ref{PW}. Soit $H$ un sous-groupe compact ouvert de $G$ tel que $F$ soit  $H$-invariante. Soit 
$P=MU$ un sous-groupe parabolique anti-standard de $G$   et $\O$ l'orbite inertielle d'une  reprŽsentation lisse, cuspidale et irrŽductible de $M$. 
Donc $$F(P, \si)= [Id \otimes (i^G_P\si) (e_H) ] F(P, \si)$$
On applique (\ref{zeta}) et on trouve que $F(P, \si)$ est Žgal ˆ: 
$$  \sum _{t \in W(M\vert G \vert M), t {\cal  O} =\O}([Id \otimes A(P,  t^{-1}. P^-, \si) \l(t^{-1}) \zeta(P,
t\si)  ] [Id \otimes A(P, t. P, t \si) \l(t) ]F(P, \si)$$
Mais $t\in W(M\vert G\vert M)$ implique que $w_{t.P}= t^{-1}$. Alors,   d'aprs (\ref{fpw}),   on a: $$(Id \otimes \l(t) )F(P, \si) = F(t.P,
t\si) 
$$
et d'aprs (\ref{fpab}),  on voit que: $$(Id \otimes A(P, t.P, t\si)) F(t.P, t\si)=(
B(t.P, P, t\si) \otimes Id) F (P,t \si ) .  $$
Donc: \beq \label{ieh} F(P, \si) = \sum _{t \in W(M\vert G \vert M), t {\cal  O} =
\cal O}[ Bt.P, P,t\si)\otimes A(P,  t^{-1}.P^-, \si) \l(t^{-1})  ]\zeta
(P, t
\si) F(P,t \si).\eeq
Soit $s$ l'unique ŽlŽment de $W^G$ tel $m:= st \in M_0 \cap K$. On utilise (\ref{bm}), avec $P_1=P$ et  $\si$ remplacŽ par $s^{-1} \si$. Comme $P$  est anti-standard,  $\s$ et $\t$ se rŽduisent ˆ l'identitŽ et l'on a: 
 \beq  \label{bwp=} B(P, s.P, \si)= \si'(m^{-1}) B(t.P,
P, s^{-1}
\si).\eeq 
On pose $ts=m'\in M_0 \cap K$. Donc $s^{-1}= m'^{-1}t$. D'aprs la Proposition \ref{B} (ii),
on a:  
$$  B(t.P, P, s^{-1} \si)= (t\si')(m') B(t.P, P, t
\si)(t\si')(m'^{-1}).
$$
Comme$(t\si)(m')= \si(m)$, on en  dŽduit:
\beq \label{btp} B(t.P, P, s^{-1} \si)= \si'(m) B(t.P, P, t\si)
\si'(m^{-1}) . \eeq
 Gr\^ace ˆ (\ref{bwp=}) et (\ref{btp}), on dŽduit de (\ref{ieh}): 
\ber \label{fpsi=}$F(P, \si)$ est Žgal ˆ: $$\sum _{t \in W(M\vert G\vert M), t{\cal  O} = \cal
O}   [ (B(P, s.P,\si) \si'(m)) \otimes (A(P,  s.P^-,
\si)
\l(t^{-1}))]\zeta (P,t \si) F(P,t
\si). $$\eer

\begin{lem}\label{zetaF} Pour $(\si, E)$ objet de $\O$, on note:  
\beq \label{PhiO} \Phi_\O (\si)= (Id \otimes  A(P^-, P, \si) ^{-1} ) \zeta (P, \si) F(P, \si).  \eeq  (i)  Pour $(\si, E)$ objet de $\O$, on a:
 \beq \label{fhat= F} \hat{f^{sh} _{\Phi_\O}}(P, \si)= F(P, \si).  \eeq 
 (ii) Si  $(\si_1, E_1)$ est une reprŽsentation lisse,  cuspidale et  irrŽductible du sous-groupe de LŽvi,  $M'$,  d'un sous-groupe parabolique standard de $G$, $P'=M'U'$,   dont l'orbite inertielle, $\O'$,  est telle que  $(M, \O)$ n'est pas conjuguŽe ˆ $(M',\O')$, on a:  
$$  \hat{f^{sh} _{\Phi_\O}}(P', \si_1)=0.$$
(iii) Si  au contraire $(M', \O')$ est conjuguŽe ˆ $(M, \O)$ on a, pour $(\si_1,E_1)$ objet de $\O'$. 
$$ \hat{f^{sh} _{\Phi_\O}}(P', \si_1)= F(P', \si_1). $$
\end{lem}
\dem
 Pour $z\in ZB(G)$, on note $Fi'= \rho_\bullet (z) F $, et  $\Phi'_\O$ la fonction rationnelle dŽduite de $F'$, comme $  \Phi_\O$ l'est de $F$.  Alors  $\Phi'_\O= \rho_\bullet (z) \Phi_\O$ et: $$ \hat{f ^{sh}_ {\Phi'_\O} }(P_1, \si_1)= (i^G_{P'} \si_1)(z) \hat{f ^{sh}_ {\Phi'_\O} }(P_1, \si_1).$$
De (\ref{pz}) Lemme \ref{zphi} et de ce qui prŽcde, on dŽduit qu'il suffit de prouver le Lemme lorsque $\Phi$ est polynomiale et  trs rŽgulire, ce que l'on suppose dŽsormais.  Alors, on peut dŽplacer le contour d'intŽgration dans la dŽfinition du paquet d'ondes dŽcalŽ et l'on a ${\hat f}^{sh}_{\Phi}= f_\phi$.
On applique le  
 ThŽorme \ref{fphihat}.  On en dŽduit (ii). Montrons (i). Soit $\si$ est un objet de $\O$. Toujours d'aprs le ThŽorme \ref{fphihat}, on a:   \beq
\label{sumis} \hat{f}_\Phi^{sh}(P,
\si)=
\sum_{ s
\in W(M\vert G\vert M), s^{-1}
\si \in O} I_s\eeq
o \ber  \label{is}  $I_s$  est Žgal au produit de  $ j (P, P^-, s.P, \si) $ par:  $$ [ 
 B(P, s.P,
\si) \otimes A(P,
s.P,\si)\l(s) ] (A(P^-, P, s^{-1} \si) ^{-1} \otimes Id ) \zeta (s^{-1} \si)
F(P, s^{-1} \si ). $$\eer 
Soit $t \in W^G$ tel que  $ts=m'$ et $st=m$ soient ŽlŽments de $M_0 \cap K$. Donc $s^{-1}=m'^{-1}t$ et $t\si(m') = \si(m)$.  On  dŽduit du Lemme
\ref{fpo}, avec $m$ changŽ en $m'^{-1}$ et $\si$ en $t\si$
, que:  
 \beq \label{fst} F(P, s^{-1} \si)= (\l(m'^{-1}) \otimes  \si'(m)) F(P,
t\si).
\eeq 
Les relations de type (\ref{fonct})  pour les intŽgrales d'entrelacement
et la fonction $\zeta$ montrent que: 
\beq \label{ast} A(P^-, P, s^{-1} \si) =  \l(m')^{-1} A(P^-, P, t\si) 
\l(m'). \eeq 
\beq \label{zst} \zeta(s^{-1} \si)= \l(m'^{-1}) \zeta (t\si)\l(m'). \eeq 
Tenant compte de (\ref{fst}),  (\ref{ast}), (\ref{zst}), on
dŽduit de la formule (\ref{is}) pour $I_s$, aprs des simplifications
Žvidentes que:  
\ber  \label{is1}  $I_s$  est Žgal au produit de  $ j (P, P^-, s.P, \si) $ par:  $$ [  B(P, s.P, \si)
\si'(m)  \otimes (A(P,
s.P,\si)\l(s)   \l(m'^{-1}) A(P^-, P, t\si)^{-1}  ] \zeta(P, t \si) F(P, t \si)  . $$
\eer    
En tenant compte de (\ref{fpsi=}),  (\ref{sumis}) et de la relation 
prŽcŽdente, il suffit,  pour prouver (\ref{fhat= F}), de  vŽrifier: 
\beq \label{muA} j (P, P^-, s.P, \si) A(P,
s.P,\si)\l(s)   \l(m'^{-1}) A(P^-, P, t\si)^{-1}= A(P, s. P^{-}, \si)\l(t^{-1}). \eeq
On a   $sm'^{-1}= t^{-1} $,  $t^{-1}.P= s.P, t^{-1}.P^-= s.P^-$. Gr\^ace ˆ (\ref{xA})
 on a:
$$\l(t^{-1}) A(P^-, P, t\si) \l(t) = A(s P^-, s.P, \si).$$
Donc notant $A$ le premier membre de  (\ref{muA}), on a:
$$A=  j (P, P^-, s.P, \si)A(P, s.P, \si)  A(s.P^-, s.P, \si)^{-1} \l(t^{-1}).$$
D'aprs  (\ref{aaa}) on a:  $$A(s.P^-, s.P, \si)
A(s.P, s.P^-, \si)=
j(s.P^-, s.P, s.P^-, \si)Id. $$
Donc:  $$ A= j (P, P^-, s.P, \si) j(s.P^-,
s.P, s.P^-, \si)^{-1} A(P, s.P, \si) A(s.P, s.P^-, \si)\l(t^{-1}). $$
Utilisant encore la dŽfinition des fonctions $j$ (cf. (\ref{aaa})), on
trouve:
$$A= j (P, P^-, s.P, \si) j (P, s.P, s.P^-)j(s.P^-,
s.P, s.P^-)^{-1} A(P, s.P^-,\si)\l (t^{-1})$$
et  on voit que:   $$j (P, P^-, s.P, \si) j (P, s.P, s.P^-, \si)j(s.P^-,
s.P, s.P^-, \si)^{-1} =1$$ donc $A= A(P, s.P^-,\si)\l(t^{-1})$ comme dŽsirŽ, ce qui
achve de prouver (\ref{muA}). Ceci achve de prouver (i).
(iii)  rŽsulte des relations (\ref{fpw}) et (\ref{fpab}) satisfaites par $F$ et les transformŽes de Fourier.\qed
{\bf Fin de la preuve du ThŽorme \ref{PW} }\\
Soit $\cal E$ un ensemble de reprŽsentants des classes de conjugaisons de couples $(M, \O)$, o $P=MU $ est un  sous-groupe parabolique anti-standard de $G$, et $\O$ est l'orbite inertielle d'une reprŽsentation,  lisse,  irrŽductible et  cuspidale de $M$. Soit  $$f= \sum_{(M, \O) \in \cal{E}} f_{\Phi_\O}^{sh},$$
la somme ne comportant  qu'un nombre fini de termes non nuls car il n'y a qu'un nombre fini de $(M, \O) $ tel que $\Phi_\O$ soit non nulle, d'aprs la condition (ii) du ThŽorme \ref{PW}. Alors, d'aprs la Propostion \ref{supcomp} , $f \in \ccu$ et d'aprs le Lemme prŽcŽdent $f$ admet $F$ comme transformŽe de Fourier-Whittaker. Pour achever la preuve du ThŽorme \ref{PW}, il ne reste plus qu'ˆ prouver la Proposition suivante.  \qed 
\begin{prop}
 La transformŽe de Fourier-Whittaker, dŽfinie sur $\ccu$,  est injective.
\end{prop}
\dem
Rappelons le contenu du  Corollaire 1 de la Proposition \ref{bernstein} de  la section \ref{app} , dont on retient les notations notamment la dŽfinition (\ref{pifxi}):  \\
Soit $f\in \ccu$. Si   pour toute reprŽsentation    lisse unitaire irrŽductible, $\pi$,  et tout $\xi \in Wh(\pi)$, $\pi' (f^{*})$ est nulle, alors $f $ est nulle.\\
Soit $P=MU$ un sous-groupe parabolique anti-standard de $G$ et  $(\si, E)$ une  reprŽsentation lisse, unitaire,  cuspidale et irrŽductible  de $M$. On remarque que, d'aprs la dŽfinition de la transformŽe de Fourier--Whittaker:
$$(\hat{f}(P, \si), \eta \otimes v) = \overline{\langle  \pi'(f^*) \xi(P, \si, \eta ), v\rangle  },  \eta \in Wh(\si), v \in i^G_PE$$
Si la transformŽe de Fourier-Whittaker de $f \in \ccu$ est nulle, on en dŽduit que $i^G_P \si (f^*) \xi= 0$ pour tout ŽlŽment de $Wh(i^G_P \si)$. 
Par polynomialitŽ, cette identitŽ s'Žtend ˆ $ \si$ reprŽsentation lisse  cuspidale et irrŽductible. \\Mais toute reprŽsentation lisse irrŽductible de $G$,  $(\pi, V)$,  apparait comme une sous-reprŽsentation d'une reprŽsentation $i^G_P \si$ avec $\si$ lisse, cuspidale irrŽductible. Par ailleurs,  l'exactitude du foncteur qui ˆ $\pi$ associe $Wh(\pi)$ montre que tout ŽlŽment de $Wh(\pi)$ est la restriction ˆ $V$ d'un ŽlŽment de $Wh(i^G_P \si)$. On dŽduit de ce qui prŽcde que $f$ est nulle, comme dŽsirŽ. Ceci achve la preuve de la Proposition et Žgalement du ThŽorme \ref{PW}. \qed

\begin{prop} Soit $\cal E$ un ensemble de reprŽsentants des classes de conjugaisons de couples $(M, \O)$, o $P=MU $ est un  sous-groupe parabolique anti-standard de $G$, et $\O$ est l'orbite inertielle d'une reprŽsentation,  lisse,  irrŽductible et  cuspidale de $M$.
Alors pour  $f \in \ccu$, notant $F$ sa transformŽe de Fourier-Whittaker et adoptant les notations du Lemme \ref{zetaF}, on a:
$$f= \sum_{(M, \O) \in \cal{E}} f_{\Phi_\O}^{sh},$$
la somme ne comportant  qu'un nombre fini de termes non nuls. 
\end{prop}
\dem D'aprs la fin de la preuve du ThŽorme \ref{PW}, les deux membres sont des ŽlŽments de $\ccu$ qui ont la mme transformŽe de Fourier-Whittaker. Donc ils sont Žgaux d'aprs la Proposition prŽcŽdente. \qed
\section{\label{LapidMao} Sur une conjecture de Lapid et Mao}
\subsection{Fonctionnelles de Whittaker de carrŽ intŽgrable et critre de Casselman}
Une  fonctions mesurable, $f$,   sur $G$ telles que $f(ug)= \chi(u)f(g)$ pour  $u \in U_0$ and $g \in G$, et telle que : $$\Vert fÊ\Vert_{\lu}:=  (\int_{U_0Ê\backslash G} \vert f(g)\vert ^2 dg)^{1/2}.$$ sera dite fonction de Whittaker de carrŽ intŽgrable. L'espace des classes  modulo l'Žquivalence presque partout de fonctions de Whittaker de carrŽ intŽgrable dŽfinit un espace de Hilbert, $\lu$, 
 sur lequel $G$ agit continument  et unitairement par reprŽsentation rŽgulire droite $\rho$. On introduit de m\^eme $(\rho, L^2(A_GU_0 \backslash G, \psi)$.
 Soit $(\pi, V)$ une reprŽsentation lisse irrŽductible  de $G$  admettant un caractre central unitaire. On dit que  $\pi$ est de carrŽ intŽgrable (resp. est une sŽrie discrte) si ses coefficients lisses sont de carrŽ intŽgrable sur $A_G\backslash G$ (resp. sur $G$). On dit que $\xi\in Wh(\pi)$ est de carrŽ intŽgrable (resp. discrte) si pour tout $v\in V$, $c_{\xi, v} $ est ŽlŽment de $L^2(A_GU_0 \backslash G, \psi)$  (resp. $\lu$).  Montrons \ber \label{auvu} Une reprŽsentation lisse irrŽductible $(\pi,V)$ de $G$ possŽde une forme linŽaire non nulle $\xi\in Wh(\pi)$ discrte (resp. de carrŽ intŽgrable)  si et seulement $(\pi, V)$ apparait comme sous-reprŽsentation irrŽductible de $(\rho, \lu)$ (resp.$(\rho, L^2(A_GU_0 \backslash G, \psi))$ \eer Traitons le cas des formes discrtes, celui des formes de carrŽ intŽgrable Žtant semblable.
Si $\xi$ est discrte et non nulle, on dŽfinit un produit scalaire  invariant sur $V$ par: 
 \beq  \label{prods}(v, v'):= \int _{ U_0 \backslash G} c_{\xi, v}(g)  \overline{ c_{\xi, v'} (g)}dg, v, v'\in V. \eeq 
Donc la reprŽsentation est unitaire et l'application $v\mapsto c_{\xi, v} $ est un entrelacement isomŽtrique de $V$ dans $\lu$. RŽciproquement, si $(\pi, V)$ est une sous-reprŽsentation de $(\rho, \lu)$, la mesure de Dirac en $1_G$ est un ŽlŽment non nul de $Wh(\pi)$. \\
 Soit $(\pi, V)$ une reprŽsentation admissible de $G$. Pour $\chi\in Hom(A_G, \C^*)$ on pose $$V_\chi
 = \{v \in V\vert Il \> existe \> d\in \N\> tel \> que \> (\pi(a)-Ê\chi(a))^d v= 0, a \in A_G\}. $$
 Si $\xi \in Wh( \pi)$, on note $\xi_\chi$ la restriction de $\xi$ ˆ $V_\chi$. On appelle exposant de $\pi$ (resp. $\xi$)  un caractre $\chi$ tel que $V_\chi$ (resp. $\xi_\chi$) soit non nul. On note ${\cal E}xp  (\pi)$ (resp. ${\cal E}xp (\xi)$) l'ensemble des exposants de $\pi$ (resp. $\xi$).  On a 
 $  {\cal E}xp  (\xi)\subset {\cal E}xp (\pi) $.
 \begin{prop} \label{cass}
Soit $(\pi, V)$ une reprŽsentation admissible de $G$ et $\xi\in Wh(\pi)$. Les conditions suivantes sont Žquivalentes:
\\(i) $\xi$ est de carrŽ intŽgrable.
\\(ii) pour tout sous-groupe parabolique standard $P=MU$ de $G$ et tout $\chi \in {\cal E}xp (\xi_P)$, on a $Re \chi \in $$^- \a_P^{G*}$,  o  $  ^- \a_P^{G*}$ est l'ensemble des $\chi \in \a_P^{G^*}$ qui sont combinaisons linŽaires ˆ coefficients strictement nŽgatifs  des racines simples de $A_0 $ dans $P$.\\ (iii) pour tout sous-groupe parabolique standard maximal  $P=MU$ de $G$ et tout $\chi \in {\cal E}xp (\pi_P)$, on a $\chi \in $$^- \a_P^{G*}$. 
\end{prop} 
\dem Soit $\Lambda$ un rŽseau contenu dans $A_0$ et tel que $A_0=( A_0\cap K ) \Lambda$. Le noyau de l'application $H_{M_0}$ (cf. (\ref{HG})) est Žgal ˆ $M_0\cap K$ et l'image de $\Lambda $ par  $H_{M_0}$ est d'indice fini dans l'image de $H_{M_0}$. Soit $I$ un ensemble d'antŽcŽdents dans $M_0$ de reprŽsentants du quotient de  l'image de  $H_{M_0}$ par l'image de $\Lambda$. 
 De l'ŽgalitŽ $G=U_0M_0K$  on dŽduit :
\beq \label{iw} G= U_0\Lambda I K. \eeq 
Soit $H$ un sous-groupe compact ouvert distinguŽ  de $K$. On voit facilement, en utilisant (\ref{iwas}),  qu'il existe des constantes $C', C''>0$ telles que:$$ C' \delta_{P_0}(\lambda^{-1}) \leq vol( U_0\backslash U_0 \lambda i xH) \leq C''\delta_{P_0}(\lambda^{-1}), \lambda \in \Lambda, i\in I, x\in K.$$ On en dŽduit que $\xi$ est de carrŽ intŽgrale si seulement si,   pour tout  $v\in V$, la restriction $c_v$ de $c_{\xi, v}$ ˆ $\Lambda$  est de carrŽ intŽgrable modulo $\Lambda\cap A_G$,  pour la mesure qui charge chaque point $\lambda \in \Lambda\cap  A_G $ de la masse $ \delta_{P_0}(\l)^{-1}$.   Par translation, $c_v$ possde cette propriŽtŽ si et seulement si c'est vrai pour $c_{\pi(a)v}$ pour un ŽlŽment $a$ de $\Lambda$. D'aprs (\ref{fa-}) , on peut donc se limiter aux $c_v$  qui sont ˆ support dans $A_0^-\cap \Lambda$.  Alors on procde comme dans la preuve du critre analogue pour les groupes [C], ThŽorme 4.4.6 en utilisant les propriŽtŽs du terme constant (cf. section \ref{termeconstant}). 
 \qed
\begin{prop} 
Soit $(\pi, V)$ une reprŽsentation admissible de carrŽ intŽgrable  de $G$ et $\xi\in Wh(\pi)$, alors $\xi$ est de carrŽ intŽgrable.
\end{prop}
\dem En effet pour tout 
 sous-groupe parabolique standard $P$, les exposants de $\xi_P$ sont des exposants de $\pi_P$ (voir ci-dessus). Alors le corollaire rŽsulte de la Proposition prŽcŽdente jointe au ThŽorme 4.4.6 de [C] (resp. ˆ la Proposition III.3.2 de [W]). \qed
 \subsection{L'analogue $p$-adique d'un rŽsultat de Wallach}
 La preuve du ThŽorme suivant est  analogue ˆ celle de son analogue rŽel donnŽe par Wallach (cf. [Wal], ThŽorme 14.12.1).
 \begin{theo}\label{supp}
 Soit $( \pi,H)$ une reprŽsentation unitaire irrŽductible de $G$ appartenant au support de la dŽcomposition en reprŽsentations irrŽductibles  de $G$ dans $\lu$. Alors la reprŽsentation lisse de $G$ dans l'espace $V$ des vecteurs de $H$ fixŽs par un sous-groupe compact ouvert est tempŽrŽe.
\end{theo} 
Raisonnant comme dans le dŽbut de la preuve de [Wall], ThŽorme 14.11. 4, on voit qu'il suffit de prouver le Lemme suivant.
\begin{lem}
Soit $f\in \lu $invariante ˆ droite par un sous-groupe compact ouvert $H$ de $K$. On note $vol(H)$ la mesure de $H$ pour la mesure de Haar sur $K$ de masse totale 1. Alors on a: 
$$\vert (\rho(g) f, f)\vert \leq vol(H)^{-1} \Vert f \Vert ^2_{\lu} \Xi(g), g \in G. $$ 
\end{lem}
\dem
On voit facilement, gr\^ace ˆ l'invariance de $f$ sous $H$ que:
$$\vert f(g)\vert^2 \leq vol(H)^{-1} \int_{K} \vert f(gk) \vert ^2  dk.$$
On pose $f_1(g)= sup_{k \in K}  \vert f(gk) \vert $. On a donc:

$$\vert f_1(g)\vert^2 \leq vol(H)^{-1} \int_{K} \vert f(gk) \vert ^2  dk.$$ 
Par intŽgration sur $U_0\backslash G$, on en dŽduit:
$${\Vert f_1 \Vert^2}_{\lu }\leq vol(H)^{-1} {\Vert f\Vert ^2}_{\lu}.$$
On voit aussi facilement que:  $$\vert (\rho(g) f,f)\vert \leq \vert (\rho(g)f_1, f_1)\vert.$$ 
On peut donc  se rŽduire ˆ prouver  l'inŽgalitŽ du Lemme en supposant $H=K$. Dans ce cas on procde comme dans la fin de la preuve du Lemme 15.1.1 de [Wall]. On doit cependant changer les intŽgrales sur $A$ en des intŽgrales sur $M_0$, changer $a$ en $m\in M_0$ et $a^{-\rho}$ en $\delta_P^{-1/2}(m)$.\qed
  \subsection{ Fonctionnelle de Whittaker de carrŽ intŽgrable  sur une reprŽsentation irrŽductible}
 Au vu de (\ref{auvu}), le  ThŽorme suivant rŽsout positivement une conjecture de Lapid et Mao (cf. [LM], conjecture 3.5). Nadir Matringe (cf. [Ma], Corollaire 3.1),  a obtenu indŽpendamment ce rŽsultat pour certains groupes. 
\begin{theo}\label{LM}
Soit $(\pi, V)$ une reprŽsentation  lisse irrŽductible de $G$ S'il existe $\xi \in Wh(\pi)$  de carrŽ intŽgrable (resp. discrte)   non nulle, alors $\pi$ est de carrŽ intŽgrable (resp. une sŽrie discrte) de $G$.
\end{theo}
\dem
Supposons $\xi$ de carrŽ intŽgrable et non nulle. Alors $(\pi,V)$ est unitaire d'aprs (\ref{auvu}).   On note $(\pi, H)$ la reprŽsentation de $G$ obtenue par complŽtion de $V$. Montrons qu'elle est contenue dans le support de $(\rho, \lu)$ (c'est trivial si $\xi$ est discrte).\\ Notons $ G^1$ le noyau de $H_G$.   le groupe 
$G^1A_G$ est d'indice fini dans $G$. Alors pour tout $v\in V$, $c_{\xi,v}\in L^{2}(U_0 \backslash G^{1})$. Donc le support de $(\pi_1,H)$ contient un ŽlŽment  du  support  de $L^{2}(U_0 \backslash G^{1})$. Par induction, le support l'induite de $(\pi_1,H)$ de $G^1$ ˆ $G$ contient un ŽlŽment du support de $(\rho, \lu)$.  Mais  cette induite  se dŽcompose en une intŽgrale hilbertienne des reprŽsentations $(\pi\otimes \chi, H)$, o $\chi$ dŽcrit l'ensemble, $X(G)_u$,  des caractres non ramifiŽs unitaires de $G$. Donc il existe $\chi\in X(G)_u$ tel que $(\pi\otimes \chi, H)$ soit ŽlŽment  support de $(\rho, \lu)$. Mais $(\rho\otimes\chi^{-1}, \lu)$ est Žquivalente $(\rho, \lu)$, l'opŽrateur de multiplication par $\chi$ Žtant un entrelacement unitaire. Donc $(\pi, V)$ est ŽlŽment du  support de $(\rho, \lu)$. \\
 D'aprs le ThŽorme \ref{supp}, $(\pi,V)$ est tempŽrŽe. C'est donc un facteur direct d'une induite ˆ partir d'un sous-groupe parabolique anti-standard $P=MU$ de $G$ d'une reprŽsentation de carrŽ intŽgrable, $(\si, E)$, de $M$. Comme $V$ est un facteur direct, il suffit   montrer que $Wh( i_P^G  \si   )$ n'a pas d'ŽlŽment non nul de carrŽ intŽgrable sauf si $P=G$. Supposons qu'il en existe un et notons le encore $\xi$. On note encore $\pi$ la reprŽsentation $i^G_P \si$ et $V$ son espace. On suppose que  $P$ est diffŽrent de $G$.\\ D'abord, d'aprs le ThŽorme \ref{jacquetint},  $\xi$ est Žgal ˆ $\xi(P, \si, \eta)$ pour un ŽlŽment non nul, $Ê\eta$, de  $Wh( \si)$. Soit $e\in E$ et 
soit $H$ un sous-groupe compact ouvert de 
$G$ contenu dans $K$ possŽdant  une factorisation d'Iwahori par rapport ˆ $(P, P^-)$ (cf. (\ref{iwa}))  et tel que $e$ soit invariant par $H_M$.  On suppose en outre que $H$ est
assez petit, de sorte  que  $H_{U^-}$  soit contenu dans $Ker \psi$. 
On
considre application de $G$ dans $E$, 
 $v_{e, \si}^{P,H}$ dŽfinie par (\ref{vsi}).  Comme $e $ est $H_M$-invariant, 
$v_{e, \si}^{P,H}$ est invariante ˆ droite par $H$. C'est un ŽlŽment de
$i^G_{P}E$ ˆ support dans $U^-P$. 
Notons $vol(H_{U^-})= \int_{H_{U^-}} du^-$, 
o $du^-$ est la mesure de Haar sur $U^-$ choisie en Ê(\ref{fintu}). Alors (cf \ref{xiv})), on voit que: 
\ber  $$\langle  \xi, v_{e, \si}^{P,H}\rangle  = vol(H_{U^-}) \langle \eta,e\rangle  . $$\eer
On note $\chi$ le caractre central de $\si$, qui est unitaire puisque $\si$ est de carrŽ intŽgrable. 
on voit facilement que pour $a\in A_M$,   $ \pi(a) v_{e, \si}^{P,H}= \delta_P^{1/2}  \chi (a)v _{e, \si}^{P,aHa^{-1}} $. Pour $a\in A_0^-$, $aH_{U^-} a^{-1}$ est contenu dans $H_{U^-}$ puisque $P$ est anti-standard. Par application de la formule prŽcŽdente et en tenant compte de l'ŽgalitŽ $vol (aU^-a^{1})= \delta_P^ {-1 } (a)
vol (U^-)$, on trouve:$$\langle  \xi, \pi(a) v_{e, \si}^{P,H}\rangle  = \chi(a)\delta_P^{1/2}(a)vol(H_{U^-}) \langle \eta,e \rangle  , a \in A_0^-\cap A_M.$$
Notons $v=  v_{e, \si}^{P,H}$.
Utilisant les dŽfinitions on voit d'aprs  (\ref{ffp}) et (\ref{tctcoef}), que, pour $\ep>0$ assez petit:
 $$c_{\xi_P, v_P}(a)=  \chi(a) vol(H_{U^-}) \langle \eta,e \rangle   , a \in A_0^-(P,<\ep ) \cap A_M.$$
Comme $ c_{\xi_P, v_P}$ est une fonction  $A_M$-finie sur $M$, on dŽduit de l'ŽgalitŽ prŽcŽdente:
 $$c_{\xi_P, v_P}(a)=  \chi(a) vol(H_{U^-}) \langle \eta,e \rangle , a \in A_M.$$
 Choisissant $e$ tel que $ \langle \eta,e \rangle$ soit bon nul, on voit que la restriction de $\xi_P$ ˆ $V_{P, \chi}$ est non nulle. Donc $\chi \in {\cal E}xp ( \xi_P)$. De plus, $\chi$ Žtant unitaire, $Re \chi $ est nul. Joint ˆ la Proposition \ref{cass}, cela contredit le fait   $\xi$ est de carrŽ intŽgrable. Cette contradiction  achve de prouver que $P=G$ et que 
 $\pi$ est de carrŽ intŽgrable. \\Si $\xi$ est discrte, $A_G$ est trivial. Donc $(\pi,V)$ est une sŽrie discrte de $G$. \\
  \qed
\section{\label{app} Appendice: Adaptation d'un rŽsultat de Joseph Bernstein ˆ notre contexte}
Soit $(\pi, V)$ une reprŽsentation lisse de $G$,  $\xi \in Wh(\pi)$ et $ f \in \ccu$.  On note  $f^*$ la fonction sur $G$ dŽfinie par  $ f^*(g) = \overline{f(g^{-1})}$ pour $ g \in G$. On dŽfinit $\pi' (f^{*}) \xi  \in \check{V}$, par:
\beq \label{pifxi}  \langle \pi'(f^*) \xi, v\rangle  =  \int_{G/U_0} f^*(g)\langle  \pi'(g) \xi, v\rangle  dg. \eeq
La reprŽsentation rŽgulire droite de $G$ dans $\lu$ se dŽcompose en une intŽgrale hilbertienne    de reprŽsentations unitaires irrŽductibles de $G$,  $(\int_Z^\oplus\pi_z d\mu(z), \int_Z^\oplus H_z d\mu(z))$. On note, pour $z\in Z $,$(\pi_z^\infty, H_z^{\infty})$  la reprŽsentation lisse de $G$, $\pi_z^\infty$, dans l'espace $H_z^{\infty} $,  des vecteurs lisses de $H_z$,  i.e. fixŽs par un sous-groupe compact ouvert de $G$. 
\begin{prop}  \label{bernstein} Pour $\mu$-presque tout $z\in  Z$, il existe un morphisme de $G$-modules, $\beta_z$, entre $\ccu$ et  $H_z^\infty$ et il existe  $\xi_z \in Wh( \pi^\infty _z)$ tels que:
\\(i)  Pour tout $ f \in \ccu$, $f= \int_Z^\oplus\beta_z( f)  d\mu(z).$
\\(ii)  Pour $\mu$-presque tout $z\in Z$, on a $$(\beta_z ( f), v)= 
 \overline{\langle \pi'_z( f^{*} )\xi_z,  v \rangle }, f \in \ccu, v\in H^{\infty}_z.$$
\end{prop}
\dem 
Pour tout sous-groupe compact ouvert $H$ de $G$ et tout sous-ensemble de $G$, $\Omega$, compact modulo l'action ˆ gauche de $U_0$, on note $\ccu^H_\Omega$ l'espace des ŽlŽments de $\ccu$ invariants ˆ droite par $H$ et ˆ support dans $\Omega$, que l'on munit de la topologie de la convergence uniforme. On note   $\lu^{H}_\Omega$ l'adhŽrence dans $\lu$ de $\ccu^H_\Omega$. Comme $H$ est ouvert, $\Omega$ est contenu dans la rŽunion d'un nombre fini d'orbites de $H$ dans $U_0\backslash G$. Donc  cet espace est  de dimension finie. On munit $\ccu$ de la topologie limite inductive des $\ccu^H_\Omega$. 
Montrons que:  \ber \label{fine} L'injection de $\ccu$ dans $\lu$ est ''fine'' dans le sens de [B1] section 1.4.\eer 
Pour cela, d'aprs [B1], section 1.6,  Lemme 2 et ThŽorme 1.5, il suffit de prouver que pour tout $H$ et $\Omega$ comme ci-dessus l'injection, $i$,  de  $\lu^{H}_\Omega $ dans $\lu$  est de Hilbert-Schmidt, ce qui est clair puisque  cet espace est de dimension finie.  Ceci achve de prouver  de prouver (\ref{fine}). \\
 Alors (i)   rŽsulte  des propriŽtŽs des applications ''fines'' (cf. [B1], section 1.4). 
 \\ On note $\alpha_z$  la restriction  ˆ  $H_z^{\infty}$ de l'adjoint de $\beta_z$ , lorsque $\beta_z$ est dŽfini. Son image est contenue dans l'espace des  vecteurs lisses du dual hermitien de $\ccu$.  Cet espace s'identifie ˆ $\cu$  par l'application qui ˆ $\phi \in \cu$ associe la forme antilinŽaire sur $\ccu$ dŽfinie par $f\mapsto \int_{U_0\backslash G} \phi(g) \overline{f(g)} dg$. On dŽfinit alors $\xi_z (v):= (\alpha_z (v))(1_G)$. On a    $\xi_z \in Wh( \pi_z^{\infty})$ . De la dŽfinition de $\alpha_z$ et $\xi_z$,  on dŽduit que: $$
 (\beta_z(f), v)= (f, \alpha_z (v))= \int_{U_0\backslash G} f(g) \overline{\xi_z (\pi_z(g)v)}dg. $$
 Alors (ii) rŽsulte de la dŽfinition (\ref{pifxi}).\qed 
 \begin{cor}Soit $f\in \ccu$, si pour toute reprŽsentation unitaire irrŽductiblle lisse de $G$ et $\xi \in Wh(\pi)$, $ 
 \pi(f^*) \xi=0$, alors $f$ est nulle. 
 \end{cor}
 \begin{rem} Les seules choses utilisŽes ici sont que $G$ est un groupe localement compact totalement discontinu, que $U_0$ est un sous-groupe fermŽ de $G$ tel que $U_0 \backslash G$ admette une mesure invariante et que $\psi$ est un caractre lisse de $U_0$. 
 \end{rem}
 \section{RŽfŽrences}

 \noindent [A]  Arthur J., A local trace formula. Publ. Math.  Inst. Hautes \' Etudes Sci.   No. 73  (1991), 5--96.

 \noindent[B1] Bernstein J.,  On the support of the Plancherel measure,  J. Geom. Phys.  5  (1988),   663--710 (1989). 
 
 \noindent[B2]  Bernstein J., Representations of p-adic groups. Lectures given at Harvard
University, Fall 1992, Notes by K. E. Rummelhart.
 
 \noindent [B3] Bernstein J.,  Second adjointness theorem for representations of reductive p-adic groups, unpublished
manuscript. 

\noindent [BZ] Bernstein J., Zelevinsky V., Induced representations of reductive $p$-adic groups I,  Annales de l'Ecole Normale SupŽrieure 10 (1977) 441-472.

\noindent [BlD] Blanc P., Delorme P.,  Vecteurs distributions $H$-invariants de
reprŽsentations induites pour un espace symŽtrique rŽductif $p$-adique
$G/H$, Ann. Inst. Fourier, 58 (2008)  213--261.

\noindent[BoT] Borel A., Tits J., Groupes rŽductifs, Publ. Math.  Inst. Hautes \' Etudes Sci. 27 (1965) 55--150. 

\noindent [BrD] Brylinski J-L.,  Delorme P.,
Vecteurs distributions $H$-invariants pour
 les sŽries principales gŽnŽralisŽes d'espaces symŽtriques rŽductifs et prolongement mŽromorphe,
d'intŽgrales d'Eisenstein. Invent. Math.
109 (1992) 619--664. 

\noindent [Bu] Bushnell, C., Representations of reductive $p$-adic groups: localization of Hecke algebras and applications.  J. London Math. Soc. 63  (2001),  no. 2, 364--386. 
 
 \noindent [BuHen] Bushnell C.,  Henniart  G., 
Generalized Whittaker models and the Bernstein center.
Amer. J. Math. 125 (2003),  513--547.

\noindent [C] Casselman W.,  Introduction  to the theory of admissible representations of $p$-adic reductive groups,
http://www.math.ubc.ca/$\tilde{}$ cass/research.html. 

\noindent[CS]  Casselman, W., Shalika, J.,  The unramified principal series of $p$-adic groups. II. The Whittaker function. Compositio Math. 41 (1980), 207--231.

\noindent [DeliB] Deligne P.,  Le ``centre'' de Bernstein rŽdigŽ par Pierre Deligne. Travaux en Cours,  Representations of reductive groups over a local field,  1--32, Hermann, Paris, 1984.

\noindent [D1] Delorme  P., Espace des coefficients de reprŽsentations admissibles d'un groupe
rŽductif $p$-adique,  131--176, Progr. Math., 220, Birkhauser Boston, Boston, MA, 2004.

\noindent [D2] Delorme P., The Plancherel formula on reductive symmetric spaces from the point of view of the Schwartz space.  Lie theory,  135--175, Progr. Math., 230, Birkhauser Boston, Boston, MA, 2005.

\noindent[D3] Delorme P., Constant term of smooth $H_\psi$-invariant functions.  Trans. Amer. Math. Soc.  362  (2010),  933--955.

 \noindent [H] Heiermann V.,  Une formule de Plancherel pour l'algbre de Hecke d'un groupe rŽductif $p$-adique. Comment. Math. Helv. 76 (2001),  388--415.
 
 \noindent [Hu] Humphreys J. E, Linear algebraic groups, Graduate Text In Math. 21, Springer, 1981.
 
 \noindentÊ[K] Knapp A., Representation theory of semisimple groups. An overview based on examples. Reprint of the 1986 original. Princeton Landmarks in Mathematics. Princeton University Press, Princeton, NJ, 2001
 
 \noindent [LM] Lapid E.,  Mao Z.,  On the asymptotics of Whittaker functions.  Represent. Theory  13  (2009), 63--81.

\noindent[Ma]  Matringe N., Derivatives and asymptotics of Whittaker functions. Preprint.  

\noindent [M]  Michael E.,  Selected selections theorems, Amer. Math. Monthly 63 (1956) 223-283.

\noindent [R] Rodier F., Modles de Whittaker des reprŽsentations admisssibles des groupes rŽductifs p-adiques quasi-dŽployŽs, manuscript  non publiŽ.

\noindent [Sh] Shahidi F.,  On certain $L$-functions.  Amer. J. Math.  103  (1981)  297--355.

\noindent [T] Tits J., ReprŽsentations linŽaires irrŽductibles  d'un groupe rŽductif sur un corps quelconque, J.  Reine Angew.   247 (1971),  196--220.

 \noindent [W] Waldspurger J.-L.,  La formule de Plancherel pour les groupes
$p$-adiques (d'aprs Harish-Chandra), J. Inst. Math. Jussieu  2  (2003),  235--333.

 \noindent [Wall] Wallach N., Real reductive groups II. Pure and Applied Mathematics, 132-II. Academic Press, Inc., Boston, MA, 1992.
 
  \noindent [War] Warner G., Harmonic analysis on semi-simple Lie groups. I. Die Grundlehren der mathematischen Wissenschaften, Band 188. Springer, New York-Heidelberg, 1972.
 \\\\
 Patrick Delorme\\Institut de Math\'ematiques de Luminy, UMR 6206 CNRS,
Universit\'e  de la M\'editerran\'ee, 163 Avenue de Luminy,
 13288 Marseille Cedex 09, France
  \\Institut Universitaire de France
\\delorme@iml.univ-mrs.fr}\\
L'auteur a bŽnŽficiŽ du soutien du programme ANR-08-BLAN-01 pendant l'Žlaboration de ce travail. 
 \end{document}